# FUCHSIAN GROUPS, QUASICONFORMAL GROUPS, CONICAL LIMIT SETS

## Peter W. Jones

## Lesley Ward

ABSTRACT We construct examples showing that the normalized Lebesgue measure of the conical limit set of a uniformly quasiconformal group acting discontinuously on the disc may take any value between zero and one. This is in contrast to the cases of Fuchsian groups acting on the disc, conformal groups acting discontinuously on the ball in dimension three or higher, uniformly quasiconformal groups acting discontinuously on the ball in dimension three or higher, and discrete groups of biholomorphic mappings acting on the ball in several complex dimensions. In these cases the normalized Lebesgue measure is either zero or one.

## 1. Introduction

A classical result of Hopf says that the normalized one-dimensional Lebesgue measure of the conical limit set of any finitely generated Fuchsian group is either zero or one. This has been extended to infinitely generated Fuchsian groups, to conformal and to uniformly quasiconformal groups acting discontinuously on the unit ball in (real) dimension three or higher, and to discrete groups of complex hyperbolic isometries of the unit ball in several complex variables. In this article it is shown that this dichotomy does not hold for uniformly quasiconformal groups acting on the unit disc in the complex plane. For each $\lambda$ with $0 \leq \lambda \leq 1$ we explicitly construct a quasiconformal group, acting on the disc, whose conical limit set has normalized one-dimensional Lebesgue measure equal to $\lambda$. The main step is to construct a doubling measure supported on the conical limit set of a certain Fuchsian group, where this conical limit set has Lebesgue measure zero.

A group $G$ of homeomorphisms of the unit disc, in dimension two, or of the unit ball, in dimension greater than two, is said to act *discontinuously* if each point in the disc or ball has a neighbourhood $U$ such that only finitely many of the images $g(U)$, $g \in G$, intersect $U$. The *limit set* $L(G)$ of such a group is the set of accumulation points of the orbit of the origin under the action of $G$. The same limit set is obtained using any point in the disc in place of the origin. The discontinuity of $G$ implies that $L(G)$ is a subset of the unit circle or sphere. A point $x \in L(G)$ is a *conical limit point* of $G$ if there is a sequence of orbit points which converges to $x$ inside a Euclidean cone with vertex at $x$, axis perpendicular to the circle or sphere, and opening angle less than $\pi/2$, so that the sides of the cone are not tangent to the circle or sphere. The *conical limit set* $L_c(G)$ of $G$ is the set of conical limit points of $G$.

Jones was partially supported by NSF grant #DMS-92-13595. Ward's research at MSRI was supported in part by NSF grant #DMS-90-22140.

A *Fuchsian group* is a group of Möbius transformations which acts discontinuously on the unit disc $\mathbb{D}$ in the complex plane. For simplicity we assume that the group has no elliptic elements. Hopf [H] showed that the conical limit set of a finitely generated Fuchsian group has normalized Lebesgue measure zero or one. This depends on whether the Poincaré series for the group converges or diverges at the exponent 1:

$$|L_c(G)| = \begin{cases} 0, & \text{if } \sum_{g \in G}(1 - |g(0)|) < \infty; \\ 1, & \text{if } \sum_{g \in G}(1 - |g(0)|) = \infty; \end{cases} \tag{1.1}$$

the result also holds for infinitely generated Fuchsian groups.

For Fuchsian groups, and for groups of Möbius transformations acting discontinuously on the $n$-dimensional ball when $n \geq 3$, there are several dichotomies which are equivalent to the zero-one dichotomy for the normalized $(n-1)$-dimensional measure of the conical limit set. The equivalence of the following conditions was proved by Sullivan [S1, S2]:

1) $|L_c(G)| = 0$;
2) The Poincaré series converges at the exponent $n - 1$: $\sum_{g \in G}(1 - |g(0)|)^{n-1} < \infty$;
3) Green's function for the Laplace-Beltrami operator exists for the hyperbolic manifold $\mathbb{B}^n/G$; and
4) the geodesic flow on $\mathbb{B}^n/G$ is transient.

Otherwise, the conical limit set has full Lebesgue measure in the $(n-1)$-sphere; the Poincaré series diverges at the exponent $n-1$; there is no Green's function for $\mathbb{B}^n/G$ (that is, the integral of the heat kernel for $\mathbb{B}^n/G$ does not converge); and the geodesic flow on $\mathbb{B}^n/G$ is recurrent and completely ergodic. See Nicholls' book [N2] for further information. Analogous results hold in several complex variables [Kam1, Kam2, MW].

A homeomorphism $\varphi: \Omega \to \Omega'$ of complex domains is called $K$-*quasiconformal* if it is in the Sobolev class $W^{1,2}_{\text{loc}}$ and its directional derivatives satisfy

$$\max_\alpha |\partial_\alpha \varphi(z)| \leq K \min_\alpha |\partial_\alpha \varphi(z)| \tag{1.2}$$

for almost every $z \in \Omega$. Geometrically, this means that for almost every $z \in \Omega$, infinitesimal circles centred at $z$ are mapped by $\varphi$ to infinitesimal ellipses centred at $\varphi(z)$, whose eccentricities are uniformly bounded below by $1/K$. A *quasiconformal group* is a group of $K$-quasiconformal maps acting discontinuously on the disc, for some fixed $K \geq 1$.

One can define quasiconformal mappings in dimension $n \geq 3$ by generalizing the geometric definition above, using spheres and ellipsoids instead of circles and ellipses. Garnett, Gehring, and Jones [GGJ], and independently Tukia [T3], showed that results analogous to Hopf's on the measure of the conical limit set and the convergence of the Poincaré series at the exponent $n-1$ hold for uniformly quasiconformal groups when the dimension is at least three: If $G$ is a group of $K$-quasiconformal mappings acting discontinuously on the $n$-dimensional unit ball, where $n \geq 3$, then

$$|L_c(G)| = \begin{cases} 0, & \text{if } \sum_{g \in G}(1 - |g(0)|)^{n-1} < \infty; \\ 1, & \text{if } \sum_{g \in G}(1 - |g(0)|)^{n-1} = \infty. \end{cases} \tag{1.3}$$

In this article we show that, in sharp contrast to the higher dimensional case, there is no zero-one dichotomy for the normalized Lebesgue measure of the conical limit set of a quasiconformal group acting on the two-dimensional unit disc.

**Theorem 1.1.** *For each number* $\lambda \in [0, 1]$, *there is a quasiconformal group* $\Gamma$ *acting on the unit disc such that the normalized one-dimensional Lebesgue measure of the conical limit set of* $\Gamma$ *is* $\lambda$.

A *doubling measure* $\mu$ on the circle is a positive measure satisfying $\mu(\widetilde{I}) \leq c\,\mu(I)$, for some constant $c > 0$, whenever $I$ is an arc of the circle and $\widetilde{I}$ is the arc with the same centre and twice the length. Equivalently, $\mu(I) \leq c'\mu(J)$ for some uniform constant $c' > 0$ whenever $I$ and $J$ are adjacent arcs of equal length.

Theorem 1.1 is a consequence of the following result:

**Theorem 1.2.** *There exist a Fuchsian group* $G$ *and a doubling measure* $\mu$ *such that the conical limit set of* $G$ *has one-dimensional Lebesgue measure zero, and* $\mu$ *is supported on the conical limit set of* $G$.

Our result is perhaps rather surprising given that Sullivan [S2] and independently Tukia [T1] have proven that every $K$-quasiconformal group acting on the disc is conjugate by a quasiconformal mapping to a Möbius group. The corresponding statement is false in higher dimensions [T2]. On the other hand, boundary values of quasiconformal homeomorphisms of the disc (in other words, quasisymmetric maps) have much less regularity than is true in higher dimensions. If $\varphi : \mathbb{B}^n \to \mathbb{B}^n$ is $K$-quasiconformal and $n \geq 3$, the restriction of $\varphi$ to $\partial \mathbb{B}^n$ is also $K$-quasiconformal and hence lies in the Sobolev space $W^{1,n-1}(\partial \mathbb{B}^n)$. In particular $\varphi|_{\partial \mathbb{B}^n}$ takes sets of $(n-1)$-dimensional Lebesgue measure zero to sets of $(n-1)$-dimensional Lebesgue measure zero. Quasisymmetric mappings do not have this property. It is this lack of regularity for quasisymmetric mappings that we exploit to construct our examples.

Theorem 1.1 is an easy consequence of Theorem 1.2. This argument, carried out in Section 3, uses little more than the Beurling-Ahlfors theorem: A homeomorphism $f$ of the circle is the restriction of a quasiconformal mapping of the disc to itself if and only if the distributional derivative $\frac{\partial f}{\partial \theta}$ on the circle is a doubling measure.

Theorem 1.2 gives some information about the possible geometry of limit sets and of sets which support doubling measures. Roughly speaking, the support of a doubling measure must be rather evenly distributed. We give a Fuchsian group whose conical limit set is evenly enough distributed to support a doubling measure, even though it has zero Lebesgue measure. (In fact, the Hausdorff dimension of this conical limit set is strictly less than one.) We construct a doubling measure that is tailored so its support is in this small set. Such constructions were first carried out by Kahane [K], and our method follows his philosophy, though the details are necessarily more complicated.

As José Fernández has pointed out to us, the Patterson-Sullivan measure associated to the group $G$ is a natural candidate for the doubling measure of Theorem 1.2. Unfortunately, certain technical difficulties are encountered in this approach, and we have not yet understood how to overcome them.

It is interesting to note that the result of Theorem 1.2 cannot be achieved if the quasisymmetric map associated to the doubling measure $\mu$ is in the Teichmüller space of the group $G$, or if $G$ is finitely generated. The first statement follows from the fact that the existence of Green's function is a quasiconformal invariant. If $G$ is finitely generated and of the first kind, then $L_c(G)$ has full measure in the unit circle. If $G$ is finitely generated and of the second kind, then $\partial \mathbb{D} \setminus L_c(G)$ contains an arc, and so $L_c(G)$ cannot support a doubling measure.

Section 4 contains a rather full outline and summary of our proof of Theorem 1.2. The rest of the paper is organized as follows.

In Section 2 we give some definitions and notation. In Section 3, we prove Theorem 1.1 as a corollary of Theorem 1.2. In Section 5 we establish some basic properties of the Fuchsian group $G$ used in our proof of Theorem 1.2. In Section 6, we prove that a certain simple construction yields doubling measures which are supported on small sets. Sections 7 to 11 contain the construction of our doubling measure, up to the specification of some parameters. In Section 7 we fix a set of fundamental domains in the universal cover $\mathbb{D}$ of the Riemann surface $\mathbb{D}/G$. In Section 8 we define a collection of Whitney intervals in the Riemann surface, and outline some properties of its lift to the universal cover. In Section 9 we construct a grid of intervals in the circle, on which our doubling measure will be built. Section 10 contains estimates which control the geometric distortion caused by the covering map. In Section 11 we define a sequence of density functions which, by the lemma in Section 6, yields a doubling measure. In Section 12 we define some auxiliary functions and outline the rest of the proof. Sections 13 and 14 contain estimates of the expectation and second moment, with respect to the doubling measure, of the auxiliary functions. In Section 15, an argument involving a constrained random walk on the Riemann surface, together with the estimates of the previous two sections, proves that the measure is supported on the conical limit set of $G$.

It is a pleasure for the second author to acknowledge the support of MSRI.

## 2. Definitions and Notation

The conical limit set is also known as the *radial* or *non-tangential* limit set. Hedlund [He] introduced conical limit points in connection with his study of horocyclic transitive points. The conical limit set may be characterized in terms of spherical caps. A *spherical cap* on a point $x$ in the unit disc $\mathbb{D}$ is an arc of the unit circle $\partial \mathbb{D}$ of the form

$$\text{Cap}\,(x,a) = \left\{ y \in \partial \mathbb{D} \mid |y - x| \leq a\,(1 - |x|) \right\}, \tag{2.1}$$

where $1 < a < \infty$. A *non-tangential cone* on a point $y \in \partial \mathbb{D}$ is a subset of $\mathbb{D}$ of the form

$$\text{Cone}\,(y,b) = \left\{ x \in \mathbb{D} \mid |y - x| \leq b\,(1 - |x|) \right\}, \tag{2.2}$$

where $1 < b < \infty$. Clearly $y \in \text{Cap}\,(x,a)$ if and only if $x \in \text{Cone}\,(y,a)$. Non-tangential cones are comparable to the Euclidean cones discussed in the introduction, in the sense that each such Euclidean cone contains a non-tangential cone and is contained in a non-tangential cone. It follows that $y$ is the non-tangential limit of $\{g_j(0)\}_{j=1}^{\infty}$ if and only if there is an $a$ such that $y$ lies in infinitely many spherical caps $\text{Cap}\,(g_j(0), a)$. Taking a countable union over wider and wider opening angles of the cones, the conical limit set is given by

$$L_c(G) = \bigcup_{l=2}^{\infty} \bigcap_{k=1}^{\infty} \bigcup_{j=k}^{\infty} \text{Cap}\,(g_j(0), l). \tag{2.3}$$

A set $K$ in the complex plane is *uniformly perfect* if there is a constant $c > 0$ such that for each $z_0 \in K$ and for all $r$ such that $0 < r < \operatorname{diam} K$,

$$K \cap \{z \mid cr \le |z - z_0| \le r\} \ne \emptyset. \tag{2.4}$$

(Here $\operatorname{diam} K$ is the diameter of $K$.) In other words, there is an upper bound on the moduli of annuli lying in the complement of $K$.

The *Poincaré series* for a discrete group $G$ is the series

$$\sum_{g \in G} (1 - |g(0)|)^s, \tag{2.5}$$

where $s$ is positive. There is a critical exponent $\delta = \delta(G)$, called the *exponent of convergence*, such that the series converges for all $s > \delta$ and diverges for all $s < \delta$.

Let $\Omega$ be a domain in $\overline{\mathbf{C}}$, let $E$ be a Borel subset of $\partial\Omega$, and let $z$ be a point in $\Omega$. The *harmonic measure at $z$ of $E$ in $\Omega$* is the Perron solution $w(z) = w(z, E, \Omega)$ of the Dirichlet problem in $\Omega$ for the boundary values $\mathbf{1}_E$. As a function of $z$, $w(z, E, \Omega)$ is harmonic in $\Omega$. For fixed $z$, $w(z, E, \Omega)$ is a probability measure on $\partial\Omega$; $w(z, E, \Omega)$ is the probability that a Brownian traveller from $z$ first hits $\partial\Omega$ in the set $E$. Harmonic measure is monotonic in the domain $\Omega$: if $z \in \Omega_1 \subset \Omega_2$ and $E \subset \partial\Omega_1 \cap \partial\Omega_2$, then $w(z, E, \Omega_1) \le w(z, E, \Omega_2)$. It is also monotonic in $E$: if $E_1 \subset E_2 \subset \partial\Omega$, then $w(z, E_1, \Omega) \le w(z, E_2, \Omega)$.

The hyperbolic metric on the disc $\mathbb{D}$ is given by the element of arclength $ds = \frac{2}{1-|z|^2}|dz|$. The Möbius transformations $g : \mathbb{D} \to \mathbb{D}$ are the isometries for this metric. They are of the form

$$g(z) = e^{i\theta}\, \frac{z - z_0}{1 - \overline{z}_0 z} \tag{2.6}$$

where $\theta \in [0, 2\pi]$ and $z_0 \in \mathbb{D}$. The hyperbolic geodesics are the *orthocircular arcs*: that is, arcs of circles which meet the unit circle $\partial\mathbb{D}$ at right angles. The hyperbolic distance from $0$ to $a \in \mathbb{D}$ is

$$d_{\mathrm{hyp}}(0, a) = \int\limits_0^{|a|} \frac{2}{1 - |z|^2}\, |dz| = \log\left[\frac{1 + |a|}{1 - |a|}\right]. \tag{2.7}$$

Let $\Omega$ be a domain in $\overline{\mathbf{C}}$ whose universal covering space is the disc, and let $\pi : \mathbb{D} \to \Omega$ be the covering map. The hyperbolic metric on $\mathbb{D}$ can be projected via $\pi$ to a metric on $\Omega$ given by $\lambda_\Omega(w)|dw| = \lambda_\Omega(\pi(z))|\pi'(z)||dz| = \frac{2}{1-|z|^2}|dz|$, where $w = \pi(z)$. The boundary $\partial\Omega$ of $\Omega$ is uniformly perfect if and only if there is a constant $c_\Omega > 0$ such that the function $\lambda_\Omega$ satisfies

$$\frac{c_\Omega}{\operatorname{dist}(w, \partial\Omega)} \le \lambda_\Omega(w) \le \frac{2}{\operatorname{dist}(w, \partial\Omega)}. \tag{2.8}$$

See [BP, Po2]. Here dist denotes Euclidean distance.

In dimension two, there is a connection between quasiconformal maps and doubling measures. Let $\varphi : \mathbf{C} \to \mathbf{C}$ be a $K$-quasiconformal map which preserves the upper half plane. Then $\varphi$ maps $\mathbb{R}$ to itself, and the restriction $f$ of $\varphi$ to $\mathbb{R}$ is an increasing homeomorphism which satisfies

$$\frac{1}{M} \leq \frac{f(x+t) - f(x)}{f(x) - f(x-t)} \leq M \tag{2.9}$$

with a constant $M$ depending on $K$, for all $x \in \mathbb{R}$ and all $t > 0$. Functions with this property are called *M-quasisymmetric*. Quasisymmetry is a necessary and sufficient condition for an increasing homeomorphism $f$ of $\mathbb{R}$ to be the boundary values of a quasiconformal mapping which preserves the upper half plane [BA]. Note that different quasiconformal mappings may have the same boundary values.

Let $f : \mathbb{R} \to \mathbb{R}$ be an increasing homeomorphism. Define a measure $\mu$ on $\mathbb{R}$ by setting $\mu([a,b]) = f(b) - f(a)$ for intervals $[a,b]$. Then $f$ is $M$-quasisymmetric if and only if $\mu$ satisfies $1/M \leq \mu(I)/\mu(J) \leq M$ whenever $I$ and $J$ are adjacent intervals of equal length. In other words, $f$ is $M$-quasisymmetric if and only if $\mu$ is doubling.

To summarize, there is a many-to-one correspondence between quasiconformal self-mappings of the upper half plane and doubling measures on the real line. There is a similar correspondence between quasiconformal self-mappings of the disc and doubling measures on the circle.

Some references for the material above are [N1] and [L].

Two quantities $A$ and $B$ are *comparable*, denoted $A \sim B$, if there is a constant $c > 0$ such that $\frac{1}{c} B \leq A \leq c B$. We also write $A \overset{c}{\sim} B$ if $A$ and $B$ are comparable with constant $c$. Throughout the paper we normalize Lebesgue measure so that the unit $n$-sphere has measure one.

In the figures we denote the unit circle $\partial \mathbb{D}$ by $S^1$.

### 3. Proof of Theorem 1.1 assuming Theorem 1.2

Fix a number $\lambda$ with $0 \leq \lambda \leq 1$. Let $\mu$ be the measure given by Theorem 1.2, and take $c > 0$ such that $\mu(2I) \leq c\,\mu(I)$ for all arcs $I$ in the unit circle $\partial\mathbb{D}$. Set $\nu(\cdot) = \lambda\mu(\cdot) + (1 - \lambda)\,|\cdot|$, where $|\cdot|$ denotes normalized Lebesgue measure on $\partial\mathbb{D}$. Then $\nu$ is a doubling measure on the circle, since for each arc $I$ in $\partial\mathbb{D}$,

$$\begin{aligned}
\nu(2I) &= \lambda\,\mu(2I) + (1 - \lambda)\,|2I| \\
&\leq c\,\lambda\,\mu(I) + 2\,(1 - \lambda)\,|I| \\
&\leq \max{(c, 2)}\,\nu(I).
\end{aligned} \tag{3.1}$$

The measure $\nu$ assigns the fraction $\lambda$ of its mass to the conical limit set of $G$.

Define $f : \partial\mathbb{D} \to \partial\mathbb{D}$ so that $\nu([a, b]) = f(b) - f(a)$ for all arcs $[a, b]$ in $\partial\mathbb{D}$. $f$ is a quasisymmetric mapping of the circle to itself, and $|f(E)| = \nu(E)$ for any measurable subset $E$ of the circle. The map $f$ may be extended to a quasiconformal mapping $\varphi$ of the unit disc onto itself [BA, DE] and further to a quasiconformal self-mapping of $\overline{\mathbf{C}}$. Conjugating the Fuchsian group $G$ of Theorem 1.2 by $\varphi$ yields a quasiconformal group $\Gamma = \varphi \circ G \circ \varphi^{-1}$ which acts on the disc.

It remains to show that the conical limit set of $\Gamma$ is the image under $\varphi$ of the conical limit set of the Fuchsian group $G$. This is a special case of a much more general fact; see [M] for example. We give a direct proof. It is sufficient to show that a quasiconformal map $\varphi : \mathbb{D} \to \mathbb{D}$ maps any non-tangential cone into another non-tangential cone, perhaps with larger opening angle. Since $\varphi^{-1}$ is also quasiconformal, this implies that $L_c(\Gamma) = \varphi(L_c(G))$.

Define the *cross ratio* of any four distinct points $\alpha$, $\beta$, $\gamma$, and $\delta$ in $\overline{\mathbf{C}}$ by

$$\tau = |\alpha, \beta, \gamma, \delta| = \frac{|\alpha - \gamma|}{|\alpha - \delta|} \cdot \frac{|\beta - \delta|}{|\beta - \gamma|}, \tag{3.2}$$

with the convention that if $\delta = \infty$,

$$|\alpha, \beta, \gamma, \infty| = \frac{|\alpha - \gamma|}{|\beta - \gamma|}. \tag{3.3}$$

Let $\varphi : \overline{\mathbf{C}} \to \overline{\mathbf{C}}$ be $K$-quasiconformal. Let $\varphi(\tau) = |\varphi(\alpha), \varphi(\beta), \varphi(\gamma), \varphi(\delta)|$. Väisälä has shown [V] that there is an increasing function $\Phi_K : (0, \infty) \to (0, \infty)$, depending only on $K$, such that $\varphi(\tau) \leq \Phi_K(\tau)$ for all $\tau$.

Let $z$ be a point in Cone $(y, a) = \{z \in \mathbb{D} \mid |z - y| \leq a\,(1 - |z|)\}$, where $y$ is a point in $\partial\mathbb{D}$ and $1 < a < \infty$. Let $w$ be the point in $\partial\mathbb{D}$ such that $1 - |\varphi(z)| = |\varphi(w) - \varphi(z)|$.

Suppose $\varphi^{-1}(\infty) = \infty$. Then

$$
\begin{aligned}
\frac{|\varphi(y) - \varphi(z)|}{1 - |\varphi(z)|} &= \frac{|\varphi(y) - \varphi(z)|}{|\varphi(w) - \varphi(z)|} \\
&= |\varphi(y), \varphi(w), \varphi(z), \infty| \\
&\leq \Phi_K(|y, w, z, \infty|) \\
&= \Phi_K\left(\frac{|y - z|}{|w - z|}\right) \\
&\leq \Phi_K\left(\frac{|y - z|}{1 - |z|}\right) \\
&\leq \Phi_K(a).
\end{aligned}
\tag{3.4}
$$

Here the second last inequality holds because $|w - z| \geq 1 - |z|$ and $\Phi_K$ is increasing. Therefore $\varphi(z)$ lies in Cone $(\varphi(y), \Phi_K(a))$.

Now suppose $\varphi^{-1}(\infty) \neq \infty$. Let $w$ and $y$ be any two points in the unit circle. Since $\varphi^{-1}(\infty)$ is not in the closed unit disc, we have

$$
|w - \varphi^{-1}(\infty)| \leq \frac{2 + d}{d} |y - \varphi^{-1}(\infty)|
\tag{3.5}
$$

where $d = \text{dist}(\varphi^{-1}(\infty), \overline{\mathbb{D}}) > 0$. Then

$$
\begin{aligned}
\frac{|\varphi(y) - \varphi(z)|}{1 - |\varphi(z)|} &= |\varphi(y), \varphi(w), \varphi(z), \infty| \\
&\leq \Phi_K(|y, w, z, \varphi^{-1}(\infty)|) \\
&= \Phi_K\left(\frac{|y - z|}{|w - z|} \cdot \frac{|w - \varphi^{-1}(\infty)|}{|y - \varphi^{-1}(\infty)|}\right) \\
&\leq \Phi_K\left(\frac{|y - z|}{1 - |z|} \cdot \frac{2 + d}{d}\right) \\
&\leq \Phi_K(a \cdot \frac{2 + d}{d}).
\end{aligned}
\tag{3.6}
$$

Hence $\varphi(z)$ lies in Cone $(\varphi(y), \Phi_K(a \cdot \frac{2+d}{d}))$.

It follows that $L_c(\Gamma) = \varphi(L_c(G))$, and $|L_c(\Gamma)| = |\varphi(L_c(G))| = \nu(L_c(G)) = \lambda$. In other words $\Gamma$ is a quasiconformal group, acting on the unit disc, whose conical limit set supports the fraction $\lambda$ of the mass of Lebesgue measure on the unit circle. This completes the proof of Theorem 1.1. $\qquad\square$

## 4. Outline of the proof of Theorem 1.2

**Theorem 1.2.** *There exist a Fuchsian group $G$ and a doubling measure such that the conical limit set $L_c(G)$ of $G$ has one-dimensional Lebesgue measure zero, and the doubling measure is supported on $L_c(G)$.*

In this section we present an outline of the proof of Theorem 1.2. Let $K \subset [0, 1]$ be the classical ternary Cantor set, and let $\Omega = \overline{\mathbf{C}} \setminus K$ be the complement of $K$ in the extended complex plane. Let $G$ be the Fuchsian group, acting on the disc $\mathbf{D}$, which uniformizes $\Omega$; in other words $G$ is the covering group of $\Omega$ and $\mathbf{D}/G$ is conformally equivalent to $\Omega$. Let $\pi : \mathbf{D} \to \Omega$ be the covering map, normalized so that $\pi(0) = \infty$. The conical limit set $L_c(G)$ is the set of non-tangential accumulation points of the orbit $\{g(0)\}_{g \in G}$. For this group $G$, the conical limit set has measure zero. (See Section 5.) Our aim is to construct a doubling measure which is supported on $L_c(G)$.

We begin with a model example of a doubling measure $\mu$, obtained as the limit of a sequence of measures $\mu_n$ of the form $d\mu_n = F_n \cdots F_1 \, dx$, such that $\mu$ is supported on a set of Lebesgue measure zero. This type of example is due to Kahane [K]. Divide the unit interval $I_0$ into 5 equal subintervals $I_{1,1}, \ldots, I_{1,5}$, $1 \leq k \leq 5$, numbered from left to right, and let

$$F_1(x) = \begin{cases} 1, & \text{for } x \in I_{1,1} \cup I_{1,5}; \\ 1/2, & \text{for } x \in I_{1,2} \cup I_{1,4}; \text{ and} \\ 2, & \text{for } x \in I_{1,3}. \end{cases} \qquad (4.1)$$

Then

.  The large fraction 2/5 of the area under the graph of $F_1$ lies above the middle fifth $I_{1,3}$ of $I_0$.

.  $F_1$ has mean value one on $I_0$;

.  $F_1$ is a positive function whose values lie between 1/2 and 2; and

.  $F_1 \equiv 1$ on subintervals of length $|I_0|/5$ at each end of $I_0$.

Let $F_n(x) = F_{n-1}(5x \bmod 1)$, for $n \geq 2$. On each of the $5^{n-1}$ subintervals $I_{n-1,j}$ of $I_0$ of length $5^{-n+1}|I_0|$, $F_n$ is a dilation of $F_1$, with the same properties on $I_{n-1,j}$ as $F_1$ has on $I_0$. We refer to the intervals $I_{n,j}$ as 5-*ary intervals*. The effect of $F_n$ is to concentrate the large fraction 2/5 of the mass assigned to each $I_{n-1,j}$ onto the middle fifth of $I_{n-1,j}$.

The measure $\mu_n$ given by $d\mu_n = F_n \cdots F_1 \, dx$ satisfies $1/4 \leq \mu_n(I)/\mu_n(J) \leq 4$ for all adjacent 5-ary intervals $I$ and $J$ of equal length. Therefore $\mu_n$ is doubling. The measures $\mu_n$ converge weakly to a unique limit $\mu$, which is also doubling. The support of this measure is a set $E$ of Hausdorff dimension less than one, concentrated near the centres of the 5-ary intervals which appear in the construction. (More precisely, this model set $E$ consists of those points $x$ in $I_0$ such that the proportions of 0's, 1's, 2's, 3's, and 4's in the base five expansion of $x$ are asymptotically equal to 1/5, 1/10, 2/5, 1/10, and 1/5 respectively. Its Hausdorff dimension is $1 - \frac{1}{5} \frac{\log 2}{\log 5} \approx 0.914$.)

It is convenient to recast this construction in terms of a random walk on a tree. We use a 5-ary tree as a discrete model of the disc. Put a vertex at the origin, with five equally spaced edges emanating from it. Put a vertex at the end of each edge; call these the *children* of the original vertex. Add five more edges emanating from each new vertex.

Repeat, with the new edges always pointing outwards towards the circle. Each branch of the tree corresponds to a point of the circle; one can think of the tip of the branch landing somewhere on the circle.

We distinguish a proper subset $V$ of vertices of the tree. Each vertex $v$ in the tree has five children, $v_1, \ldots, v_5$; let $V$ consist of all middle children $v_3$.

Consider a particle performing a random walk on the vertices of the tree, with the following constraints:

  . the particle may jump from one vertex to another if and only if the second is a
    child of the first: in other words they are connected by an edge and the second
    vertex is *below* the first (i.e. closer to the circle); and
  . when the particle is at the vertex $v$, it jumps to the children $v_k$, $1 \leq k \leq 5$,
    of $v$ with *jump probabilities* $p_k$, $1 \leq k \leq 5$, respectively, where $p_1 = p_5 = 1/5$,
    $p_2 = p_4 = 1/10$, and $p_3 = 2/5$.

The particle always moves away from the origin, towards the tips of the branches.

We recover the 5-ary intervals by projecting the tree radially outwards onto the circle, so that each vertex corresponds to an interval in the circle. (The construction of the measure $\mu$ given above can be done on the unit circle instead of the unit interval.) Let the origin correspond to the whole circle. Project the five vertices which are the children of the origin to five intervals of equal length, which together cover the circle. And so on: for each vertex $v$, project the children of $v$ to five equal subintervals of the interval corresponding to $v$. The collection of all these intervals forms a 5-adic *grid* of nested intervals in the circle; we refer to these intervals as *grid intervals*. We say a grid interval is in the $n^{\text{th}}$ *layer* of the grid if it is the projection of a vertex which is separated from the origin by $n$ edges in the tree.

Each path which the particle may take, given by a sequence of vertices $v_1 \rightarrow v_2 \rightarrow \cdots$ in the tree, corresponds to a point $x$ in the circle, and also to a nested, decreasing sequence of grid intervals $I_1 \supset I_2 \supset \cdots \ni x$ in the circle.

Clearly the jump probabilities $\{p_k\}$ and the functions $\{F_n\}$ defined on the grid intervals encode the same information. We can recover the functions $\{F_n\}$ as follows. Let $v$ be a vertex which is $n-1$ edges away from the origin. Let $v_k$, $1 \leq k \leq 5$, be the children of $v$, let $p_k$ be the probability of jumping from $v$ to $v_k$, and let $I$ and $I_k$ be the grid intervals which are the projections of $v$ and $v_k$ respectively. Then, on $I$, $F_n$ is simply the function which is constant on each $I_k$ and satisfies

$$p_k = \text{Prob} \, (v \rightarrow v_k) = \frac{\int_{I_k} F_n \, dx}{\int_I F_n \, dx} \qquad (4.2)$$

for $1 \leq k \leq 5$.

In this setting, the measure $\mu$ which is the limit of the $\mu_n$'s defined by $d\mu_n = F_n \cdots F_1 \, dx$ is the hitting measure for the random walk on the tree: the measure of a subset $A$ of the circle is exactly the probability that a particle which starts from the origin will hit the circle somewhere in $A$. Note that the support of $\mu$ contains those points in the circle $\partial \mathbb{D}$ which correspond to paths in the tree which pass through infinitely many vertices from $V$.

As noted above, the support of $\mu$ is a small set $E$ which is distributed in a very regular way. If we had a Fuchsian group whose conical limit set were exactly this set $E$, then we would have established Theorem 1.2. However, it is not apparent how to build a Fuchsian group with pre-determined conical limit set. Instead, we start with the particular Fuchsian group described above, and modify the Kahane-type construction to produce a doubling measure whose support $E$ is contained in $L_c(G)$.

In Section 6 we generalize the Kahane-type construction, showing that one obtains a doubling measure under rather mild conditions on the functions $\{F_n\}$. In particular, the 5-adic grid of intervals in the example above may be replaced by one where each interval $I_{n-1,j}$ at the $(n-1)^{\text{th}}$ level is divided into finitely or infinitely many intervals $I_{n,k}$ at the $n^{\text{th}}$ level, such that $F_n$ is constant on each $I_{n,k}$. It suffices (Lemma 6.3) to assume that

. $F_n$ has mean value one on each $I_{n-1,j}$ at the $(n-1)^{\text{th}}$ level;

. the functions $F_n$ are positive and uniformly bounded away from zero and infinity; and

. $F_n \equiv 1$ on subintervals at each end of each interval $I_{n-1,j}$, whose lengths are at least some definite fraction (independent of $n$ and $j$) of the length of $I_{n-1,j}$.

Then the measures $\mu_n$ given by $d\mu_n = F_n \cdots F_1 \, dx$ converge to a doubling measure $\mu$.

The idea is now to distinguish a collection $V$ of grid intervals in a new, less regular grid, such that points in $\partial \mathbb{D}$ which lie in infinitely many of the intervals in $V$ are actually conical limit points, and to define the functions $F_n$ so that they satisfy the three hypotheses above, are large on the intervals in $V$, and are small elsewhere. The result is a doubling measure which assigns full measure to $L_c(G)$.

In the setting of random walks on trees, our goal is to define a tree in the disc, in which the orbit points $\{g(0)\}_{g \in G}$ appear as a proper subset $V$ of the vertices, and to define the jump probabilities for the edges in this tree so that the subset of paths which hit $V$ infinitely often carries full measure. Here the measure on the space of paths is that which gives, to each family of paths with a common initial segment, measure equal to the product of the jump probabilities for the edges in that segment.

The branches of the tree which contain infinitely many vertices from $V$ correspond to points in the circle which lie in the limit set of $G$. This is almost by definition: for each such point $x$ in the circle there is a sequence of infinitely many orbit points which approaches $x$. When we implement this abstract procedure by specifying the vertices and edges in the tree, we will get something more: these points $x$ will actually be *conical* limit points of $G$. This is shown in Section 15.

We make a minor simplification of the goal stated above. We impose on the particle the additional constraint that:

. whenever the particle reaches a vertex in $V$, it stops there.

Suppose the particle is started at some vertex $v$, with this stopping condition. Let $V(v)$ be the subset of vertices in $V$ which can be reached from $v$. In other words, $V(v)$ consists of those $w \in V$ below $v$ such that there are no other vertices from $V$ between $v$ and $w$. Then it suffices to define the jump probabilities so that, from each starting vertex $v$, the particle reaches $V(v)$ with probability one. We show below how to do this; then by removing the stopping condition we obtain the measure we require.

In the setting of the grid of intervals, this simplification means that it is sufficient to define the functions $F_n$ so that, beginning with Lebesgue measure on any grid interval $I$, we obtain a doubling measure supported on the maximal grid intervals in $I$ which correspond to orbit points $g(0)$ in the tree.

In the rest of this section we describe:

. a tree in the disc, in which the orbit points $\{g(0)\}_{g \in G}$ appear as a proper subset $V$ of the vertices;

. the adjacencies in the tree (which pairs of vertices are connected by edges);

. a projection $P$ from vertices to grid intervals in the circle;

. a definition of $F_n$ (equivalently, of the jump probabilities); and

. a model calculation showing that this definition of $F_n$ implies that the resulting measure $\mu$ is supported on $L_c(G)$.

We now implement the scheme described above, making extensive use of the well-understood geometry of the complement of the Cantor set $K$.

Let $D$ be the closure of a simply connected domain in the disc such that $\pi(D)$ is either the upper or lower half plane. We call such a region $D$ a *half fundamental domain* for $\Omega$. The boundary of $D$ consists of a closed subset $\partial D \cap \partial \mathbb{D}$ of the unit circle, and a union of disjoint arcs of circles which meet the unit circle at right angles. We call these *orthocircular arcs*. The covering map $\pi$ maps $\partial D \cap \partial \mathbb{D}$ onto the Cantor set $K$, and it maps the orthocircular arcs in $\partial D$ onto the components of $\overline{\mathbb{R}} \setminus K$. The normal fundamental domain $\mathcal{F}$ for $\Omega$ is the union of the two half fundamental domains whose boundaries contain the origin. Its images $\{g(\mathcal{F})\}_{g \in G}$ tile the disc. We use instead the tiling in which each fundamental domain $g(\mathcal{F})$, except $\mathcal{F}$ itself, is divided into the two half fundamental domains which it contains.

Divide each orthocircular arc in the boundaries of the half fundamental domains in this tiling into infinitely many intervals, all of comparable hyperbolic length. These intervals will be the vertices of the tree. A convenient way to divide the arcs is to make a Whitney decomposition of the components of $\overline{\mathbb{R}} \setminus K$, in the Riemann surface $\Omega$, and then to pull it back to the disc via the branches of the inverse of the covering map.

To this end, let $L$ be a component of $[0,1] \setminus K$. $L$ is an open interval of length $3^{-k}$, for some $k \geq 1$. Partition $L$ into subintervals $J$ so that there are $2^j$ $J$'s of length $3^{-j}|L|$, for each $j \geq 1$, and so that each $J$ satisfies the following Whitney condition:

$$|J| = 2 \operatorname{dist}(J, K). \tag{4.3}$$

This can be done by putting the two $J$'s of length $|J| = 3^{-1}|L|$ in the middle of $L$, and then adjoining the smaller $J$'s at each end in decreasing order of size. The Euclidean length of each $J$ is a negative power of 3.

For the component $L = \overline{\mathbb{R}} \setminus [0,1]$ of $\overline{\mathbb{R}} \setminus K$, we make a similar decomposition, except that there is one interval $J_\infty$ which contains the point at infinity. Fix a small positive number $\sigma$, and let $J_\infty = \overline{\mathbb{R}} \setminus (-\sigma, 1+\sigma)$. Partition the intervals $(-\sigma, 0]$ and $[1, 1+\sigma)$ into $J$'s satisfying (4.3) so that the length of each $J$ is a negative power of 3 and the lengths decrease as the $J$'s approach $[0,1]$. Make the convention that $|J_\infty| = 1/3$.

Let the vertices of the tree be the images $I$ of these intervals $J$ under all branches of the inverse of the covering map. Also, let $V$ consist of those intervals $I$ in the disc such that $\pi(I) = J_\infty$. The orbit $\{g(0)\}_{g \in G}$ is mapped by the covering map to the point at infinity. So each interval $I$ such that $\pi(I) = J_\infty$ contains an orbit point, and each orbit point appears in the tree as a vertex in the set $V$.

The random walk on the tree in the disc can be pushed forward by the covering map to a random walk in $\Omega = \overline{\mathbf{C}} \setminus K$, with the intervals $J$ as vertices. Our goal becomes: to define the jump probabilities for this random walk so that if the particle starts at any initial vertex, it reaches $J_\infty$ with probability one.

The idea is that $J_\infty$ is the largest interval (in the Euclidean sense), and we define the $\{p_k\}$ so that from any interval, the most likely outcome is that the particle jumps to a much larger interval (in the Euclidean sense).

For each interval $J$, we specify the intervals $J'$ to which the particle may jump from $J$. The particle may jump only to the intervals $J'$ in a particular segment $E_J$ of $\overline{\mathbb{R}}$ near to $J$. See Figure 4 in Section 9 below. This segment $E_J$ is much larger than $J$. We define $E_J$ as follows for most grid intervals $J$, and we call these $J$ *standard*. The Euclidean length of $E_J$ is $|E_J| = 2 \cdot 3^N |J|$, where $N$ is a fixed large integer, *independent of $J$*. $E_J$ comprises one of the scaled copies of the Cantor set $K$, of length $3^N |J|$, which arise in the construction of $K$, together with an adjacent gap in $K$, also of length $3^N |J|$. By a *gap* in $K$, we mean an interval which is a component of $\overline{\mathbb{R}} \setminus K$. The segment $E_J$ contains infinitely many intervals $J'$. Counting, we see that the two largest have Euclidean length $|J'| = 3^{N-1} |J|$, and there are $2^k$ intervals $J'$ in $E_J$ of length $|J'| = 3^{N-k} |J| = 2^{-1} \cdot 3^{-k} |E_J|$, for each $k \geq 1$.

For non-standard grid intervals $J$, the segments $E_J$ take different forms. In particular, from some intervals $J$ it is possible to jump directly to $J_\infty$. (Full definitions of the segments $E_J$ for standard and non-standard intervals are made in Section 9.)

There are now infinitely many intervals $J'$ accessible from $J$. Pulling back to the disc, we obtain a tree in which each vertex has infinite valence. We define the adjacencies in the tree, in other words which pairs of vertices are joined by edges, as follows. Let $I$ be a vertex in the tree, and let $J = \pi(I)$. Pull $E_J$ back to the disc by the branch of $\pi^{-1}$ which takes $J$ to $I$. One can show (Section 10) that this pre-image $\pi^{-1}(E_J)$ lies roughly below $I$, that is, between $I$ and the nearest part of the unit circle. The Whitney intervals $J'$ in $E_J$ are pulled back to vertices $I' = \pi^{-1}(J')$ in $\pi^{-1}(E_J)$. Let these $I'$ be the children of $I$; put an edge between $I$ and each such $I'$.

We also use the segments $\pi^{-1}(E_J)$ to define the grid intervals corresponding to the vertices of the tree. We define the projection $P$ from the vertices to the grid intervals as follows. Since the endpoints of $E_J$ lie in $K$, the endpoints of $\pi^{-1}(E_J)$ lie in the unit circle. Let $P(I)$ be the arc of the unit circle which has the same endpoints as $\pi^{-1}(E_J)$ and which lies below $\pi^{-1}(E_J)$; in other words $P(I)$ is the shorter arc of the unit circle between the endpoints of $\pi^{-1}(E_J)$. This interval $P(I)$ is the grid interval corresponding to the vertex $I$. The projection $P$ is almost radial, but not quite.

We define some auxiliary functions $X_i$ on the circle which keep track of the size of the Whitney intervals $J$ along a path in $\Omega$ travelled by the particle. (See Section 12.) Let $I_0$ be a grid interval in the circle, corresponding to a vertex $v_0$ in the tree. Let $x$ be a point in $I_0$. It corresponds to a path $v_0 \to v_1 \to v_2 \to \cdots$ in the tree; to a nested, decreasing sequence of grid intervals $I_0 \supset I_1 \supset I_2 \supset \cdots \ni x$ in the circle; and to a path $J_0 \to J_1 \to J_2 \to \cdots$ in $\Omega$. Define

$$X_1(x) = \log_3 \left[ \frac{|J_1|}{|J_0|} \right]. \qquad (4.4)$$

The logarithm is simply in order to make $X_1$ integer-valued; recall that the Euclidean lengths of the intervals $J$ are all powers of 3. $X_1$ is positive if and only if $J_1$ is larger than $J_0$. Define

$$X_i(x) = \begin{cases} \log_3 \left[ \frac{|J_i|}{|J_{i-1}|} \right], & \text{if } J_1, \ldots, J_{i-1} \neq J_\infty; \\ 1, & \text{otherwise,} \end{cases} \qquad (4.5)$$

for $i \geq 1$, with the convention that $|J_\infty| = 1/3$. This function measures the difference in size of the intervals where the particle is before and after the $i^{\text{th}}$ jump, unless the particle has already reached $J_\infty$. If the particle reaches $J_\infty$, then it stays there, and all subsequent $X_i$'s simply record a 1.

Let

$$S_k(x) = \sum_{i=1}^{k} X_i(x). \qquad (4.6)$$

Then $S_k(x) \to +\infty$ if and only if the path in $\Omega$ which corresponds to $x$ reaches $J_\infty$ (and therefore stops there).

We wish to ensure that a particle starting from any vertex in the tree reaches $V$ with probability one. Let $v$ be a vertex, and let $I_0$ be the corresponding grid interval. Let $V(I_0)$ denote the union of the maximal grid intervals in $I_0$ which correspond to $J_\infty$. These maximal grid intervals are the projections to the circle of those vertices $v'$ in $V$ which can be reached from $v$ by a particle which stops when it first hits $V$. Then

$$V(I_0) = \{ x \in I_0 \mid S_k(x) \to +\infty \}, \qquad (4.7)$$

and our problem is reduced to defining the functions $\{F_n\}$ in such a way that $S_k(x) \to +\infty$ a.e. $d\mu$ on $I_0$, where $\mu$ is the limit of the measures $\mu_n$, and $d\mu_n = F_n \cdots F_1\, dx$.

To do this, it is sufficient to define the $\{F_n\}$ so that the following estimates hold:

$$EX_i = \frac{1}{\mu(I)} \int_I X_i(x)\, d\mu \geq c_1, \qquad (4.8)$$

and

$$EX_i^2 = \frac{1}{\mu(I)} \int_I X_i^2(x)\, d\mu \leq c_2, \qquad (4.9)$$

with uniform positive constants $c_1$ and $c_2$, for all grid intervals $I$ in the $(n-1)^{\text{th}}$ layer of the grid, for all $n \geq 1$. Here $X_i$ is the auxiliary function, defined in (4.5) above, for a particle which has not yet reached $J_\infty$ and which makes its $(i-1)^{\text{th}}$ jump from the vertex corresponding to $I$. In other words, the expectations of the $X_i$'s with respect to $\mu$ are uniformly large, and the second moments are uniformly bounded. Then an appeal to the strong law of large numbers for martingales (see Section 15) yields the result that $S_k(x) \to +\infty$ a.e. $d\mu$ on $I_0$, as required.

In the remainder of this section, we show how to define $F_n$ for those grid intervals $I$ which correspond to standard intervals $J$. (Full definitions of $F_n$ are given in Section 11, for standard intervals, and in Section 14, for non-standard intervals.) We also give a model calculation for the estimate (4.8) on the expectation of $X_i$.

Let $I$ be a grid interval, from the $(n-1)^{\text{th}}$ layer of the grid, whose corresponding interval $J$ in $\overline{\mathbb{R}} \setminus K$ is standard. The grid intervals $I'$ in $I$ from the next $(n^{\text{th}})$ layer of the grid correspond to those vertices in the tree which are accessible from the vertex coresponding to $I$. Therefore these $I'$ are exactly those grid intervals which correspond to the intervals $J'$ in $E_J$ in the Riemann surface $\Omega$. In other words, we pull back the segment $E_J$ from $\overline{\mathbb{R}}$ to the disc via the appropriate branch of the inverse of the covering map, and then project it to the interval $I$ in the circle using the projection $P$ defined above.

The key idea is that the geometry of the segment $E_J$ is not distorted very much by this operation. In fact, we now assume that it is not distorted at all, so that the $n^{\text{th}}$ layer of grid intervals in $I$ has exactly the same structure as $E_J$. Namely, one half of $I$ contains an exact scaled copy $\widetilde{K}$ of the ternary Cantor set $K$, of length $|I|/2$, and the other half of $I$ corresponds to a gap in $K$. Again, we are assuming that the intervals $I'$ in $I \setminus \widetilde{K}$ are exactly Whitney intervals with $|I'| = 2\,\text{dist}\,(I', \widetilde{K})$, just like the Whitney intervals $J'$ in $E_J$. In particular, $I \setminus \widetilde{K}$ contains $2^k$ intervals $I'$ of length $|I'| = 2^{-1} \cdot 3^{-k}|I|$, for each $k \geq 1$.

Much of the paper (especially in Sections 8, 10, 11, 13, and 14) is devoted to proving estimates which control the distortion of the segments $E_J$ by $P \circ \pi^{-1}$. These estimates rely on the fact that on the boundary of a chord-arc domain, one can use harmonic measure to estimate arclength (via $A_\infty$-equivalence). For the remainder of this section, however, we simply assume that there is no distortion.

Let $B_1$ denote the union of the two largest grid intervals $I'$ in $I$ from the $n^{\text{th}}$ layer of the grid; they satisfy $|I'| = |I|/6$. Fix a small positive number $\varepsilon$, independent of $I$. Define the function $F_n$ on $I$ so that:

.  $F_n \equiv (1-\varepsilon)\,|I|/|B_1|$ on $B_1$;

.  $F_n \equiv 1$ on small segments $I_l$ and $I_r$ at each end of $I$, where $|I_l|/|I|$ and $|I_r|/|I|$ are fixed quantities independent of $I$, and the endpoints of $I_l$ and $I_r$ are endpoints of intervals from the $n^{\text{th}}$ layer of the grid; and

.  $F_n \equiv \delta$ on the rest of $I$, where $\delta$ is a small number chosen to ensure that $F_n$ has mean value 1 on $I$.

The force of the first condition is that $\mu(B_1)/\mu(I) = 1-\varepsilon$; the effect of $F_n$ is to concentrate almost all the mass of $I$ onto $B_1$. This is because of the way that $\mu$ is constructed: the $\mu$-measure of an interval from the $n^{\text{th}}$ layer of the grid is simply its $\mu_n$-measure, since all

$F_l$ with $l > n$ have mean value 1 on the grid interval. So

$$\begin{aligned}
\frac{\mu(B_1)}{\mu(I)} &= \frac{\mu_n(B_1)}{\mu_{n-1}(I)} \\
&= \frac{\int_{B_1} F_n(x) \cdots F_1(x) \, dx}{\int_I F_{n-1}(x) \cdots F_1(x) \, dx} \\
&= (1 - \varepsilon) \, \frac{|I|}{|B_1|} \, \frac{\int_{B_1} F_{n-1}(x) \cdots F_1(x) \, dx}{\int_I F_{n-1}(x) \cdots F_1(x) \, dx} \\
&= 1 - \varepsilon,
\end{aligned} \tag{4.10}$$

since the product $F_{n-1} \cdots F_1$ is constant on $I$.

In Sections 11 and 14 we show that these functions $F_n$ are uniformly bounded away from zero and infinity, which is a hypothesis of Lemma 6.3 on doubling measures.

The second condition above is another hypothesis of Lemma 6.3. However, in order to simplify the model calculation below, we temporarily forget that $F_n = 1$ near the ends of $I$, and assume that $F_n = \delta$ on all of $I \setminus B_1$.

Let

$$\begin{aligned}
B_k &= \{x \in I \mid X_i(x) = N - k\} \\
&= \bigcup \{I' \subset I \mid |I'| = 2^{-1} \cdot 3^{-k}|I|\},
\end{aligned} \tag{4.11}$$

for $k \geq 1$. We know that $B_k$ consists of exactly $2^k$ intervals $I'$, each of length $2^{-1} \cdot 3^{-k}|I|$, and so

$$|B_k| = \frac{1}{2} \left(\frac{2}{3}\right)^k |I|. \tag{4.12}$$

The function $X_i$ takes the constant value $N - k$ on $B_k$, for $k \geq 1$. The function $F_n$ is also constant on each $B_k$; we are assuming that $F_n \equiv \delta$ on $I \setminus B_1$, and forgetting that actually $F_n = 1$ on the parts of some $B_k$'s which are very near the endpoints of $I$. We know that $\mu(B_k) = \mu_n(B_k)$, and $\mu(I) = \mu_{n-1}(I)$. Therefore, in computing the ratio $\mu(B_k)/\mu(I)$, we need only take account of the value of $F_n$ on $B_k$; the $F_l$'s with $l \neq n$ are irrelevant. We see that $\mu(B_k)/\mu(I)$ is just $\delta\,|B_k|/|I|$, if $k \geq 2$. Therefore

$$
\begin{aligned}
EX_i &= \frac{1}{\mu(I)} \int_I X_i(x)\,d\mu \\
&= \frac{1}{\mu(I)} \sum_{k=1}^{\infty} \int_{B_k} X_i(x)\,d\mu \\
&= \sum_{k=1}^{\infty} (N-k)\,\frac{\mu(B_k)}{\mu(I)} \\
&= (N-1)\,\frac{\mu(B_1)}{\mu(I)} + \delta \sum_{k=2}^{\infty} (N-k)\,\frac{|B_k|}{|I|} \\
&\geq (N-1)\,(1-\varepsilon) - \delta \sum_{k=2}^{\infty} k\,\frac{|B_k|}{|I|} \\
&= (N-1)\,(1-\varepsilon) - \delta \sum_{k=2}^{\infty} k\,\frac{1}{2}\left(\frac{2}{3}\right)^k.
\end{aligned}
\tag{4.13}
$$

The last series converges, so we can ensure that $EX_i \geq c_1 > 0$ by choosing $N$ sufficiently large with respect to $\delta$.

When we correct for the fact that $F_n = 1$, not $\delta$, near the endpoints of $I$, we get some other terms involving similar convergent series, and we can still make the whole argument work by choosing $N$ large. More crucially, we prove estimates (see Sections 10, 13, and 14, especially Lemmas 13.1 and 14.1) showing that, even allowing for the distortion of lengths produced by the inverse of the covering map and by the projection from the disc to the circle, the ratio $|B_k|/|I|$ is controlled by $\beta^k$, where $\beta$ is smaller than 1 and independent of $I$, so that the relevant series converge and a calculation like the one above can be done.

The estimate $EX_i^2 \leq c_2$ is similar. The estimates (4.8) and (4.9) are established in Section 13 for standard intervals, and in Section 14 for non-standard intervals.

## 5. Choice of the Fuchsian group

Let $K \subset [0,1]$ be the classical ternary Cantor set. Let $\Omega = \overline{\mathbf{C}} \setminus K$ be the complement of $K$ in the extended complex plane. $\Omega$ is an infinitely connected planar Riemann surface. Since $K = \partial \Omega$ contains more than two points, the universal covering space of $\Omega$ is the unit disc $\mathbb{D}$, and $\Omega$ is conformally equivalent to $\mathbb{D}/G$, where $G$ is a Fuchsian group. We choose this infinitely generated Fuchsian group $G$ as our example.

In this section we show that the conical limit set $L_c(G)$ of $G$ has one-dimensional Lebesgue measure zero, and moreover the Hausdorff dimension of $L_c(G)$ satisfies $1/2 \leq \dim(L_c(G)) < 1$. The remainder of the paper is devoted to the construction of a doubling measure which is supported on $L_c(G)$.

Since $K$ has positive logarithmic capacity, Green's function exists for $\Omega$, and the result cited in the introduction on the measure of the conical limit set of a Fuchsian group implies that $|L_c(G)| = 0$. Fernández [F] has shown that if the boundary $\partial\Omega$ of a planar domain $\Omega$ is uniformly perfect, then the exponent of convergence of the corresponding Fuchsian group is strictly less than one. The Cantor set $K$ is uniformly perfect. For Fuchsian groups the exponent of convergence is equal to the Hausdorff dimension of the conical limit set [Pa1, S1]. (See also [BJ], where this is proved in much greater generality.) We may therefore conclude that in fact $\dim(L_c(G)) < 1$.

A Fuchsian group $G$ is of *fully accessible type* if there is a measurable fundamental set $B \subset \partial\mathbb{D}$ for the action of $G$ on $\partial\mathbb{D}$. In other words, $B$ contains no two $G$-equivalent points, and $\partial\mathbb{D}$ is a.e. equal to the union of the $G$-images of $B$: $\sum_{g \in G} |g(B)| = 1$. $G$ is of *accessible type* if there is a measurable set $B \subset \partial\mathbb{D}$, containing no two $G$-equivalent points, such that $\sum_{g \in G} |g(B)| > 0$.

The normal fundamental domain $\mathcal{F}_0$ (see Section 7) for $\Omega = \overline{\mathbf{C}} \setminus K$ has the property that $|\partial\mathcal{F}_0 \cap \partial\mathbb{D}| = 0$. We outline the well-known proof. The covering map $\pi$ takes $\partial\mathcal{F}_0$ to $[0,1]$. Since $\partial\mathcal{F}_0$ is a rectifiable curve, $\pi$ preserves sets of zero length. (This is the F. and M. Riesz theorem; see [Po3].) Since $K$ has Lebesgue measure zero, $\partial\mathcal{F}_0 \cap \pi^{-1}(K) = \partial\mathcal{F}_0 \cap \partial\mathbb{D}$ has measure zero.

This implies that the group $G$ is not of accessible type [Po1]. On the other hand, Patterson [Pa2] showed that if $\delta(G) < 1/2$ then $G$ is of fully accessible type. It follows that our group $G$ must have $\delta(G) = \dim(L_c(G)) \geq 1/2$.

## 6. A Lemma on Doubling Measures

In this section we follow the philosophy first laid out by Kahane [K] to construct doubling measures on an interval $I_0$. Our measure is the limit of a sequence of measures whose densities are step functions. In Lemma 6.3 we show that mild conditions on these step functions ensure that the limit measure is doubling.

Simple constructions of this type can yield doubling measures which are supported on very small sets. In particular there are examples of such measures which are supported on sets of arbitrarily small, but positive, Hausdorff dimension. Such a measure is equivalent to the derivative of a quasisymmetric mapping that takes a set of small Hausdorff dimension to a set of full one-dimensional Lebesgue measure.

The two interrelated ingredients of this construction are a grid of nested subintervals of $I_0$ and a sequence of suitable density functions.

**Definition 6.1.** A *grid* of subintervals of an interval $I_0$ is a collection $\mathcal{H} = \bigcup_{n=1}^{\infty} \mathcal{H}_n$ of subintervals of $I_0$ satisfying

   i) $\mathcal{H}_0 = I_0$;

   ii) for each $n \geq 0$, the intervals in $\mathcal{H}_n$ have disjoint interiors, and $|I_0 \setminus \bigcup_{I \in \mathcal{H}_n}| = 0$; and

   iii) for each $n \geq 1$, for each interval $J \in \mathcal{H}_n$ there is an interval $I \in \mathcal{H}_{n-1}$ such that $J \subset I$.

The collections $\mathcal{H}_n$ are called the *layers* of intervals in the grid $\mathcal{H}$.

**Definition 6.2.** A function $F$ defined on a interval $I$ is a *step function* if it is constant on each interval in some finite or infinite partition of $I$. We say $F$ is $(\delta, \eta)$-*suitable for $I$* if $F$ is a step function satisfying

   i) $\frac{1}{|I|} \int_I F(x)\, dx = 1$;

   ii) $0 < \delta \leq F(x) \leq 1/\delta$ for all $x \in I$; and

   iii) $F(x) \equiv 1$ on subintervals $I_l$ and $I_r$ of $I$ such that $I_l$ has the same left endpoint as $I$, $I_r$ has the same right endpoint as $I$, $|I_l| \geq \eta|I|$, and $|I_r| \geq \eta|I|$.

See Figure 1.

Fix an interval $I_0$, and numbers $\delta$ and $\eta$ such that $0 < \delta \leq 1$ and $0 < \eta \leq 1/2$. We construct a sequence $\{\mu_n\}$ of measures on $I_0$, whose densities are of the form $d\mu_n = F_n(x) \cdots F_1(x)\, dx$. We define the functions $F_n$ simultaneously with a related grid $\mathcal{H}$ of subintervals of $I_0$. Each $F_n$ is $(\delta, \eta)$-suitable for each interval from the $(n-1)^{\text{th}}$ layer $\mathcal{H}_{n-1}$ of the grid.

Figure 1. $F(x)$ is $(\delta, \eta)$-suitable for $I = [0, 1]$.

Let $\mu_0$ be Lebesgue measure on $I_0$. Let $F_1$ be any function which is $(\delta, \eta)$-suitable for $I_0$. Define $\mu_1$ by $d\mu_1 = F_1(x)\, dx$. Since $F_1$ has mean value one on $I_0$, $\mu_1$ has the same total mass as $\mu_0$. Take any partition of $I_0$ into subintervals $J$ such that $F_1$ is constant on each $J$. Note that we allow $F_1$ to take the same value on adjacent $J$'s; the $J$'s need not be maximal. Also, there may be finitely or infinitely many $J$'s. Let $\mathcal{H}_1$, the first layer of intervals in the grid, be the collection of subintervals in this partition.

We define the measures $\mu_n$ inductively for $n \geq 2$. Suppose we have already defined functions $F_1, \ldots, F_{n-1}$, measures $\mu_1, \ldots, \mu_{n-1}$, and the layers $\mathcal{H}_1, \ldots, \mathcal{H}_{n-1}$ of the grid. In particular, $F_{n-1}$ is constant on each interval $I \in \mathcal{H}_{n-1}$. Let $F_n$ be any function which is $(\delta, \eta)$-suitable for each interval $I \in \mathcal{H}_{n-1}$. Define $\mu_n$ by

$$d\mu_n = F_n(x)\, d\mu_{n-1} = F_n(x) \cdots F_1(x)\, dx. \tag{6.1}$$

Again, the fact that $F_n$ has mean value one implies that $\mu_n$ has the same total mass as $\mu_0$. To define the next layer $\mathcal{H}_n$ of the grid, we partition each $I \in \mathcal{H}_{n-1}$ into subintervals $J$ such that $F_n$ is constant on each $J$. Let $\mathcal{H}_n$ be the collection of all these subintervals $J$.

Note that the collection $\mathcal{H} = \bigcup_{n=1}^\infty \mathcal{H}_n$ of subintervals of $I_0$ does form a grid according to the definition above: each interval $J \in \mathcal{H}_n$ is contained in some $I \in \mathcal{H}_{n-1}$, and each layer $\mathcal{H}_n$ is a partition of $I_0$, so $|I_0 \setminus \bigcup_{I \in \mathcal{H}_n} I| = 0$ and the intervals in $\mathcal{H}_n$ have disjoint interiors. Also, on each $J \in \mathcal{H}_n$, $\mu_n$ is a constant multiple of Lebesgue measure: $d\mu_n = (\mu_n(J)/|J|)\, dx$ on $J$. Furthermore, the total mass assigned to an interval $J \in \mathcal{H}_n$ does not change after the $n^{\text{th}}$ step: $\mu_n(J) = \mu_{n+1}(J) = \cdots$.

**Lemma 6.3.** *Fix an interval $I_0$ and numbers $\delta$ and $\eta$ with $0 < \delta \leq 1$ and $0 < \eta \leq 1/2$. Let $\{\mu_n\}$ be a sequence of measures and $\mathcal{H}$ a grid of subintervals of $I_0$ as described above, so that for each $n \geq 1$, $d\mu_n = F_n(x) \cdots F_1(x)\, dx$ and $F_n(x)$ is $(\delta, \eta)$-suitable for each interval* $I \in \mathcal{H}_{n-1}$. Then the sequence $\{\mu_n\}$ converges weakly to a measure $\mu$ on $I_0$ which has the same total mass as $\mu_0$ and which is doubling with a constant depending only on $\delta$ and $\eta$.

**Proof.** The $\mu_n$ are a sequence of measures, all with the same finite total mass, on a compact interval $I_0$. Therefore there is a subsequence converging weakly to a measure $\mu$ on $I_0$ with the same mass. It is an exercise to see that this limit measure is unique; in other words the full sequence $\mu_n$ converges to $\mu$. The $\mu$-measure of any interval $I \in \mathcal{H}_n$ is given by $\mu(I) = \mu_n(I)$.

Let $J$ and $K$ be adjacent intervals, of equal length, contained in $I_0$. If, for each $n$, the function $F_n$ is constant on $J \cup K$, then $\mu(J) = \mu(K)$. Otherwise, let $m$ be the first integer such that $F_m$ is not constant on $J \cup K$. In particular,

$$\mu_m(J) \overset{1/\delta}{\sim} \mu_{m-1}(J) = \mu_{m-1}(K) \overset{1/\delta}{\sim} \mu_m(K). \tag{6.2}$$

We show that $\mu(J) \sim \mu_m(J)$ and $\mu(K) \sim \mu_m(K)$.

We may assume that $F_m$ is not constant on $J$. Let $x \in J$ be a point where the value of $F_m$ changes. Then for each $n \geq m + 1$ there is a neighbourhood of $x$ on which $F_n \equiv 1$.

Let $J_r$ be the part of $J$ to the right of $x$. If $J_r$ does not contain any whole interval from $\mathcal{H}$, then $F_n \equiv 1$ on $J_r$ for all $n \geq m + 1$, and so $\mu(J_r) = \mu_m(J_r)$. Otherwise, let

$$\mathcal{L}_m = \bigcup_j \{ I_{m,j} \in \mathcal{H}_m \mid I_{m,j} \subset J_r \}. \tag{6.3}$$

Then $\mu(\mathcal{L}_m) = \mu_m(\mathcal{L}_m)$. (Therefore if $J_r$ consists only of whole intervals from $\mathcal{H}_m$, in other words if $J_r = \mathcal{L}_m$, then $\mu(J_r) = \mu_m(J_r)$ and we are done.)

The function $F_{m+1}$ is constant on each interval in $\mathcal{H}_{m+1}$. Suppose there is an interval $I_{m+1,j} \in \mathcal{H}_{m+1}$ which meets $J_r \setminus \mathcal{L}_m$ and on which $F_{m+1} \neq 1$. The interval $I_{m,k} \in \mathcal{H}_m$ which contains $I_{m+1,j}$ is adjacent to $\mathcal{L}_m$, or has $x$ as its left endpoint in case $\mathcal{L}_m$ is empty, and it contains $J_r \setminus \mathcal{L}_m$. Let $I$ be the largest subinterval of $I_{m,k}$, with the same left endpoint as $I_{m,k}$, on which $F_{m+1} \equiv 1$. $I$ is a union of intervals in $\mathcal{H}_{m+1}$, $|I| \geq \eta\, |I_{m,k}|$, and $I \subset J_r \setminus \mathcal{L}_m \subset I_{m,k}$. Then

$$\begin{aligned}
\mu(J_r \setminus \mathcal{L}_m) &\geq \mu_{m+1}(I) \\
&= \mu_m(I) \\
&= \frac{\mu_m(I_{m,k})}{|I_{m,k}|}\, |I| \\
&\geq \eta\, \mu_m(I_{m,k}) \\
&\geq \eta\, \mu_m(J_r \setminus \mathcal{L}_m),
\end{aligned} \tag{6.4}$$

and

$$\begin{aligned}
\mu(J_r \setminus \mathcal{L}_m) &\leq \mu_m(I_{m,k}) \\
&\leq \frac{\mu_m(I_{m,k})}{|I_{m,k}|}\, \eta^{-1}\, |I| \\
&= \eta^{-1}\, \mu_m(I) \\
&\leq \eta^{-1}\, \mu_m(J_r \setminus \mathcal{L}_m).
\end{aligned} \tag{6.5}$$

Therefore $\mu(J_r \setminus \mathcal{L}_m) \overset{1/\eta}{\sim} \mu_m(J_r \setminus \mathcal{L}_m)$, and so $\mu(J_r) \overset{1/\eta}{\sim} \mu_m(J_r)$.

On the other hand, suppose there is no interval in $\mathcal{H}_{m+1}$ which meets $J_r \setminus \mathcal{L}_m$ and on which $F_{m+1} \neq 1$. Then $F_{m+1} \equiv 1$ on $J_r \setminus \mathcal{L}_m$. If $J_r \setminus \mathcal{L}_m$ contains an interval $I_{m+1,j}$ from $\mathcal{H}_{m+1}$ on which $F_{m+1} \equiv 1$, then the $\mu$-measure of $I_{m+1,j}$ is not only equal to its $\mu_{m+1}$-measure but also to its $\mu_m$-measure. Let

$$\mathcal{L}_{m+1} = \bigcup_j \left\{ I_{m+1,j} \in \mathcal{H}_{m+1} \ \Big| \ I_{m+1,j} \subset J_r \setminus \mathcal{L}_m; \ F_{m+1} \equiv 1 \text{ on } I_{m+1,j} \right\}. \tag{6.6}$$

be the union of all such intervals. $\mathcal{L}_{m+1}$ is an interval to the right of the interval $\mathcal{L}_m$ and adjacent to it, if both are non-empty. Now $\mu(\mathcal{L}_{m+1}) = \mu_m(\mathcal{L}_{m+1})$. If $J_r = \mathcal{L}_m \cup \mathcal{L}_{m+1}$, then $\mu(J_r) = \mu_m(J_r)$ and we are done. Otherwise, we have reduced the problem to estimating the $\mu$-measure of $J_r \setminus (\mathcal{L}_m \cup \mathcal{L}_{m+1})$.

Apply this argument to $J_r \setminus \bigcup_{n=m}^{m+i} \mathcal{L}_n$ for $i \geq 1$, where $\mathcal{L}_{m+i}$ is defined inductively by

$$\mathcal{L}_{m+i} = \bigcup_j \left\{ I_{m+i,j} \in \mathcal{H}_{m+i} \ \Big| \ I_{m+i,j} \subset J_r \setminus \left( \bigcup_{n=m}^{m+i-1} \mathcal{L}_n \right); \ F_{m+i} \equiv 1 \text{ on } I_{m+i,j} \right\}. \tag{6.7}$$

Suppose that for some integer $i \geq 1$ there is an interval $I_{m+i+1,j} \in \mathcal{H}_{m+i+1}$ which meets $J_r \setminus \bigcup_{n=m}^{m+i} \mathcal{L}_n$ and on which $F_{m+i+1} \neq 1$. Let $l$ be the first such integer. As in the case $l = 1$ above, let $I_{m+l,k}$ be the interval from $\mathcal{H}_{m+l}$ which contains $I_{m+l+1,j}$ and let $I$ be the largest subinterval of $I_{m+l,k}$ on which $F_{m+l+1} \equiv 1$. Then $I \subset J_r \setminus \bigcup_{n=m}^{m+l} \mathcal{L}_n \subset I_{m+l,k}$ and so, by the analogues of (6.4) and (6.5), $\mu(J_r) \overset{1/\eta}{\sim} \mu_m(J_r)$, and we are done.

On the other hand, suppose there is no such integer $l \geq 1$. If the intervals $\mathcal{L}_{m+i}$, $i \geq 1$, exhaust $J_r$, then $\mu(J_r) = \mu_m(J_r)$, and we are done. Otherwise, $J_r \setminus \bigcup_{n=m}^{\infty} \mathcal{L}_n$ is a non-empty interval which has the same right endpoint as $J_r$ and which contains no whole interval from $\mathcal{H}$. So $F_{m+i} \equiv 1$ on $J_r \setminus \bigcup_{n=m}^{\infty} \mathcal{L}_n$ for all $i \geq 1$. Therefore $\mu(J_r \setminus \bigcup_{n=m}^{\infty} \mathcal{L}_n) = \mu_m(J_r \setminus \bigcup_{n=m}^{\infty} \mathcal{L}_n)$, and so $\mu(J_r) = \mu_m(J_r)$. This is the last possible case; we have shown that $\mu(J_r) \overset{1/\eta}{\sim} \mu_m(J_r)$ in all cases.

The same argument shows that $\mu(J_l) \overset{1/\eta}{\sim} \mu_m(J_l)$, where $J_l$ is the part of $J$ to the left of $x$. Therefore $\mu(J) \overset{1/\eta}{\sim} \mu_m(J)$.

If the function $F_m$ is not constant on the adjacent interval $K$ then $\mu(K) \overset{1/\eta}{\sim} \mu_m(K)$ by the argument above. Also, if $F_n \equiv 1$ on $K$ for all $n \geq m+1$, then $\mu(K) \overset{1/\eta}{\sim} \mu_m(K)$. Otherwise, either $F_n$ is constant (although not necessarily $\equiv 1$) on $K$ for each $n \geq m+1$, or we let $l \geq m+1$ be the first integer such that $F_l$ is not constant on $K$. In these cases it is enough to estimate the number of integers $n \geq m+1$ for which $F_n \equiv$ constant $\neq 1$ on $K$. We show that this number is bounded by a constant $c = c(\eta)$. Then in the first case $\mu(K) \overset{c'}{\sim} \mu_m(K)$, and in the second $\mu(K) \overset{1/\eta}{\sim} \mu_l(K) \overset{c'}{\sim} \mu_m(K)$, where the constants $c'$ depend on $c$ and on $\delta$.

Suppose there are at least $k$ integers $n_1 < \cdots < n_k$ greater than $m$ for which $F_{n_j} \equiv$ constant $\neq 1$ on the interval $K$. We may assume that $K$ is to the right of $J$. For each $j$ with $1 \leq j \leq k$, there is an interval $I_j \in \mathcal{H}_{n_j}$ which contains $K$ and whose left endpoint lies in $J$. These intervals are nested: $I_1 \supset \cdots \supset I_k$. Let $L_j$ (respectively $R_j$) be the largest subinterval of $I_j$, with the same left (respectively right) endpoint as $I_j$, on which $F_{n_j} \equiv 1$, and let $M_j = I_j \setminus (L_j \cup R_j)$. Then $|M_j| \leq (1 - 2\eta)|I_j|$, and $I_j \subset M_{j-1}$ for $2 \leq j \leq k$. See Figure 2.

Figure 2. Graphs of $F_{n_1+1}(x)$ and $F_{n_2+1}(x)$.

Also, $L_1 \subset J$ and $\eta |I_1| \leq |L_1| \leq |J|$. So

$$|K| \leq |M_k| \leq (1 - 2\eta) |I_k| \leq (1 - 2\eta) |M_{k-1}|. \tag{6.8}$$

After $k$ steps we see that

$$|K| \leq (1 - 2\eta)^k |I_1| \leq (1 - 2\eta)^k \eta^{-1} |J| = (1 - 2\eta)^k \eta^{-1} |K|. \tag{6.9}$$

Therefore $k \leq \log \eta / \log (1 - 2\eta)$.

We have shown that $\mu(J) \overset{c}{\sim} \mu(K)$, where $c$ depends only on $\delta$ and $\eta$, for every pair $J$, $K$ of adjacent intervals of equal length. In other words $\mu$ is a doubling measure, and the doubling constant depends only on $\delta$ and $\eta$. $\qquad\square$

## 7. Fundamental domains for $\Omega$

We begin the construction of the doubling measure by describing the fundamental domains for $\Omega = \overline{\mathbf{C}} \setminus K$, where $K \subset [0,1]$ is the classical ternary Cantor set. $G$ is the Fuchsian group such that $\mathbb{D}/G$ is conformally equivalent to $\Omega$. $G$ acts as the covering group of $\Omega$ on the universal covering space $\mathbb{D}$ of $\Omega$. Let $\pi : \mathbb{D} \to \Omega$ be the covering map; it is a many-to-one conformal mapping which takes each $G$-orbit in $\mathbb{D}$ to a single point in $\Omega$. Normalize $\pi$, by composition with a Möbius transformation of the disc, so that $\pi(0) = \infty$.

A *fundamental domain* for $\Omega$ is a domain $\mathcal{F}$ in $\mathbb{D}$ such that $\overline{\mathcal{F}}$ contains at least one point from every $G$-orbit, and $\mathcal{F}$ contains no two $G$-equivalent points. The domain $\Omega = \overline{\mathbf{C}} \setminus K$ is a *Denjoy domain*; that is, it is the complement in $\overline{\mathbf{C}}$ of a closed linear set. Denjoy domains have fundamental domains which have a particularly simple form; namely, they are *orthocircular*. This means that the boundary is of the form $\partial \mathcal{F} = E \cup \bigcup_{n=1}^{\infty} A_n$, where $E$ is a closed subset of the unit circle $\partial \mathbb{D}$, and the $A_n$'s (called *orthocircular arcs*) are disjoint arcs of circles which meet $\partial \mathbb{D}$ at right angles. The covering map takes the set $\partial \mathcal{F} \cap \partial \mathbb{D}$ to $\partial \Omega$, and it takes the orthocircular arcs comprising $\partial \mathcal{F} \setminus \partial \mathbb{D}$ to the components of $\overline{\mathbb{R}} \setminus \partial \Omega$. The images $\{g(\mathcal{F})\}_{g \in G}$ are orthocircular domains which tile the disc. (See [RR]; they consider the case where the boundary $\partial \Omega$ has positive length, but a Denjoy domain $\Omega$ with $|\partial \Omega| = 0$ also has an orthocircular fundamental domain $\mathcal{F}$, with $|\partial \mathcal{F} \cap \partial \mathbb{D}| = 0$.)

We refer to the images in the disc of the upper and lower half planes under the branches of $\pi^{-1}$ as *half fundamental domains* for $\Omega = \overline{\mathbf{C}} \setminus K$. Each half fundamental domain $D$ is orthocircular. To fix ideas, suppose that $\pi(D)$ is the upper half plane. The boundary $\overline{\mathbb{R}}$ of the upper half plane lifts via $\pi^{-1}$ to the boundary of $D$. The Cantor set $K$ lifts to the closed subset $\partial D \cap \partial \mathbb{D}$ of the circle, and the open intervals which are the components of $\overline{\mathbb{R}} \setminus K$ lift to the orthocircular arcs in the boundary of $D$.

We now fix a tiling of the disc by half fundamental domains for $\Omega$, together with one fundamental domain for $\Omega$. Since $\pi(0) = \infty$, one of the preimages under $\pi$ of the open interval $\overline{\mathbb{R}} \setminus [0,1]$ is a diameter $\gamma$ of $\mathbb{D}$. Let $D$ be the preimage of the upper half plane whose boundary contains $\gamma$. Let $D'$ be the reflection of $D$ in $\gamma$. $D'$ is a preimage of the lower half plane. Then $D \cup \gamma \cup D'$ is a fundamental domain for $\Omega$; it contains the origin; and it is symmetric with respect to the diameter $\gamma$. We denote this fundamental domain by $\mathcal{F}_0$; it is, by symmetry, the normal fundamental domain for $\Omega$. Its image under the covering map is $\pi(\mathcal{F}_0) = \overline{\mathbf{C}} \setminus [0,1]$. Each open interval in $[0,1] \setminus K$ lifts to two orthocircular arcs in $\partial \mathcal{F}_0$.

Tile the remainder of the disc, $\mathbb{D} \setminus \mathcal{F}_0$, by the half fundamental domains which are the $G$-images of $D$ and $D'$. These half fundamental domains are indexed in a natural way by the number of orthocircular boundary arcs separating them from the origin. Write $\Omega_{1,j}$ for any half fundamental domain which is adjacent to $\mathcal{F}_0$; that is, which is separated from the origin $0$ by a single orthocircular arc. Similarly, write $\Omega_{n,j}$ for a half fundamental domain which is separated from $0$ by $n$ orthocircular arcs in the boundaries of the half fundamental domains in the fixed tiling of $\mathbb{D}$.

A note on orientation: we think of the direction from a point in the disc towards the unit circle as *down*, and of the direction towards the origin as *up*. Terms like *below* and *above* are to be understood in the same way; *below* means closer to the unit circle. Given a half fundamental domain $\Omega_{n,j}$, we refer to the large orthocircular arc, denoted $A_{n,j}$, which separates $\Omega_{n,j}$ from 0 as the *upper part of the boundary of $\Omega_{n,j}$*, and to $\partial \Omega_n \setminus A_{n,j}$ as the *lower part of the boundary of $\Omega_{n,j}$*.

Let $\mathcal{A}_1 = \partial \mathcal{F}_0 \setminus \partial \mathbb{D}$ be the union of the orthocircular arcs in the boundary of $\mathcal{F}_0$. Each arc $A_{1,j} \subset \mathcal{A}_1$ is the upper part of some $\partial \Omega_{1,j}$. For each $n \geq 2$, let $\mathcal{A}_n$ be the union, over all $j$, of the orthocircular arcs in the lower parts of the boundaries of the $\Omega_{n-1,j}$:

$$\mathcal{A}_n = \left[ \bigcup_j \partial \Omega_{n-1,j} \setminus \partial \mathbb{D} \right] \setminus \mathcal{A}_{n-1}. \tag{7.1}$$

The union $\mathcal{A} = \bigcup_{n=1}^{\infty} \mathcal{A}_n$ is the collection of all orthocircular boundary arcs which appear in the tiling of $\mathbb{D} \setminus \mathcal{F}_0$ by half fundamental domains.

We show in Section 10 that the orthocircular boundary arcs are uniformly hyperbolically separated: the hyperbolic distance between any two of these arcs is greater than a fixed positive constant. It follows from this observation that the boundaries of the half fundamental domains are chord-arc curves, with a uniform chord-arc constant. We will not need this fact, but we will need to show (Sections 13 and 14) that certain subsets of the half fundamental domains are chord-arc with uniform chord-arc constants.

## 8. Whitney decompositions

In this section we make a Whitney decomposition of the components of $\overline{\mathbb{R}} \setminus K$, and show that it lifts via $\pi^{-1}$ to a Whitney decomposition of the orthocircular arcs in the boundaries of the half fundamental domains in $\mathbb{D}$. The intervals in this lifted decomposition are the vertices in the tree described in Section 4.

The Cantor set $K$ is formed by removing the open middle third $\left( \frac{1}{3}, \frac{2}{3} \right)$ from $\left[ 0, 1 \right]$, then removing the open middle thirds $\left( \frac{1}{9}, \frac{2}{9} \right)$ and $\left( \frac{7}{9}, \frac{8}{9} \right)$ from the closed intervals $\left[ 0, \frac{1}{3} \right]$ and $\left[ \frac{2}{3}, 1 \right]$ which remain, and so on. We refer to the closed intervals $[0,1]$, $\left[ 0, \frac{1}{3} \right]$, $\left[ \frac{2}{3}, 1 \right]$, $\left[ 0, \frac{1}{9} \right]$, $\left[ \frac{2}{9}, \frac{1}{3} \right]$, ... which appear in this procedure as *closed construction intervals* of $K$.

A *Whitney decomposition* of an open interval $L$ in $\mathbb{R}$ is a partition of $L$ into closed intervals $J$, with disjoint interiors, such that the Euclidean length of each $J$ is comparable to the Euclidean distance from $J$ to the nearest endpoint of $L$. In other words, there is a constant $c > 0$ such that $|J| \overset{c}{\sim} \operatorname{dist}(J, \mathbb{R} \setminus L)$ for each $J$. This also implies that for each point $x$ in such an interval $J$, the distance from $x$ to $\mathbb{R} \setminus L$ is comparable to the distance from $J$ to $\mathbb{R} \setminus L$. The intervals $J$ in a Whitney decomposition are called *Whitney intervals*. We will also deal with Whitney decompositions of orthocircular arcs; these are defined analogously, so that the length of each Whitney interval in the arc is comparable to its distance from the nearest end of the arc.

Write $\overline{\mathbb{R}} \setminus K$ as a countable union of disjoint open intervals; let $L$ denote any one of these intervals. We now fix a Whitney decomposition of each component $L$ of $\overline{\mathbb{R}} \setminus K$.

First consider a component $L = (a, b)$ of $[0, 1] \setminus K$. $L$ is an open interval of length $3^{-k}$, for some $k \geq 1$. Partition $L$ into subintervals of length $3^{-n-k}$, $n \geq 1$, as follows. For

$n \geq 1$, let

$$J_n = \left[ b - 3^{-n+1} \frac{|L|}{2}, \ b - 3^{-n} \frac{|L|}{2} \right], \qquad \text{and}$$

$$J_{-n} = \left[ a + 3^{-n} \frac{|L|}{2}, \ a + 3^{-n+1} \frac{|L|}{2} \right]. \tag{8.1}$$

Then $J_1$ is the interval whose left endpoint is the midpoint of $L$ and which extends two-thirds of the way towards the right endpoint of $L$; for $n \geq 2$, $J_n$ is the interval whose left endpoint is the right endpoint of $J_{n-1}$ and which extends two-thirds of the way towards the right endpoint of $L$; and for $n \geq 1$, $J_{-n}$ is the reflection of $J_n$ in the midpoint of $L$. The intervals $\{J_{\pm n}\}_{n \geq 1}$ form a Whitney decomposition of $L$, $|J_{\pm n}| = 2 \operatorname{dist}(J_{\pm n}, K)$, and $|J_{\pm n}| = 3^{-n} |L| = 3^{-n-k}$.

Let $N$ be a large integer; its value will be fixed below. It is convenient to amalgamate the $2N$ intervals $J_{-N}, \ldots, J_N$ which are closest to the midpoint of $L$ into a single interval: Let

$$J_c = J_{-N} \cup \cdots \cup J_N. \tag{8.2}$$

This interval has length $|J_c| = (1 - 3^{-N})|L|$, and its distance from $K$ is $2^{-1} \cdot 3^{-N} |L|$, so it satisfies $|J_c| = 2 (3^N - 1) \operatorname{dist}(J_c, K)$.

The intervals $\{J_{\pm n}\}_{n \geq N+1}$ and $J_c$ form a Whitney decomposition of $L$, with constant $2 (3^N - 1)$.

Now consider the remaining component $L = \overline{\mathbb{R}} \setminus [0,1]$. Fix the small number $\sigma = (2 \cdot 3^{N+2})^{-1}$; the reasons for this particular choice will become apparent later. Let

$$J_\infty = \overline{\mathbb{R}} \setminus (-\sigma, 1 + \sigma). \tag{8.3}$$

$J_\infty$ is an interval of large but finite hyperbolic length, which contains the point at infinity. For $n \geq N + 3$, let

$$J_n = \left[ -\frac{1}{2} \frac{1}{3^{n-1}}, \ -\frac{1}{2} \frac{1}{3^n} \right], \qquad \text{and}$$

$$J_{-n} = \left[ 1 + \frac{1}{2} \frac{1}{3^n}, \ 1 + \frac{1}{2} \frac{1}{3^{n-1}} \right]. \tag{8.4}$$

Then $J_{N+3} = [-\sigma, -\sigma/3]$, and for $n \geq N + 4$, $J_n$ is the interval whose left endpoint is the right endpoint of $J_{n-1}$ and which extends two-thirds of the way towards 0. For $n \geq N + 3$, $J_{-n}$ is the reflection of $J_n$ through the midpoint of $[0,1]$. The intervals $\{J_{\pm n}\}_{n \geq N+3}$ form a Whitney decomposition of $[-\sigma, 0) \cup (1, 1 + \sigma]$, in which $|J_{\pm n}| = 3^{-n} = 2 \operatorname{dist}(J_{\pm n}, K)$.

Make the convention that

$$|J_\infty| = 1/3. \tag{8.5}$$

With this convention, the intervals $\{J_{\pm n}\}_{n \geq N+3}$ and $J_\infty$ form a Whitney decomposition of $L = \overline{\mathbb{R}} \setminus [0,1]$, with constant $2 \cdot 3^{N+3}$.

We use the phrase *the Whitney decomposition of $\overline{\mathbb{R}} \setminus K$* to denote the collection of Whitney intervals $\{J_{\pm n}\}_{n \geq N+1}$ and $J_c$ in the components of $[0,1] \setminus K$, together with the Whitney intervals $\{J_{\pm n}\}_{n \geq N+3}$ and $J_\infty$ in $\overline{\mathbb{R}} \setminus [0, 1]$. This decomposition will remain fixed for the rest of the construction. We call the intervals $\{J_{\pm n}\}$ *standard* and the intervals $J_c$ and $J_\infty$ *non-standard*.

The inverse $\pi^{-1}$ of the covering map lifts the components of $\overline{\mathbb{R}} \setminus K$ to the orthocircular arcs $\bigcup_n \mathcal{A}_n$ which appear in the boundaries of the half fundamental domains in our tiling of $\mathbb{D}$. $\pi^{-1}$ lifts the Whitney decomposition of $\overline{\mathbb{R}} \setminus K$ to a decomposition of $\bigcup_n \mathcal{A}_n$ into closed intervals with disjoint interiors.

The main result of this section is that this is a Euclidean Whitney decomposition of the orthocircular arcs $\bigcup_n \mathcal{A}_n$.

**Lemma 8.1.** *The Whitney decomposition of $\overline{\mathbb{R}} \setminus K$ lifts via $\pi^{-1}$ to a Whitney decomposition of the orthocircular arcs $\bigcup_n \mathcal{A}_n$.*

**Proof.** Since $\partial\Omega = K$ is uniformly perfect, there is a constant $c_\Omega > 0$ such that the function $\lambda_\Omega$ satisfies

$$\frac{c_\Omega}{\text{dist}\,(w, K)} \leq \lambda_\Omega(w) \leq \frac{2}{\text{dist}\,(w, K)} \tag{8.6}$$

for all $w$ in $\Omega$, where the element of hyperbolic arclength in $\Omega$ is $ds = \lambda_\Omega(w)\,|dw|$. Combined with the relation $|J| = 2\,\text{dist}\,(J, K)$ for standard Whitney intervals $J$, this implies that $\ell_{\text{hyp}}(J) \sim 1$ with a uniform constant for all standard $J$ in $\overline{\mathbb{R}} \setminus K$. Similarly, since the ratio $|J_c|/\text{dist}\,(J_c, K)$ is the same for all non-standard intervals $J_c$, $\ell_{\text{hyp}}(J_c) \sim 1$ for all such $J_c$. Also, the hyperbolic length of $J_\infty$ depends only on the constant $\sigma$.

These observations imply that the intervals $I$ in the preimage of the Whitney decomposition of $\overline{\mathbb{R}} \setminus K$ are all of comparable hyperbolic length, because the covering map $\pi$ preserves hyperbolic length. So $1/r \leq \ell_{\text{hyp}}(I) \leq r$ for all $I$, where the constant $r$ depends on the uniformly perfect constant of $K$, on the constant of the Whitney decomposition of $\overline{\mathbb{R}} \setminus K$, and on the choice of $\sigma$ in the definition of $J_\infty$.

It follows that for each $I$,

$$\frac{1}{2r}(1 - |z_1|) \leq |I| \leq 2re^r\,(1 - |z_1|), \tag{8.7}$$

where $z_1$ is the endpoint of $I$ closest to $\partial\mathbb{D}$. Therefore the intervals $I$ form a Whitney decomposition of the orthocircular arcs $\bigcup_n \mathcal{A}_n$. $\qquad\square$

## 9. Construction of palm leaves and the grid of intervals

In this section we construct a grid $\mathcal{H}$ of subintervals of the unit circle, of the kind described in Section 6. We already have a collection of Whitney intervals which form a Whitney decomposition of the orthocircular arcs $\bigcup_n \mathcal{A}_n$ in the boundaries of the half fundamental domains in the tiling of the unit disc. The idea is to project these Whitney intervals onto the unit circle, using a projection map $P$ defined below. The projections of the intervals in $\mathcal{A}_n$ form the $n^{\text{th}}$ layer $\mathcal{H}_n$ of the grid. Heuristically, the projection $P$ is almost radial projection. Let $\hat{}$ denote the inverse of the projection $P$: given a grid interval $I$, $\hat{I}$ denotes the unique Whitney interval in $\bigcup_n \mathcal{A}_n$ such that $P(\hat{I}) = I$.

Recall that $\mathcal{A}_1 = \partial \mathcal{F}_0 \setminus \partial \mathbb{D}$ is the collection of orthocircular arcs in the boundary of the fundamental domain $\mathcal{F}_0$ which contains the origin; and that for $n \geq 2$, $\mathcal{A}_n$ consists of the collection of orthocircular boundary arcs, from the tiling of $\mathbb{D} \setminus \mathcal{F}_0$ by half fundamental domains, which are separated from 0 by exactly $n-1$ other orthocircular boundary arcs.

In order to define the projection $P$, we make a construction in $\Omega = \overline{\mathbf{C}} \setminus K$ and lift it to the disc via $\pi^{-1}$. An outline follows. Given any Whitney interval $\hat{I} \subset \mathcal{A}_n$, we must specify which of the intervals $\hat{I}_{n+1,k} \subset \mathcal{A}_{n+1}$ should project via $P$ to subintervals of $P(\hat{I})$. (In terms of the tree described in Section 4, whose vertices are the Whitney intervals in $\bigcup_n \mathcal{A}_n$, we are now specifying the adjacencies between vertices; in other words which pairs of vertices are connected by edges.) Let $A_{n,j}$ be the orthocircular arc containing $\hat{I}$, and let $\Omega_{n,j}$ be the half fundamental domain below $A_{n,j}$. The covering map $\pi$ maps $\hat{I}$ to a Whitney interval $J$ in $\overline{\mathbb{R}} \setminus K$; it maps $A_{n,j}$ to the component $L$ of $\overline{\mathbb{R}} \setminus K$ which contains $J = \pi(\hat{I})$; it maps $\Omega_{n,j}$ to either the upper or lower half plane; and it maps $\partial \Omega_{n,j}$ to $\overline{\mathbb{R}}$. See Figure 3.

To each interval $J$ in the Whitney decomposition of $L$ we associate an interval $E_J \subset \overline{\mathbb{R}} \setminus L$, in such a way that the $E_J$ for all $J$ in $L$ have disjoint interiors and their union is $\overline{\mathbb{R}} \setminus L$. (The rest of this section contains precise definitions of these intervals $E_J$.) These intervals $E_J$ are of a certain form; in particular, for most $J$ the Euclidean length $|E_J|$ is a large fixed multiple of $|J|$, and the endpoints of $E_J$ always lie in the Cantor set $K$.

Let $\tilde{I}$ be the segment of $\partial \Omega_{n,j}$ such that $\pi(\tilde{I}) = E_J$. $\tilde{I}$ lies in the lower part of the boundary of $\Omega_{n,j}$, and its endpoints lie in the unit circle. Let $P(\hat{I})$ be the arc of the unit circle which has the same endpoints as $\tilde{I}$ and which lies below $\tilde{I}$. With this definition the intervals $\hat{I}_{n+1,k} \subset \mathcal{A}_{n+1}$ whose projections $P(\hat{I}_{n+1,k})$ are contained in $P(\hat{I})$ are precisely those intervals $\hat{I}_{n+1,k}$ which are contained in $\tilde{I}$.

The projections $\{P(\hat{I}_{n,j})\}_j$ of the Whitney intervals $\hat{I}_{n,j}$ in $\mathcal{A}_n$ constitute the $n^{\text{th}}$ layer $\mathcal{H}_n$ of the grid of subintervals of the circle. With the definition of $P$ outlined above, each interval $P(\hat{I}_{n+1,k})$ in $\mathcal{H}_{n+1}$ is contained in some $P(\hat{I}_{n,j})$ from the previous layer $\mathcal{H}_n$. Also, since the intervals $E_J$ have pairwise disjoint interiors for all $J$ in each $L$, the $P(\hat{I}_{n,j})$ at the $n^{\text{th}}$ level have pairwise disjoint interiors. Finally, by the remarks in Section 5 above, the part of $\partial \mathbb{D}$ below the orthocircular arcs in $\mathcal{A}_n$ has full Lebesgue measure in $\partial \mathbb{D}$ for each $n \geq 1$, so $\bigcup_{\hat{I} \subset \mathcal{A}_n} P(\hat{I})$ has full Lebesgue measure in $\partial \mathbb{D}$ for each $n \geq 1$. In other words, the union of the intervals $P(\hat{I})$ in $\mathcal{H}_n$ covers $\partial \mathbb{D}$, up to a set of measure zero. Thus the layers $\mathcal{H}_n$ of intervals form a grid $\mathcal{H}$ of subintervals of $\partial\mathbb{D}$, according to the definition in Section 6.

If $\widehat{R}$ is a segment of the boundary of some half fundamental domain $\Omega_{n,j}$, consisting of a collection of whole Whitney intervals $\{\widehat{I}_{n+1,k}\}_k$, possibly together with a subset $E$ of $\partial\Omega_{n,j} \cap \partial\mathbb{D}$, then we define $P(\widehat{R})$ to be $E \cup \bigcup_k P(\widehat{I}_{n+1,k})$.

Figure 3. Fan of palm leaves $\Lambda_J$ based at the Whitney
intervals $J = \pi(\widehat{I})$ in $L$; and a preimage in $\mathbb{D}$ of one leaf $\Lambda_J$.

We now give the details of the construction of the grid $\mathcal{H}$ of intervals, in particular defining the intervals $E_J$ associated to the Whitney intervals $J$ in $\overline{\mathbb{R}} \setminus K$. Let $N >> 1$ be the same large integer as in the previous section. Until the last part of the paper we regard $N$ as fixed, and we make all our geometric constructions using this fixed value.

Let $L$ be a component of $\overline{\mathbb{R}} \setminus K$. We define a decomposition of the upper half plane into infinitely many regions $\Lambda_J$, indexed by the Whitney intervals $J$ in $L$. See Figure 3. We refer to these regions as *palm leaves based at the Whitney intervals $J$ in $L$*, and to the collection of these regions as a *fan of palm leaves, based at $L$*. The boundary of the leaf $\Lambda_J$ based at the Whitney interval $J$ consists of two intervals in $\overline{\mathbb{R}}$, $J$ itself and $E_J$, together with two non-intersecting semicircles which join the endpoints of $J$ to the endpoints of $E_J$ and which meet $\mathbb{R}$ at right angles. We call $J$ the *base* and $E_J$ the *tip* of the leaf $\Lambda_J$. The union over all $J$ in $L$ of the intervals $E_J$ covers $\overline{\mathbb{R}} \setminus L$.

Once the entire construction of the fans of palm leaves is complete, we reflect them in the real axis to obtain analogous fans in the lower half plane.

We distinguish two cases: when the interval $L$ is a component of $[0,1] \setminus K$, and when $L = \overline{\mathbb{R}} \setminus [0,1]$.

**Case 1:** Let $L$ be a component of $\overline{\mathbb{R}} \setminus K$ which lies in $[0,1]$ and has length $|L| = 3^{-l}$ for some integer $l \geq 1$. Write $L = (a, a+3^{-l})$. Let $K_l = [a-3^{-l}, a]$ and $K_r = [a+3^{-l}, a+2 \cdot 3^{-l}]$ be the closed construction intervals of $K$ of length $|K_l| = |K_r| = |L| = 3^{-l}$ such that $K_l$ is immediately to the left of $L$ and $K_r$ is immediately to the right of $L$. The Whitney intervals in $L$ are enumerated from left to right as $\ldots, J_{-N-2}, J_{-N-1}, J_c, J_{N+1}, J_{N+2}, \ldots$ Recall that $|J_{\pm n}| = 3^{-n}|L|$, for $n \geq N+1$, and $|J_c| = (1 - 3^{-N})|L|$.

We refer to the Whitney intervals $J_{\pm n}$ in $L$ such that $n \geq N+1$ as *standard* intervals. For each standard interval $J$ we define a leaf $\Lambda_J$ based at $J$ such that the tip $E_J$ of the leaf has length $|E_J| = 2 \cdot 3^N |J|$. We begin with the standard Whitney interval $J_{N+1}$, which has length $|J_{N+1}| = 3^{-N-1}|L|$ and is in the right half of $L$.

Let $E_{J_{N+1}}$ be the closed interval which has the same right endpoint as $K_r$ and whose length is two-thirds the length of $K_r$. $E_{J_{N+1}}$ consists of a closed construction interval of $K$ and the closure of an adjacent open interval in $[0,1] \setminus K$ of the same length. Now $|E_{J_{N+1}}| = \frac{2}{3}|K_r| = \frac{2}{3}|L| = \frac{2}{3} \cdot 3^{N+1}|J_{N+1}| = 2 \cdot 3^N |J_{N+1}|$. Join the left endpoint of $J_{N+1}$ to the right endpoint of $E_{J_{N+1}}$ by a semicircle in the upper half plane. Join the right endpoint of $J_{N+1}$ to the left endpoint of $E_{J_{N+1}}$ in the same way. Let the leaf $\Lambda_{J_{N+1}}$ be the region in the upper half plane bounded by $J_{N+1}$, $E_{J_{N+1}}$, and the two semicircles. See Figure 4.

Figure 4. Standard palm leaf $\Lambda_J$ for $J = J_{N+1}$.
Not to scale: $|E_J| = 2 \cdot 3^N |J|$ and $N >> 1$.

Consider the next standard Whitney interval $J_{N+2}$; it has length $3^{-N-2}|L|$ and is immediately to the right of $J_{N+1}$. Construct the leaf $\Lambda_{J_{N+2}}$ and its tip $E_{J_{N+2}}$ exactly as for $J_{N+1}$: Let $E_{J_{N+2}}$ be the right two-thirds of $K_r \setminus E_{J_{N+1}}$, and join the endpoints of $J_{N+2}$ to those of $E_{J_{N+2}}$ by semicircles. The leaf $\Lambda_{J_{N+2}}$ bounded by $J_{N+2}$, $E_{J_{N+2}}$, and these semicircles is a copy of $\Lambda_{J_{N+1}}$, shrunk by a factor of one-third. The smaller semicircle in the boundary of $\Lambda_{J_{N+1}}$ is the larger semicircle in the boundary of $\Lambda_{J_{N+2}}$. Clearly $|E_{J_{N+2}}| = 2 \cdot 3^N |J_{N+2}|$.

Repeat this construction for the intervals $J_{N+k}$, $k \geq 1$, obtaining a sequence of leaves $\Lambda_{J_{N+k}}$ based at $J_{N+k}$ which satisfy $|E_{J_{N+k}}| = 2 \cdot 3^N |J_{N+k}|$. The leaves $\Lambda_{J_{N+k}}$ are all congruent to each other via dilations by powers of three. The intervals $E_{J_{N+k}}$, $k \geq 1$, have disjoint interiors and their union is $K_r$, the closed construction interval of length $|L|$ to the right of $L$.

Define leaves $\Lambda_{J_{-N-k}}$, $k \geq 1$, for the standard intervals at the left end of $L$, by reflecting the leaves $\Lambda_{J_{N+k}}$ through the perpendicular bisector of $L$. The leaf $\Lambda_{J_{-N-k}}$ based at $J_{-N-k}$ is the mirror image of the leaf $\Lambda_{J_{N+k}}$ based at $J_{N+k}$. The intervals $E_{J_{-N-k}}$, $k \geq 1$, cover $K_l$, the closed construction interval of length $|L|$ to the left of $L$.

We refer to the leaves we have just constructed for the standard intervals $J_{\pm(N+k)}$, $k \geq 1$, as *standard* leaves, and to their tips as *standard $E_J$'s*. We also refer to any grid interval $I$ such that $J = \pi(\widehat{I})$ is standard as a *standard* grid interval. For any standard leaf $\Lambda_J$, $|E_J| = 2 \cdot 3^N |J|$; the endpoints of $E_J$ lie in the Cantor set $K$; and $E_J$ consists of a closed construction interval of length $3^N |J|$ together with the closure of an adjacent open interval in $\overline{\mathbb{R}} \setminus K$, also of length $3^N |J|$.

It remains to define the leaf for the central interval $J_c = J_{-N} \cup \cdots \cup J_N$. Let $E_{J_c} = \overline{\mathbb{R}} \setminus (K_l \cup L \cup K_r)$. Define the leaf $\Lambda_{J_c}$ to be the region in the upper half plane bounded by $J_c$, $E_{J_c}$, and the two non-intersecting semicircles joining the endpoints of $J_c$ to those of $E_{J_c}$. The endpoints of $E_{J_c}$ lie in the Cantor set $K$. Notice that when $L = \left(\frac{1}{3}, \frac{2}{3}\right)$, $E_{J_c}$ is exactly $\overline{\mathbb{R}} \setminus [0, 1]$, since the tips of the leaves for the standard intervals in $\left(\frac{1}{3}, \frac{2}{3}\right)$ take up all of $\left[0, \frac{1}{3}\right] \cup \left[\frac{2}{3}, 1\right]$. For smaller intervals $L$, $E_{J_c}$ contains $\overline{\mathbb{R}} \setminus [0, 1]$ and part of $[0, 1]$. We use the term *non-standard* to refer to $J_c$, $\Lambda_{J_c}$, $E_{J_c}$, and any grid interval $I$ such that $\pi(\widehat{I}) = J_c$.

To summarize: given an open interval $L$ which is a component of $[0, 1] \setminus K$, we have defined a fan of palm leaves $\Lambda_J$ based at the Whitney intervals $J$ in $L$, such that:

. the $\Lambda_J$'s tile the upper half plane;

. the tips $E_J$ of the leaves are intervals whose union is $\overline{\mathbb{R}} \setminus L$;

. the $E_J$ have disjoint interiors;

. the endpoints of each $E_J$ lie in the Cantor set $K$;

. $|E_J| = 2 \cdot 3^N |J|$ for all standard $J$, i.e. for all but the central interval $J_c$ in $L$; and

. each standard $E_J$ consists of a closed construction interval of $K$ and the closure of an adjacent open interval in $\overline{\mathbb{R}} \setminus K$ of the same length.

**Case 2:** $L$ is the component $\overline{\mathbb{R}} \setminus [0, 1]$ of $\overline{\mathbb{R}} \setminus K$. Then $L$ consists of the large interval $J_\infty = \overline{\mathbb{R}} \setminus (-\sigma, 1 + \sigma)$ and Whitney intervals $\{J_{\pm n}\}_{n \geq N+3}$ in $[-\sigma, 0) \cup (1, 1 + \sigma]$, where $\sigma = (2 \cdot 3^{N+2})^{-1}$.

For the interval $J_\infty = \overline{\mathbb{R}} \setminus (-\sigma, 1 + \sigma)$, define $E_{J_\infty} = \left[\frac{1}{9}, \frac{8}{9}\right]$. Join the endpoints of $J_\infty$ to the endpoints of $\left[\frac{1}{9}, \frac{8}{9}\right]$ by non-intersecting semicircles in the upper half plane, and let the leaf $\Lambda_{J_\infty}$ be the region in the upper half plane bounded by $J_\infty$, $E_{J_\infty}$, and the two semicircles. See Figure 5(a).

Figure 5(a). Non-standard palm leaf $\Lambda_{J_\infty}$ for $J = J_\infty$.

Figure 5(b). Standard palm leaves for $J_1 = J_{N+3}$, $J_2 = J_{N+4}$ in $[-\sigma, 0)$.

The Whitney intervals in $[-\sigma, 0)$ are enumerated from left to right as $J_{N+3}, J_{N+4}, \ldots,$ in order of decreasing size. See Figure 5(b). For $J_{N+3} = [-\sigma, -\sigma/3]$, let $E_{J_{N+3}} = [3^{-3}, 3^{-2}]$. Form the leaf $\Lambda_{J_{N+3}}$ by joining the endpoints of $J_{N+3}$ to those of $E_{J_{N+3}}$ by two non-intersecting semicircles which meet $\mathbb{R}$ at right angles. Then $|J_{N+3}| = \frac{2}{3}\sigma = 3^{-(N+3)} = (2 \cdot 3^N)^{-1}|E_{J_{N+3}}|$. So $\Lambda_{J_{N+3}}$ and $E_{J_{N+3}}$ are of the standard form.

Similarly, for $J_{N+k} = [-\sigma/3^{k-3}, -\sigma/3^{k-2}]$, $k \geq 3$, let $E_{J_{N+k}} = [3^{-k}, 3^{-k+1}]$ and define $\Lambda_{J_{N+k}}$ as usual as the region in the upper half plane bounded by $J_{N+k}$, $E_{J_{N+k}}$, and the two non-intersecting semicircles joining the endpoints of $J_{N+k}$ to those of $E_{J_{N+k}}$. Then $|E_{J_{N+k}}| = 2 \cdot 3^N |J_{N+k}|$ for $k \geq 3$; the $J_{N+k}$'s are standard intervals with standard $E_{J_{N+k}}$'s and $\Lambda_{J_{N+k}}$'s; the $E_{J_{N+k}}$'s cover $[0, \frac{1}{9}]$, and the $\Lambda_{J_{N+k}}$'s tile the half-disc bounded by $[-\sigma, \frac{1}{9}]$ and the semicircle in the upper half plane which joins $-\sigma$ to $\frac{1}{9}$. See Figure 5(b).

The Whitney intervals in $(1, 1 + \sigma]$ are enumerated from right to left in order of decreasing size, as $J_{-(N+3)}, J_{-(N+4)}, \ldots$ Define leaves $\Lambda_{J_{-(N+k)}}$ for $J_{-(N+k)}$, $k \geq 3$, by reflecting the leaves $\Lambda_{J_{N+k}}$ through the line $x = \frac{1}{2}$. The leaf $\Lambda_{J_{-(N+k)}}$ based at $J_{N+k}$ is the mirror image of the leaf $\Lambda_{J_{N+k}}$ based at $J_{N+k}$. All these leaves and their $E_{J_{N+k}}$'s are standard. The $E_{J_{N+k}}$, $k \geq 3$, cover $[\frac{8}{9}, 1]$.

In the case $L = \overline{\mathbb{R}} \setminus [0, 1]$ we have defined standard leaves $\Lambda_J$ for all the Whitney intervals $J$ in $L$, except for $J_\infty = \overline{\mathbb{R}} \setminus (-\sigma, 1 + \sigma)$, which has a non-standard leaf. For all intervals $J \subset L$, the endpoints of $E_J$ lie in the Cantor set $K$. The palm leaves have the same properties as those summarized at the end of Case 1.

For each component $L$ of $\overline{\mathbb{R}} \setminus K$, reflect the fan of palm leaves based at $L$ through the real axis, obtaining a fan of palm leaves, also based at $L$, which tiles the lower half plane.

## 10. Distortion estimates for Whitney intervals

Let $\widehat{I}$ be a Whitney interval in the upper orthocircular boundary arc $A_{n,j}$ of a half fundamental domain $\Omega_{n,j}$. Let $I = P(\widehat{I})$ be the corresponding grid interval, and let $\widetilde{I}$ be the segment of $\partial\Omega_{n,j}$, below $A_{n,j}$, with the same endpoints as $I$. The purpose of this section is to show that the Euclidean lengths of $\widehat{I}$, $\widetilde{I}$, and $I$ are comparable to each other, with constants which are uniform for all $\widehat{I}$.

We prove a preliminary lemma.

**Lemma 10.1.** *The orthocircular arcs in $\bigcup_n \mathcal{A}_n$ are uniformly hyperbolically separated.*

**Proof.** It is sufficient to prove that the components of $\overline{\mathbb{R}} \setminus K$ are uniformly hyperbolically separated, since these components lift via $\pi^{-1}$ to $\bigcup_n \mathcal{A}_n$, and the conformal map $\pi^{-1}$ is a hyperbolic isometry.

Recall that the hyperbolic metric on $\Omega = \overline{\mathbf{C}} \setminus K$, given by $\lambda_\Omega(w)|dw|$, satisfies

$$\frac{c_\Omega}{\operatorname{dist}(w, K)} \leq \lambda_\Omega(w) \leq \frac{2}{\operatorname{dist}(w, K)} \tag{10.1}$$

where $c_\Omega > 0$ depends only on the uniformly perfect constant of the Cantor set $K$.

Let $L$ and $L'$ be components of $\overline{\mathbb{R}} \setminus K$, with $|L| \leq |L'|$. Let $\gamma$ be an arc joining $L$ to $L'$. Then the Euclidean length of $\gamma$ is at least $|L|$, since between $L$ and $L'$ there is a closed construction interval of $K$ of length at least $|L|$.

If each point in $\gamma$ is within distance $2|L|$ of $K$, then

$$
\begin{aligned}
\ell_{\mathrm{hyp}}(\gamma) &= \int_{\gamma} \lambda_{\Omega}(w)|dw| \\
&\geq c_{\Omega} \int_{\gamma} \frac{|dw|}{\mathrm{dist}\,(w, K)} \\
&\geq c_{\Omega} \int_{\gamma} \frac{|dw|}{2|L|} \\
&\geq c_{\Omega}/2.
\end{aligned}
\tag{10.2}
$$

If some point in $\gamma$ is not within distance $2|L|$ of $K$, then there is a segment $\gamma'$ of $\gamma$ of length $|L|$, which has one endpoint at distance $2|L|$ from $K$, and which stays within distance $2|L|$ of $K$. Then

$$
\begin{aligned}
\ell_{\mathrm{hyp}}(\gamma) &\geq \ell_{\mathrm{hyp}}(\gamma') \\
&= \int_{\gamma'} \lambda_{\Omega}(w)\,|dw| \\
&\geq c_{\Omega} \int_{\gamma'} \frac{|dw|}{\mathrm{dist}\,(w, K)} \\
&\geq \frac{c_{\Omega}}{2|L|}\,|L| \\
&= \frac{c_{\Omega}}{2}.
\end{aligned}
\tag{10.3}
$$

Therefore the hyperbolic distance between any two orthocircular arcs in $\bigcup_n \mathcal{A}_n$ is at least $c_{\Omega}/2$. $\qquad\square$

**Lemma 10.2.** *Let $I \in \mathcal{H}$ be any grid interval. Then $|\widehat{I}| \overset{c}{\sim} |\widetilde{I}| \overset{\frac{\pi}{2}}{\sim} |I|$, and the constant $c > 0$ is independent of $I$.*

We split the proof into several sublemmas, showing that $|\widehat{I}| \leq c\,|I|$ for all standard and non-standard grid intervals, and then that $|\widehat{I}| \geq c\,|I|$ for all standard and non-standard grid intervals. As usual, $c$ denotes constants which may change from line to line.

Note that if $A_{n,j}$ is any orthocircular arc, then its length is comparable with constant $\pi/2$ to the length of the arc of $\partial\mathbb{D}$ below $A_{n,j}$ which has the same endpoints as $A_{n,j}$. For each $I \in \mathcal{H}$, the endpoints of $\widetilde{I}$ are the endpoints of $I$, and they lie in $\partial\mathbb{D}$. Therefore $\widetilde{I}$ consists of a subset of $\partial\mathbb{D}$ together with a union of whole orthocircular arcs, and so $|\widetilde{I}| \sim |I|$ with constant $\pi/2$.

**Lemma 10.3.** $|\widehat{I}| \leq c\,|I|$ *for all standard intervals* $I \in \mathcal{H}$.

**Proof.** Let $A_{n,j}$ be the orthocircular arc containing $\widehat{I}$, and let $\Omega_{n,j}$ be the half fundamental domain below $A_{n,j}$. Without loss of generality, assume that $\pi(\Omega_{n,j})$ is the upper half plane. Let $\mathcal{F}$ be the fundamental domain which consists of $\Omega_{n,j}$, the arc $A_{n,j}$, and the half fundamental domain immediately above $A_{n,j}$. Let $L = \pi(A_{n,j})$. Then $\pi(\mathcal{F})$ is the union of the open upper and lower half planes together with the component $L$ of $\overline{\mathbb{R}} \setminus K$; its boundary is $\partial \pi(\mathcal{F}) = \overline{\mathbb{R}} \setminus L$.

Let $J = \pi(\widehat{I})$; $J$ is a standard Whitney interval in $\overline{\mathbb{R}} \setminus K$. As usual, let $E_J$ be the tip of the standard leaf based at $J$; $E_J = \pi(\widetilde{I})$. See Figure 6.

Figure 6. Standard grid interval $I$; definition of $z_J$.

Let $w_J$ be the midpoint of $J$ and let $w_I$ be the point in $\widehat{I}$ such that $\pi(w_I) = w_J$. Let $\alpha \leq |J|/2$ be a small number, which will be fixed later. Let $z_J$ be the point in the upper half plane, directly above the midpoint $w_J$ of $J$, such that $|z_J - w_J| = \alpha |J|$. Let $z_I$ be the point in $\Omega_{n,j}$ such that $\pi(z_I) = z_J$.

The hyperbolic distance from $z_I$ to $w_I$ is bounded away from zero and infinity, uniformly for all $\widehat{I}$. The Euclidean ball of radius $|J|/2$ centred at $w_J$, which contains $z_J$ and $J$, lies in the set $\{z \mid \mathrm{dist}\,(J, K) \leq \mathrm{dist}\,(z, K) \leq 3\,\mathrm{dist}\,(J, K)\}$, since $|J| = 2\,\mathrm{dist}\,(J, K)$. It

follows that

$$
\begin{aligned}
d_{\mathrm{hyp}}(z_I, w_I) &= d_{\mathrm{hyp}}(z_J, w_J) \\
&\leq |z_J - w_J| \, \frac{2}{\mathrm{dist}\,(J, K)} \\
&= \alpha \, |J| \, \frac{2}{|J|/2} \\
&= 4\alpha.
\end{aligned}
\tag{10.4}
$$

Also,

$$
\begin{aligned}
d_{\mathrm{hyp}}(z_J, w_J) &\geq |z_J - w_J| \, \frac{c_\Omega}{3 \, \mathrm{dist}\,(J, K)} \\
&= \frac{2}{3} \, c_\Omega \alpha.
\end{aligned}
\tag{10.5}
$$

So $\frac{2}{3} c_\Omega \alpha \leq d_{\mathrm{hyp}}(z_I, w_I) \leq 4\alpha$ for all $I$. This implies that the Euclidean distance from $z_I$ to $w_I$ is less than $c(\alpha)\,(1 - |w_I|)$, where $c(\alpha)$ depends only on $\alpha$ and decreases to zero with $\alpha$.

Since $z_I$ and $w_I$ are close in the hyperbolic metric, harmonic measure on $\partial \mathcal{F}$ with basepoint $z_I$ is close to harmonic measure on $\partial \mathcal{F}$ with basepoint $w_I$. Let $u(z) = \omega(z, E, \mathcal{F})$, where E is any Borel subset of $\partial \mathcal{F}$. The function $u(z)$ is harmonic in $\mathcal{F}$. The hyperbolic distance from $w_I$ to $\partial \mathcal{F}$ is at least $c_\Omega/2$, by Lemma 10.1. The hyperbolic ball of radius $c_\Omega/2$ centred at $w_I$ contains a Euclidean ball of radius $c'_\Omega\,(1 - |w_I|)$, centred at $w_I$, where $c'_\Omega$ depends only on $c_\Omega$. The Euclidean distance from $z_I$ to $w_I$ decreases to zero with the hyperbolic distance from $z_I$ to $w_I$. Choose the number $\alpha$, where $|z_I - w_I| = \alpha\,|J|$ and $|z_I - w_I| \leq c(\alpha)\,(1 - |w_I|)$, small enough that the Euclidean distance from $z_I$ to $w_I$ is less than $(c'_\Omega/3)\,(1 - |w_I|)$. By Harnack's inequality,

$$
\frac{1}{2} u(w_I) \leq u(z_I) \leq 2u(w_I).
\tag{10.6}
$$

In particular, setting $E = \widetilde{I}$, we obtain $\omega(z_I, \widetilde{I}, \mathcal{F}) \overset{2}{\sim} \omega(w_I, \widetilde{I}, \mathcal{F})$ for all $I$.

After these preliminaries, we now use estimates on harmonic measure based at $z_I$ and at $w_I$ to show that $|\widetilde{I}| \leq c\,(1 - |w_I|)$. This is sufficient to establish Lemma 10.3, since $|\widetilde{I}| \sim 1 - |w_I|$ with a constant independent of $I$ (by Lemma 8.1).

For standard intervals $J$, the shapes of the leaves $\Lambda_J$ and the position of $z_J$ within $\Lambda_J$ are all identical, up to dilations and reflections. Specifically, $|E_J| = 2 \cdot 3^N |J|$; the distance from $J$ to $E_J$ is also a fixed multiple of $|J|$; and $z_J$ is always at height $\alpha|J|$ above the midpoint of $J$. Therefore the harmonic measure of $E_J$ in the upper half plane $\mathbf{U}$, taken with respect to the basepoint $z_J$, is constant for all standard $J$: $\omega(z_J, E_J, \mathbf{U}) = c > 0$.

Let $\mathcal{F}_{\widetilde{I}}$ be the unit disc without the fundamental domains below $\widetilde{I}$. Let $\gamma$ be the orthocircular arc joining the endpoints of $\widetilde{I}$, and let $\mathcal{F}_\gamma$ be the disc without the "bite" below $\gamma$. So $\Omega_{n,j} \subset \mathcal{F} \subset \mathcal{F}_\gamma \subset \mathcal{F}_{\widetilde{I}}$. Then

$$
\begin{aligned}
c &= \omega(z_J, E_J, \mathbf{U}) \\
&= \omega(z_I, \widetilde{I}, \Omega_{n,j}) \\
&\leq \omega(z_I, \widetilde{I}, \mathcal{F}) \\
&\leq 2\,\omega(w_I, \widetilde{I}, \mathcal{F}) \\
&\leq 2\,\omega(w_I, \widetilde{I}, \mathcal{F}_{\widetilde{I}}) \\
&= 2\left[1 - \omega(w_I, \partial\mathcal{F}_{\widetilde{I}} \setminus \widetilde{I}, \mathcal{F}_{\widetilde{I}})\right] \\
&\leq 2\left[1 - \omega(w_I, \partial\mathcal{F}_{\widetilde{I}} \setminus \widetilde{I}, \mathcal{F}_\gamma)\right] \\
&= 2\left[1 - \omega(w_I, \partial\mathcal{F}_\gamma \setminus \gamma, \mathcal{F}_\gamma)\right] \\
&= 2\,\omega(w_I, \gamma, \mathcal{F}_\gamma).
\end{aligned}
\tag{10.7}
$$

Figure 7. Möbius transformations $\tau, \sigma : \mathbb{D} \to \mathbb{D}$.

We have used the conformal invariance of harmonic measure; the observations made above; and the fact that harmonic measure $\omega(z, E, \Omega)$ is monotonic in the domain $\Omega$: if $z \in \Omega \subset \Omega'$ and $E \subset \partial\Omega \cap \partial\Omega'$, then $\omega(z, E, \Omega) \leq \omega(z, E, \Omega')$.

We assumed here that $w_I$ was not on or below the orthocircular arc $\gamma$ which joins the endpoints of $I$. Otherwise, $|I| \geq |\gamma|/2 \geq 1 - |w_I|$ and we are done.

Let $\tau : \mathbb{D} \to \mathbb{D}$ be the Möbius transformation $\tau(z) = (z - w_I)/(1 - \overline{w}_I z)$ which takes $w_I$ to 0. See Figure 7. Then

$$
\begin{aligned}
\omega(w_I, \gamma, \mathcal{F}_\gamma) &= \omega(0, \tau(\gamma), \tau(\mathcal{F}_\gamma)) \\
&= 2\,\omega(0, \tau(I), \mathbb{D}) \\
&= \pi^{-1}|\tau(I)|.
\end{aligned}
\tag{10.8}
$$

The second step can be justified by mapping $\tau(\gamma)$ to a diameter of $\mathbb{D}$ via a Möbius transformation $\sigma : \mathbb{D} \to \mathbb{D}$ (Figure 7). Consider the probability that a Brownian traveller from $\sigma(0)$ in $\mathbb{D}$ first hits $\partial\mathbb{D}$ somewhere on $\sigma \circ \tau(I)$. By symmetry, this is exactly half the probability that the traveller first hits the boundary of the half disc $\sigma \circ \tau(\mathcal{F}_\gamma)$ somewhere on the diameter $\sigma \circ \tau(\gamma)$. That is,

$$
\begin{aligned}
\omega(0, \tau(I), \mathbb{D}) &= \omega(\sigma(0), \sigma \circ \tau(I), \mathbb{D}) \\
&= 2^{-1}\, \omega(\sigma(0), \sigma \circ \tau(\gamma), \sigma \circ \tau(\mathcal{F}_\gamma)) \\
&= 2^{-1}\, \omega(0, \tau(\gamma), \tau(\mathcal{F}_\gamma)),
\end{aligned}
\tag{10.9}
$$

as required.

We have shown that $c \leq 2\,\omega(w_I, \gamma, \mathcal{F}_\gamma)| \leq 2\pi^{-1}|\tau(\gamma)|$. Now $\tau'(z) = (1 - |w_I|^2)/(1 - w_I z)^2$, and so

$$
\begin{aligned}
|\tau(I)| &\leq |I| \max_{z \in I} |\tau'(z)| \\
&\leq |I|\, \frac{2}{1 - |w_I|}.
\end{aligned}
\tag{10.10}
$$

It follows that $|I| \geq c\,(1 - |w_I|)$, where the constant $c$ is independent of $I$, and therefore $|I| \geq c\,|\widehat{I}|$, with $c$ independent of $I$. $\qquad\square$

**Lemma 10.4.** $|\widehat{I}| \leq c\,|I|$ *for all non-standard intervals* $I \in \mathcal{H}$.

**Proof.** We use the same proof as for Lemma 10.3 above. The key point is to choose a suitable point $z_J$ in the upper half plane $\mathbf{U}$ for each $J$.

Let $J$ be a non-standard Whitney interval. $J$ is of one of two types: either $J$ is the central Whitney interval $J_c$ in an open component $L$ of $[0,1] \setminus K$, or $J = J_\infty$. In the first case, define $z_J$ as in the proof of Lemma 10.3, with the same $\alpha$. See Figure 8(a). The leaves $\Lambda_J$ for these intervals are identical up to dilations, so the harmonic measure of $E_J$ in $\mathbf{U}$, taken from $z_J$, is uniformly bounded below for all these $J$. (The bound is not necessarily the same as that in Lemma 10.3.) Also, $z_J$ is close enough to $w_J$ in the hyperbolic metric that the other estimates on harmonic measure hold with the same constants as in Lemma 10.3.

Figure 8(a). Definition of $z_J$ for non-standard $J = J_c$.

For $J = J_\infty$, fix a point $w_J \in J_\infty$ which is far from $[0, 1]$. See Figure 8(b).

Figure 8(b). Definition of $z_J$ for non-standard $J = J_\infty$.

Let $z_J$ be a point in $\mathbf{U}$, directly above $w_J$, which is close to $w_J$ in the hyperbolic metric. Then the harmonic measure of $E_J$ in $\mathbf{U}$, taken from $z_J$, is a positive constant, not necessarily the same as the lower bound in the earlier cases. Choose $z_J$ sufficiently close to $w_J$ that the other harmonic measure estimates hold with the same constants as in the proof of Lemma 10.3. $\square$

**Lemma 10.5.** $|\widehat{I}| \geq c\,|I|$ for all standard intervals $I \in \mathcal{H}$.

**Proof.** We use the same notation as in the proof of Lemma 10.3. In particular, $w_I$ is the point in $\widehat{I}$ such that $w_J = \pi(w_I)$ is the midpoint of the Whitney interval $J = \pi(\widehat{I})$. Without loss of generality assume that $w_I$ is in the left half of the orthocircular arc $A_{n,j}$ containing $\widehat{I}$. Let $a$ and $b$ be the left and right endpoints, respectively, of $I$. See Figure 9.

It is enough to show that $|b - w_I| \leq c\,(1 - |w_I|)$ with a constant $c$ independent of $I$. For $I$ is trapped between $b$ and the left endpoint $a'$ of $A_{n,j}$, which satisfies $|a' - w_I| \leq 2\,(1 - |w_I|)$. So then

$$\begin{aligned}
|I| &\leq |a' - b| \\
&\leq |a' - w_I| + |w_I - b| \\
&\leq c\,(1 - |w_I|) \\
&\leq c\,|\widehat{I}|
\end{aligned} \tag{10.11}$$

(using Lemma 8.1), with a constant $c$ independent of $I$.

Let $B$ be the shorter segment of $\partial\mathcal{F}$ between the right endpoint $b$ of $I$ and the right endpoint $b'$ of $A_{n,j}$. See Figure 9. The following calculation is justified below:

$$\begin{aligned}
|w_I - b| &\leq c\,\mathrm{dist}\,(w_I, B) \\
&\leq c\,\frac{\mathrm{dist}\,(w_I, \partial\mathcal{F})}{\omega(w_I, B, \mathcal{F})^2} \\
&\leq c\,\mathrm{dist}\,(w_I, \partial\mathcal{F}) \\
&\leq c\,(1 - |w_I|),
\end{aligned} \tag{10.12}$$

where dist denotes Euclidean distance, and $c$ denotes constants independent of $I$.

Figure 9. $\widetilde{I}$ is near $\widehat{I}$; definition of $B$.

The Euclidean distance from $w_I$ to $B$ is $|w_I - b|$, unless there is a point in an ortho-circular arc in $B$ which lies closer to $w_I$ than $b$ does. If $|b - w_I| \leq 10\,(1 - |w_I|)$ the lemma is proved. If not, $B$ is far to the right of $w_I$ and the distance from $w_I$ to any point in $B$ is at least $\frac{1}{2}|w_I - b|$, justifying the first line of (10.12).

By Beurling's Lemma [A], there is a constant $C > 0$ such that for each point $z_0$ in $\mathcal{F}$, and for all $M > 0$,

$$\omega(z_0, \{z \mid |z - z_0| \geq M \operatorname{dist}(z_0, \partial\mathcal{F})\}, \mathcal{F}) \leq CM^{-\frac{1}{2}}. \tag{10.13}$$

Setting $z_0 = w_I$ and $M = \operatorname{dist}(w_I, B)/\operatorname{dist}(w_I, \partial\mathcal{F})$, we find

$$\omega(w_I, B, \mathcal{F}) \leq \omega(w_I, \{z \mid |z - w_I| \geq \operatorname{dist}(z_0, B)\}, \mathcal{F})$$
$$\leq C \left[\frac{\operatorname{dist}(w_I, \partial\mathcal{F})}{\operatorname{dist}(w_I, B)}\right]^{\frac{1}{2}}, \tag{10.14}$$

and so

$$\operatorname{dist}(w_I, B) \leq C \frac{\operatorname{dist}(w_I, \partial\mathcal{F})}{\omega(w_I, B, \mathcal{F})^2}, \tag{10.15}$$

which is the second line of (10.12).

The harmonic measure of $B$ in $\mathcal{F}$, measured from $w_I$, is bounded below by a positive constant independent of $I$. With $z_I$ as in the proof of Lemma 10.3,

$$\begin{aligned}
\omega(w_I, B, \mathcal{F}) &\geq 2^{-1}\,\omega(z_I, B, \mathcal{F}) \\
&= 2^{-1}\,\omega(z_J, \pi(B), \pi(\mathcal{F})) \\
&\geq 2^{-1}\,\omega(z_J, \pi(B), \mathbf{U}),
\end{aligned} \tag{10.16}$$

by (10.6) and since $\pi(\mathcal{F})$ contains $\mathbf{U}$. $\pi(B)$ is one of the two components of $\overline{\mathbb{R}} \setminus (L \cup E_J)$, where $L$ is the open component of $\overline{\mathbb{R}} \setminus K$ which contains $J$. Again, since the picture of $\mathbf{U}$ with $J$, $z_J$, and $E_J$ is invariant up to dilations and reflections for all standard $J$, the harmonic measure of $\pi(B)$ in $\mathbf{U}$, taken from $z_J$, is bounded below by a constant $c > 0$ independent of $J$, for all standard $J$. This justifies the third line of (10.12). The last line of (10.12) is true because $\widehat{I}$ is a Whitney interval in $A_{n,j}$. $\qquad\square$

**Lemma 10.6.** $|\widehat{I}| \geq c\,|I|$ for all non-standard intervals $I \in \mathcal{H}$.

**Proof.** We use the same method as for Lemma 10.5. Define $z_J$ for each non-standard $J$ as in the proof of Lemma 10.4. The interval $\pi(B)$ is one of the two components of $\overline{\mathbb{R}} \setminus (L \cup E_J)$. See Figures 8(a) and 8(b). By the usual arguments we may conclude that the harmonic measure of $B$ in $\mathcal{F}$, taken from $w_I$, is bounded below by a positive constant independent of $I$. The rest of the proof applies without change. $\qquad\square$

This completes the proof of Lemma 10.2: $|\widehat{I}| \sim |\widetilde{I}| \sim |I|$ for all grid intervals $I \in \mathcal{H}$, with constants which are independent of $I$.

## 11. Definition of the density functions $F_n \cdots F_1$

In this section we define the functions $F_n$ whose products $F_n \cdots F_1$ give the densities of the measures $\mu_n$. Here we define the $F_n$ on all standard grid intervals $I \in \mathcal{H}_{n-1}$, $n \geq 1$, up to the values of certain parameters, $Q$ and $\varepsilon$, which will be fixed later. We show that these functions $F_n$ are $(\delta, \eta)$-suitable for each standard $I \in \mathcal{H}_{n-1}$, with constants $\delta$ and $\eta$ independent of $I$ and $n$. For non-standard intervals $I \in \mathcal{H}_{n-1}$, we define the functions $F_n$ and prove that they are $(\delta, \eta)$-suitable in Section 14.

We first prove a lemma giving a rough estimate of the distortion caused by the composition $P \circ \pi^{-1}$ of a branch of the inverse $\pi^{-1}$ of the covering map with the almost radial projection $P$ of the disc onto the circle.

**Lemma 11.1.** *For each standard Whitney interval $J$ in $\overline{\mathbb{R}} \setminus K$, let $S_J$ be a segment of $E_J$, and let $A_1$ and $A_2$ be the components (possibly empty) of $E_J \setminus S_J$. Suppose that $|S_J|/|E_J| = c_1$, $|A_1|/|E_J| = c_2$, and $|A_2|/|E_J| = 1 - c_1 - c_2$, where $c_1$ and $c_2$ are constants independent of $J$. Suppose also that $S_J$ is a union of whole Whitney intervals, possibly together with a subset of $K$. Let $I$ be any standard grid interval such that $\pi(\widehat{I}) = J$, and let $\widehat{R_I}$ be the segment of $\widehat{I}$ such that $\pi(\widehat{R_I}) = S_J$. Let $R_I = P(\widehat{R_I}) = P \circ \pi^{-1}(S_J)$. Then $|R_I| \geq c\,|I|$, where $c$ is a constant independent of $I$ and $J$.*

**Proof.** Without loss of generality, assume that $\pi$ maps the half fundamental domain below the orthocircular arc containing $\widehat{I}$ to the upper half plane. As in the proof of Lemma 10.3, let $w_J$ be the midpoint of $J$, and let $z_J$ be the point in the upper half plane, directly above $w_J$, such that $|z_J - w_J| = \alpha|J|$, where $\alpha$ is a small fixed number independent of $J$. Let $w_I$ be the point in $\widehat{I}$ such that $\pi(w_I) = w_J$. See Figure 6 in Section 10.

The harmonic measure of $S_J$ in the upper half plane, taken from $z_J$, is the same for all $J$, since $S_J$ is always in the same place in $E_J$. The proof of Lemma 10.3, applied to $S_J$ instead of to $E_J$, shows that $|R_I| \geq c\,(1 - |w_I|)$, where $c$ is a constant depending on $c_1$ and $c_2$. By Lemmas 8.1 and 10.2, $1 - |w_I|$ is comparable to $|I|$ with a constant independent of $I$ and $J$, which proves the lemma. $\qquad\square$

We have shown that for all standard intervals $I \in \mathcal{H}$, if $S_J$ is a fixed segment of $E_J$, where $J = \pi(\widehat{I})$, then $|R_I|/|I| = |P \circ \pi^{-1}(S_J)|/|P \circ \pi^{-1}(E_J)|$ is uniformly bounded away from zero. The same result holds for non-standard intervals $I$, if $S_J$ is defined so that its length is a fixed multiple of $|J|$ and so that it is always in the same position in $E_J$.

On each interval $I$ in the grid $\mathcal{H}$ of subintervals of the unit circle we define a function $F$ which is $(\delta, \eta)$-suitable for $I$, where the constants $\delta$ and $\eta$ are independent of $I$. Recall from Section 6 the definition of $(\delta, \eta)$-suitable: $F$ has mean value one on $I$; $0 < \delta \leq F(x) \leq 1/\delta$ on $I$; and $F \equiv 1$ on subintervals $I_l$ and $I_r$ of $I$, at the left and right ends respectively of $I$, such that $|I_l| \geq \eta|I|$ and $|I_r| \geq \eta|I|$.

The idea is as follows. Divide each $E_J$ into five segments $S_1, \ldots, S_5$ whose lengths are prescribed fractions of $|E_J|$. Pull these back to $\widetilde{I}$ via $\pi^{-1}$, then project them via $P$ to subintervals $R_j = P \circ \pi^{-1}(S_j)$, $1 \leq j \leq 5$, of $I$. In Lemma 11.1 we gave bounds on the distortion in length caused by $P \circ \pi^{-1}$. Each ratio $|R_j|/|I|$ is uniformly bounded away from zero and infinity for all $I$. Define $F$ to be constant on each $R_j$, in such a way that $F$ is large on $R_3$, identically equal to one on $R_1$ and $R_5$, and small enough on $R_2 \cup R_4$ to ensure that $F$ has mean value one on $I$. See Figure 10. We show that $F$ is uniformly bounded away from zero and infinity for all $I$.

Figure 10. Graph of $F$; $\alpha = (1-\varepsilon)\,|I|/|R_3|$.

Now we define the segments $S_j$ and $R_j$ for a standard grid interval $I$. Let $\widehat{I}$ be the Whitney interval in $\bigcup_n \mathcal{A}_n$ such that $I = P(\widehat{I})$, and let $J = \pi(\widehat{I})$. Without loss of generality we describe the case when $E_J$ is to the right of $J$. Divide $E_J$ into five segments $S_j$, $1 \le j \le 5$, numbered from left to right, as follows. The left half of $E_J$ is an open component of $\overline{\mathrm{IR}} \setminus K$ of length $3^N|J|$. Let $S_3$ be the central Whitney interval $J_c$ in this open interval. The length of $S_3$ is $|S_3| = \frac{1}{2}(1-3^{-N})|E_J|$.

Let $Q$ be a large integer; its value will be fixed below. Let $S_1$ be the interval with the same left endpoint as $E_J$ and with length $|S_1| = \frac{3}{2}3^{-Q}|J|$. $S_1$ is a union of whole Whitney intervals. Let $S_5$ be the interval with the same right endpoint as $E_J$ and with length $|S_5| = 3^{-Q}|J|$. $S_5$ is a closed construction interval of the Cantor set $K$. Let $S_2$ be the interval between $S_1$ and $S_3$, and let $S_4$ be the interval between $S_3$ and $S_5$.

Now pull back these intervals $S_j$ to $\widetilde{I}$ via the appropriate branch of $\pi^{-1}$. Then project them to $I$ via $P$, obtaining five intervals $R_j = P \circ \pi^{-1}(S_j)$, $1 \le j \le 5$, whose union is $I$. By Lemma 11.1, there are positive constants $c_j$, $1 \le j \le 5$, such that $|R_j| \ge c_j\,|I|$ for all standard intervals $I$. These constants are less than one; they depend on $N$ and $Q$ but are independent of $J$ and $I$.

Let $\varepsilon$ be a positive number, small enough that $(1-\varepsilon)|I|/|R_3| > 1$; its exact value will be fixed below. Define a step function $F$ on $I$ by

$$F(x) = \begin{cases} 1, & x \in R_1 \cup R_5; \\ (1-\varepsilon)\frac{|I|}{|R_3|}, & x \in R_3; \\ \delta, & x \in R_2 \cup R_4; \end{cases} \tag{11.1}$$

where $\delta = (\varepsilon|I| - |R_1 \cup R_5|)/|R_2 \cup R_4|$ is chosen so that $\frac{1}{|I|}\int_I F = 1$. See Figure 10.

This function $F$ takes its maximum on $R_3$, where $F(x) \equiv (1-\varepsilon)|I|/|R_3| \leq 1/c_3$; this upper bound is independent of $I$. The subintervals $R_1$ and $R_5$ at each end of $I$, on which $F \equiv 1$, each have length at least $\eta|I|$, where $\eta = \min(c_1, c_5)$ is independent of $I$.

If $d\mu_n = F_n \cdots F_1\, dx$, and each $F_m$ is defined on each $I \in \mathcal{H}_{m-1}$, $m \geq 1$, according to (11.1), then most of the mass $\mu_{n-1}(I)$ is concentrated onto the subinterval $R_3$ of $I$:

$$
\begin{aligned}
\mu_n(R_3) &= \int_{R_3} F_n(x)\, d\mu_{n-1}(x) \\
&= (1-\varepsilon)\frac{|I|}{|R_3|}\, \mu_{n-1}(R_3) \\
&= (1-\varepsilon)\frac{|I|}{|R_3|} \cdot \frac{|R_3|}{|I|}\, \mu_{n-1}(I) \\
&= (1-\varepsilon)\, \mu_{n-1}(I).
\end{aligned}
\tag{11.2}
$$

The third step is valid because $\mu_{n-1}$ is a constant multiple of Lebesgue measure on $I$.

To show that the function $F$ defined in (11.1) is $(\delta_0, \eta)$-suitable for $I$, with constants $\delta_0$ and $\eta$ independent of $I$, it only remains to give a uniform lower bound for the value $\delta$ of $F$ on $R_2 \cup R_4$. Observe that

$$
\delta = \frac{\varepsilon|I| - |R_1 \cup R_5|}{|R_2 \cup R_4|} = \left(\varepsilon - \frac{|R_1 \cup R_5|}{|I|}\right) \frac{|I|}{|R_2 \cup R_4|} \geq \varepsilon - \frac{|R_1 \cup R_5|}{|I|}.
$$

We would like to ensure that $|R_1 \cup R_5|/|I| \leq \varepsilon/2$, say, by choosing $|S_1 \cup S_5|/|E_J|$ sufficiently small (that is, by choosing $Q$ sufficiently large). Unfortunately, an application of Lemma 11.1 to $S_J = E_J \setminus (S_1 \cup S_5)$ does not guarantee that this can be done, since the lemma does not give a good estimate on the size of the constant $c$, nor on its behaviour as $c_1 \to 0$. However, the sharper estimate in Lemma 13.1 below, based on the $A_\infty$-equivalence of harmonic measure and arclength on the boundary of a chord-arc domain, does imply that $|R_1 \cup R_5|/|I|$ goes to zero with $|S_1 \cup S_5|/|E_J|$. Therefore, we may choose a large $Q$ so that $|S_1 \cup S_5|/|E_J| = \frac{5}{2} \cdot \frac{1}{2} \cdot 3^{-Q-N}$ is sufficiently small that $|R_1 \cup R_5|/|I| \leq \varepsilon/2$. Then $\delta \geq \varepsilon/2$.

We have shown that the function $F$ defined in (11.1) is $(\delta_0, \eta)$-suitable for $I$, where the constants $\delta_0$ and $\eta$ depend on $\varepsilon$, $N$, and $Q$, but are independent of $I$, for all standard grid intervals $I \in \mathcal{H}$.

In Section 14 we define the function $F$ for non-standard intervals $I \in \mathcal{H}$, and show that it is $(\delta_0, \eta)$-suitable for those $I$.

## 12. Definition of auxiliary functions $X_i$

Let $I_0$ be an interval from the $l^{\text{th}}$ layer $\mathcal{H}_l$ of the grid, such that $\pi(\widehat{I}_0)$ is not $J_\infty$. For each point $x \in I_0$ and for each $i \geq 0$, let $I_i(x)$ be the unique interval in $\mathcal{H}_{l+i}$ which contains $x$. The intervals $I_i(x)$ are nested: $I_0(x) \supset I_1(x) \supset I_2(x) \supset \cdots \supset I_i(x) \supset \cdots \ni x$. Let $\widehat{I}_i(x)$ be the Whitney interval in the disc such that $P(\widehat{I}_i(x)) = I_i(x)$. Let $J_i(x) = \pi(\widehat{I}_i(x))$. We have associated to the point $x \in I_0 \subset \partial \mathbb{D}$ a sequence $\{J_i(x)\}_{i \geq 0}$ of Whitney intervals in $\overline{\mathbb{R}} \setminus K$.

We now define auxiliary functions $X_i(x)$ which keep track of the lengths $|J_i(x)|$ of these intervals. Make the convention that

$$|J_c| = |L|/3 \tag{12.1}$$

when $J_c$ is the non-standard central Whitney interval in a component $L$ of $[0,1] \setminus K$. Let

$$X_1(x) = \log_3 \left[ \frac{|J_1(x)|}{|J_0|} \right] \tag{12.2}$$

for $x \in I_0$. For $i \geq 2$, let

$$X_i(x) = \begin{cases} \log_3 \left[ \frac{|J_i(x)|}{|J_{i-1}(x)|} \right], & \text{if } J_1(x), \ldots, J_{i-1}(x) \neq J_\infty; \\ 1, & \text{otherwise,} \end{cases} \tag{12.3}$$

for $x \in I_0$ and $i \geq 1$. Let

$$S_k(x) = \sum_{i=1}^{k} X_i(x), \tag{12.4}$$

for $k \geq 1$. If none of $J_1(x), \ldots, J_{k-1}(x)$ is the large interval $J_\infty = \overline{\mathbb{R}} \setminus (-\sigma, 1+\sigma)$, then $S_k(x) = \log_3\big(|J_k(x)|/|J_0|\big)$. Note that $X_i$ and $S_k$ are integer-valued.

Let $V(I_0)$ denote the collection of those grid intervals $I$ in $I_0$ which satisfy $\pi(\widehat{I}) = J_\infty$, and which are maximal with respect to this property. In other words, $V(I_0)$ is the collection of the maximal grid intervals $I$ in $I_0$ which are images of $J_\infty$ under branches of $P \circ \pi^{-1}$.

Observe that $V(I_0) = \{x \in I_0 \mid S_k(x) \to +\infty\}$. For if $x$ is contained in some interval $I = I_N(x)$, say, which is in $V(I_0)$, then $X_i(x) = 1$ for all $i \geq N$, and so $S_k(x)$ tends to infinity. And if $x$ is not contained in any interval $I$ in $V(I_0)$, then no $J_i(x)$ is $J_\infty$, and so $S_k(x)$ remains bounded above, by $\log_3\big[(9|J_0|)^{-1}\big]$, for all $k$. (Recall that $J_0 = \pi(\widehat{I}_0)$ is fixed in this discussion, and that the largest Whitney interval other than $J_\infty$ is the $J_c$ in $\left(\frac{1}{3}, \frac{2}{3}\right)$, which has length $1/9$ by our convention.)

In the remainder of the paper we show that for appropriate choices of the parameters in the definitions of the functions $F_n$, the measures $\mu_n$ given by $d\mu_n = F_n \cdots F_1 \, dx$ converge to a measure $\mu$ whose restriction to $I_0$ is supported on $V(I_0)$. In Sections 13 and 14 we show that for all grid intervals $I \subset I_0$, the mean of $X_i$ on $I$ (where $i$ is such that $J_{i-1} = \pi(\widehat{I})$) with respect to $\mu$ is uniformly large, and the second moment is uniformly small. In Section 15 we conclude that $S_k \to +\infty$ a.e. $(d\mu)$ on $I_0$; that is, $\mu(V(I_0)) = \mu(I_0)$.

Finally, also in Section 15, we observe that this implies the existence of a doubling measure $\mu$ supported on the set $S$ of points which lie in infinitely many grid intervals corresponding to $J_\infty$. Heuristically, the grid intervals corresponding to $J_\infty$ are near preimages $\widehat{I_\infty}$ of $J_\infty$ under $\pi$. Each preimage $\widehat{I_\infty}$ contains an orbit point $g(0)$, $g \in G$, since the covering map $\pi$ is normalized so that the orbit of 0 is $\pi^{-1}(\infty)$, and $\infty \in J_\infty$. We prove that $S$ is contained in the conical limit set of $G$, which completes the proof of Theorem 1.2.

## 13. Estimates $EX_i \geq c_1$ and $EX_i^2 \leq c_2$ for standard intervals

The main result of this section (Lemma 13.4) is that on each standard grid interval $I \in \mathcal{H}$, the function $X_i$ satisfies $\mu(I)^{-1} \int_I X_i \, d\mu \geq c_1$ and $\mu(I)^{-1} \int_I X_i^2 \, d\mu \leq c_2$, where $c_1$ and $c_2$ are positive constants independent of $I$. Here $X_i = \log_3\big(|J_i(x)|/|J_{i-1}(x)|\big)$ is the auxiliary function, defined in Section 12, for a particle which has not yet reached $J_\infty$ and which makes its $i^{\text{th}}$ jump from $J = \pi(\widehat{I})$. We begin by establishing an estimate, sharper than that in Section 11, on the distortion in length caused by the map $P \circ \pi^{-1}$.

As usual, for a standard grid interval $I$ with $I = P(\widehat{I})$, let $J = \pi(\widehat{I})$, and let $E_J = \pi(\widetilde{I})$ be the tip of the leaf $\Lambda_J$ based at $J$. See Figure 11.

Figure 11. $B$ is in $I$; $\pi(\widehat{B})$ is in $E_J$.

**Lemma 13.1.** *Let $I \in \mathcal{H}_{n-1}$, $n \geq 1$, be a standard grid interval. Let $B \subset I$ be a union of grid intervals from $\mathcal{H}_n$. Let $\widehat{B}$ be the subset of $\widetilde{I}$ such that $P(\widehat{B}) = B$. There are positive constants $c$ and $\alpha$ such that*

$$\frac{|B|}{|I|} \leq c \left[ \frac{|\pi(\widehat{B})|}{|E_J|} \right]^\alpha . \tag{13.1}$$

*The constants $c$ and $\alpha$ are independent of $B$, $I$, and $n$, but they depend on the large number $N$ such that $|E_J| = 2 \cdot 3^N |J|$ for standard intervals.*

To prove this lemma, we use Lemmas 13.2 and 13.3 below and the following remarks.

A Jordan curve $\Gamma$ is said to be *chord-arc* if for all points $x$ and $y$ on $\Gamma$ the Euclidean length $|x - y|$ of the chord between $x$ and $y$ is comparable, with a uniform constant, to the Euclidean arclength of the shorter arc of $\Gamma$ between $x$ and $y$. In other words, there is a constant $c > 0$ such that

$$|x - y| \leq \ell_\Gamma(x, y) \leq c \, |x - y| \tag{13.2}$$

for all $x$, $y$ on $\Gamma$. A domain $D$ is a *chord-arc domain* if its boundary is a chord-arc curve.

The leaves $\Lambda_J$ are chord-arc domains with chord-arc constant independent of $J$.

If $D$ is a bounded chord-arc domain, then arclength on its boundary $\partial D$ and harmonic measure based at any point $z$ in $D$ are $A_\infty$-*equivalent*. This means that there are positive constants $c_1$, $c_2$, $\alpha_1$, and $\alpha_2$ such that whenever $S$ is a segment of $\partial D$ and $E$ is a Borel subset of $S$, then

$$\frac{|E|}{|S|} \leq c_1 \left[ \frac{\omega(z, E, D)}{\omega(z, S, D)} \right]^{\alpha_1} \quad \text{and} \quad \frac{\omega(z, E, D)}{\omega(z, S, D)} \leq c_2 \left[ \frac{|E|}{|S|} \right]^{\alpha_2} . \tag{13.3}$$

See [JK]. Moreover, if $D$ is bounded, then for all points $z \in D$ which satisfy

$$\text{dist}\,(z, \partial D) \geq C \, \text{diam}\,(D), \tag{13.4}$$

the constants $c_1$, $c_2$, $\alpha_1$, and $\alpha_2$ depend only on $C$ and the chord-arc constant of $\partial D$.

In the situation of Lemma 13.1, consider the domain $D = \pi^{-1}(\Lambda_J)$ whose boundary consists of $\widehat{I}$, $\widetilde{I}$, and two arcs $\gamma_1$ and $\gamma_2$ which project via $\pi$ to the two semicircles in the boundary of the leaf $\Lambda_J$. We prove (Lemma 13.2) that $\partial D$ is chord-arc, with a constant independent of $I$, and (Lemma 13.3) that the inequality (13.4) holds, with a constant independent of $I$, for the point $z_I \in D$ defined in Section 10. Then we prove Lemma 13.1.

**Lemma 13.2.** *The boundary $\partial D$ of $D$ is chord-arc, with a chord-arc constant which is independent of $I$.*

**Proof.** We follow the proof of a similar result by González [G]. The boundary of $D$ consists of $\widehat{I}$, $\gamma_1$, $\gamma_2$, the orthocircular arcs in $\widetilde{I}$, and the subset $\widetilde{I} \cap \partial \mathbb{D}$ of the unit circle. It is sufficient to show that:

1. Wherever two of $\widehat{I}$, $\gamma_1$, $\gamma_2$, and $\widetilde{I}$ meet, the angle they form is bounded below by some $\theta_0 > 0$ which is independent of $I$; and

2. The components of $\partial \Lambda_J \setminus K$ (which are $J = \pi(\widehat{I})$; the semicircles $\pi(\gamma_1)$ and $\pi(\gamma_2)$; and the open intervals in $E_J \setminus K$) are chord-arc, in the hyperbolic metric, with a chord-arc constant which is independent of $I$.

Then the components of $\partial D \setminus \partial \mathbb{D}$ are also chord-arc in the hyperbolic metric; they are chord-arc in the Euclidean metric; and the whole boundary $\partial D$ is chord-arc in the Euclidean metric with a constant depending only on the chord-arc constant of $\Lambda_J$.

1. $\widehat{I}$ meets $\gamma_1$ and $\gamma_2$ at right angles, because $J$ meets $\pi(\gamma_1)$ and $\pi(\gamma_2)$ at right angles and $\pi^{-1}$ is conformal.

Figure 12. $\Lambda_J$ and $D = \pi^{-1}(\Lambda_J)$ are chord-arc domains.

Let $\gamma_1$ be the arc in $\partial D$ such that $\pi(\gamma_1)$ is the smaller semicircle in $\partial \Lambda_J$, and let $A$ be the orthocircular arc in the lower part of $\partial D$ which meets $\gamma_1$. See Figure 12. We show that the hyperbolic distance between any two points $z_1 \in \gamma_1$ and $z_2 \in A$ is uniformly bounded away from zero; this implies that $\gamma_1$ and $A$ do not form a cusp at their common endpoint but make a positive angle there. Take a path from $\pi(z_1)$ to $\pi(z_2)$ in $\Lambda_J$. Let $\gamma$ be the segment of this path which starts at $\pi(z_1)$ and has length $\text{Im}(z_1)/2$. Then

$$
\begin{aligned}
d_{\text{hyp}}(z_1, z_2) &= d_{\text{hyp}}(\pi(z_1), \pi(z_2)) \\
&\geq c_\Omega \int_\gamma \frac{|dw|}{\text{dist}\,(w, K)} \\
&\geq c_\Omega \, \frac{1}{\frac{3}{2}\,\text{Im}\,\pi(z_1)} \cdot \frac{\text{Im}\,\pi(z_1)}{2} \\
&\geq c_\Omega/3.
\end{aligned}
\tag{13.5}
$$

So the angle formed by $\gamma_1$ and $A$ is bounded below by some $\theta_0 > 0$, independent of $I$.

The same reasoning applied to $\gamma_2$ shows that $\gamma_2$ makes a positive angle, uniformly bounded below, with $\partial \mathbb{D}$.

2. $J$ and the open intervals in $E_J \setminus K$ are geodesic arcs in $\Omega = \overline{\mathbf{C}} \setminus K$, so they are chord-arc in the hyperbolic metric. To show that the semicircle $\pi(\gamma_1)$ is hyperbolically chord-arc, take points $z$ and $\zeta$ in $\pi(\gamma_1)$. Let $\gamma$ be the segment of $\pi(\gamma_1)$ between $z$ and $\zeta$, and let $L$ be the chord between $z$ and $\zeta$. It is enough to consider $z$ and $\zeta$ in the part of $\pi(\gamma_1)$ near (in the Euclidean metric) to an endpoint of $\pi(\gamma_1)$. Then $\gamma$ is almost vertical, and

$$
\ell_{\text{hyp}}(\gamma) \sim \int_\gamma \frac{|dw|}{\text{dist}\,(w, K)} \sim \int_\gamma \frac{|dw|}{\text{Im}\,w}
\tag{13.6}
$$

is comparable to

$$
\ell_{\text{hyp}}(L) \sim \int_L \frac{|dw|}{\text{dist}\,(w, K)} \sim \int_L \frac{|dw|}{\text{Im}\,w}
\tag{13.7}
$$

with a uniform constant.

This completes the proof of Lemma 13.2. $\qquad\blacksquare$

**Lemma 13.3.** *For each standard grid interval $I \in \mathcal{H}$, let $z_I$ be the point in $D$, near $\widehat{I}$, defined in Section 10. There is a constant $C > 0$ independent of $I$ such that*

$$
\text{dist}\,(z_I, \partial D) \geq C \,\text{diam}\,(D).
\tag{13.8}
$$

**Proof.** We compare both sides of (13.8) to $1 - |w_I|$, where $w_I$ is the point in $\widehat{I} \subset \partial D$, near $z_I$, defined in Section 10. (See Figure 6.)

The hyperbolic distance from $z_I$ to $\partial D$ is uniformly bounded away from zero. To see this, consider a geodesic arc from $z_J = \pi(z_I)$ to $\partial \Lambda_J$. Let $\gamma$ be the segment of this geodesic which has $z_J$ as one endpoint and which has Euclidean length $\alpha|J|/2$. Then

$$
\begin{aligned}
d_{\text{hyp}}(z_I, \partial D) &= d_{\text{hyp}}(\pi(z_I), \pi(\partial D)) \\
&= d_{\text{hyp}}(z_J, \partial \Lambda_J) \\
&\geq c_\Omega \int_\gamma \frac{|dz|}{\text{dist}\,(z, K)} \\
&\geq c_\Omega \cdot \frac{1}{\frac{3}{2}\alpha|J| + \frac{1}{2}|J|} \cdot \frac{\alpha|J|}{2} \\
&= c\,(\Omega, \alpha).
\end{aligned}
\tag{13.9}
$$

The hyperbolic ball of radius $c\left(\Omega, \alpha\right)$ centred at $z_I$ contains a Euclidean ball of radius $c'\left(1 - |z_I|\right)$, centred at $z_I$. Here $c'$ is independent of $I$. Since $(1 - |z_I|) \sim (1 - |w_I|)$, with a constant independent of $I$, we conclude that

$$\operatorname{dist}\left(z_I, \partial D\right) \geq c\left(1 - |w_I|\right), \tag{13.10}$$

where $c$ is independent of $I$.

To estimate $\operatorname{diam}\left(D\right)$, it is sufficient to show that

$$\operatorname{dist}\left(x, \widehat{I}\right) \leq c\left(1 - |w_I|\right) \tag{13.11}$$

for all $x \in \partial D$. For if $x$ and $y$ are in $\partial D$, then

$$|x - y| \leq \operatorname{dist}\left(x, \widehat{I}\right) + \operatorname{dist}\left(y, \widehat{I}\right) + |\widehat{I}|, \tag{13.12}$$

and we saw in Section 10 that $|\widehat{I}| \sim 1 - |w_I|$. We showed (13.11) for $x \in \widetilde{I}$, in the proof of Lemma 10.5. We prove (13.11) for $x \in \gamma_1$; the same argument works for $x \in \gamma_2$.

Let $a$ and $\widehat{a}$ be the endpoints of $\gamma_1$ which lie in $I$ and $\widehat{I}$ respectively. Using the fact that $\gamma_1$ is chord-arc,

$$\begin{aligned}
\operatorname{dist}\left(x, \widehat{I}\right) &\leq |x - \widehat{a}| \\
&\leq \ell_{\gamma_1}(x, \widehat{a}) \\
&\leq \ell_{\gamma_1}(a, \widehat{a}) \\
&\leq c\,|a - \widehat{a}| \\
&\leq |I| + |\widehat{I}| + \operatorname{dist}\left(I, \widehat{I}\right).
\end{aligned} \tag{13.13}$$

Now $|I| \leq c\left(1 - |w_I|\right)$ by Lemma 8.1; $|\widehat{I}| \leq c\,|I|$ by Lemma 10.3, and $\operatorname{dist}\left(I, \widehat{I}\right) \leq c\left(1 - |w_I|\right)$ by the proof of Lemma 10.5. So we have established (13.11), which together with (13.10) proves Lemma 13.3. □

**Proof of Lemma 13.1.** By Lemma 13.2, Lemma 13.3, and the remarks after the statement of Lemma 13.1, arclength on $\partial D$ and harmonic measure at $z_I$ in $D$ are $A_\infty$-equivalent with constants depending only on the chord-arc constant of $\Lambda_J$ and the constant $C$ from Lemma 13.3. Let $c_1$ and $\alpha_1$ be such that

$$\frac{|E|}{|S|} \leq c_1 \left[\frac{\omega(z_I, E, D)}{\omega(z_I, S, D)}\right]^{\alpha_1} \tag{13.14}$$

for all Borel subsets $E$ of segments $S$ of $\partial D$. The leaves $\Lambda_J$ are also chord-arc, and $d_{\mathrm{hyp}}(z_J, \partial \Lambda_J) \geq c > 0$, so we may choose $c_2$ and $\alpha_2$ independent of $I$ such that

$$\frac{\omega(z_J, E, \Lambda_J)}{\omega(z_J, S, \Lambda_J)} \leq c_2 \left[\frac{|E|}{|S|}\right]^{\alpha_2} \tag{13.15}$$

for all Borel subsets $E$ of segments $S$ of $\partial \Lambda_J$.

By the distortion estimates in Section 10, $|\widetilde{I}| \sim |I|$ and $|B| \sim |\widehat{B}|$. We use (13.14) and (13.15) on the sets $\widehat{B} \subset \widetilde{I} \subset \partial D$ and $\pi(\widehat{B}) \subset E_J \subset \partial \Lambda_J$:

$$
\begin{aligned}
\frac{|B|}{|I|} &\le c \, \frac{|\widehat{B}|}{|\widetilde{I}|} \\
&\le c \, c_1 \left[ \frac{\omega(z_I, \widehat{B}, D)}{\omega(z_I, \widetilde{I}, D)} \right]^{\alpha_1} \\
&= c \, c_1 \left[ \frac{\omega(\pi(z_I), \pi(\widehat{B}), \pi(D))}{\omega(\pi(z_I), \pi(\widetilde{I}), \pi(D))} \right]^{\alpha_1} \\
&= c \, c_1 \left[ \frac{\omega(z_J, \pi(\widehat{B}), \Lambda_J)}{\omega(z_J, E_J, \Lambda_J)} \right]^{\alpha_1} \\
&\le c \, c_1 \left[ c_2 \left( \frac{|\pi(\widehat{B})|}{|E_J|} \right)^{\alpha_2} \right]^{\alpha_1},
\end{aligned}
\tag{13.16}
$$

which establishes the right hand inequality in (13.1). A similar argument, again using (13.3), proves the left hand inequality in (13.1). The constants are independent of $I$ and $B$, but they depend on the $N$ such that $|E_J| = 2 \cdot 3^N |J|$, since the chord-arc constant of $\Lambda_J$ depends on this $N$. $\qquad\square$

**Lemma 13.4.** *There are positive constants $c_1$ and $c_2$ such that for all $n \ge 1$, for every standard grid interval $I \in \mathcal{H}_{n-1}$,*

$$
EX_i = \frac{1}{\mu(I)} \int_I X_i(x) \, d\mu \ge c_1;
\tag{13.17}
$$

*and*

$$
EX_i^2 = \frac{1}{\mu(I)} \int_I X_i^2(x) \, d\mu \le c_2.
\tag{13.18}
$$

*Here $X_i(x) = \log_3\big(|J_i(x)|/|J|\big)$ is the auxiliary function for a particle which makes its $i$th jump from $J = \pi(\widehat{I})$ and which has not yet reached $J_\infty$.*

**Proof.** Let $L$ be the component of $[0,1] \setminus K$ which contains $J = \pi(\widehat{I})$.

We begin by slightly modifying some definitions, in order to simplify the calculations. In each component $L'$ of $E_J \setminus K$, we amalgamated the $2N$ Whitney intervals in the centre of $L'$ into a single interval $J_c = J_{-N} \cup \cdots \cup J_N$. We now temporarily assume that we have retained the individual intervals $J_{-N}, \ldots, J_N$ instead. We also change the definitions of the region $R_3$ and of the function $F_n$. Set $S_3 = J_{-1} \cup J_1$ in the component $L'$ of $[0,1] \setminus K$ which is the left half of $E_J$, and let $R_3 = P \circ \pi^{-1}(S_3)$. (In Section 11, we had $S_3 = J_{-N} \cup \cdots \cup J_N \subset L'$. The regions $R_1$ and $R_5$ are unchanged; the regions $R_2$ and $R_4$ become larger since they now include $J_{-N} \cup \cdots \cup J_{-2}$ and $J_2 \cup \cdots \cup J_N$ respectively. )

Define

$$F_n(x) = \begin{cases} 1, & x \in R_1 \cup R_5; \\ (1-\varepsilon)\frac{|I|}{|R_3|}, & x \in R_3; \\ \delta, & x \in R_2 \cup R_4, \end{cases} \qquad (13.19)$$

using the new definition of $R_3$. (At the end of the proof we show that the estimates with these new definitions imply the estimates for the original functions $F_n$ and $X_i$.)

With these assumptions, we begin with (13.17). Let $I \in \mathcal{H}_{n-1}$ be a standard interval. Let $J = \pi(\widehat{I})$. For each point $x$ in $I$, let $I_i(x)$ be the interval in $\mathcal{H}_n$ which contains $x$. $I$, $n$, and $i$ are fixed throughout the proof. Let $J_i(x) = \pi(\widehat{I_i}(x))$. Then

$$X_i(x) = \log_3\left[\frac{|J_i(x)|}{|J|}\right]. \qquad (13.20)$$

Notice that to compute $EX_i$ or $EX_i^2$, we need only deal with $F_n$. The functions $F_l$ for $l \neq n$ are irrelevant, since for $l < n$ $F_l$ is constant on $I$, and for $l > n$ $F_l$ has mean value one on each subset of $I$ where $X_i$ is constant.

Let

$$\begin{aligned} B_k &= \{x \in I \mid X_i(x) = N - k\} \\ &= \{x \in I \mid |J_i(x)| = 3^{N-k}|J|\} \\ &= \{x \in I \mid |J_i(x)| = 2^{-1} \cdot 3^{-k}|E_J|\}, \end{aligned} \qquad (13.21)$$

for $k \geq 1$. These $B_k$'s are exactly the sets on which the integrand $X_i$ in (13.17) is constant. The union of the $B_k$'s is $I$. We estimate (13.17) by integrating over each $B_k$ separately (see (13.27) below); before that we need some preliminaries.

Let $\widehat{B_k}$ be the subset of $\widetilde{I}$ such that $P(\widehat{B_k}) = B_k$. Then $\pi(\widehat{B_k})$ is the union of all Whitney intervals in $E_J$ of length $2^{-1} \cdot 3^{-k}|E_J|$.

We count the number of Whitney intervals of this length in $E_J$. The segment $E_J$ contains one component of $\overline{\mathbb{R}} \setminus K$ of length $2^{-1}|E_J|$, and it contains $2^j$ components of $\overline{\mathbb{R}} \setminus K$ of length $2^{-1} \cdot 3^{-j-1}|E_J|$, for each $j \geq 0$. Within each component $L$ of $\overline{\mathbb{R}} \setminus K$ there are two Whitney intervals of size $3^{-l}|L|$, for each $l \geq 1$. Therefore $E_J$ contains $2^k$ Whitney intervals of length $2^{-1} \cdot 3^{-k}|E_J| = 3^{N-k}|J|$, for each $k \geq 1$. Hence

$$|\pi(\widehat{B_k})| = 2^k \cdot 2^{-1} \cdot 3^{-k}|E_J| = \frac{1}{2}\left(\frac{2}{3}\right)^k |E_J|, \qquad (13.22)$$

for $k \geq 1$.

Therefore, by Lemma 13.1, there are positive constants $c$ and $\alpha$, independent of $I$ but dependent on $N$, such that

$$\frac{|B_k|}{|I|} \leq c\left[\frac{|\pi(\widehat{B_k})|}{|E_J|}\right]^\alpha = c\left[\frac{1}{2}\left(\frac{2}{3}\right)^k\right]^\alpha = c\,\beta^k, \qquad (13.23)$$

where $\beta = (2/3)^\alpha$ is strictly less than one.

We also estimate the sizes of $B_k \cap R_1$ and $B_k \cap R_5$, for $k \geq 1$. The sets $R_j$ and $S_j$, $1 \leq j \leq 5$, were defined in Section 11. $S_1 = \pi(\widehat{R_1})$ is the segment at one end of $E_J$, contained in $E_J \setminus K$, such that $|S_1| = \frac{3}{2} \cdot 3^{-Q}|J|$. $Q$ is a large integer, independent of $I$. The largest Whitney interval in $S_1$ has length $3^{-Q}|J| = 2^{-1} \cdot 3^{-Q-N}|E_J|$. The segment $S_1$ contains exactly one Whitney interval of length $2^{-1} \cdot 3^{-k}|E_J|$, in other words one Whitney interval which lies in $\pi(\widehat{B_k})$, for each $k \geq Q + N$. By Lemma 13.1,

$$
\begin{aligned}
\frac{|B_k \cap R_1|}{|I|} &\leq c \left[ \frac{|\pi(\widehat{B_k} \cap \widehat{R_1})|}{|E_J|} \right]^\alpha \\
&= c \left[ \frac{|\pi(\widehat{B_k}) \cap S_1|}{|E_J|} \right]^\alpha \\
&= c \left[ \frac{2^{-1} \cdot 3^{-k}\,|E_J|}{|E_J|} \right]^\alpha \\
&= c\,(3^{-\alpha})^k,
\end{aligned}
\tag{13.24}
$$

for $k \geq Q + N$. The left hand side is zero when $0 \leq k \leq Q + N - 1$.

$S_5 = \pi(\widehat{R_5})$ is the segment at the other end of $E_J$ such that $|S_5| = 3^{-Q}|J|$. $S_5$ is a closed construction interval of the Cantor set $K$. The largest Whitney intervals in $S_5$ are two intervals of length $3^{-2} \cdot 3^{-Q}|J| = 2^{-1} \cdot 3^{-Q-N-2}|E_J|$. $S_5$ contains $2^{k-Q-N-1}$ Whitney intervals of size $2^{-1} \cdot 3^{-k}|E_J|$, for each $k \geq Q + N + 2$. By Lemma 13.1,

$$
\begin{aligned}
\frac{|B_k \cap R_5|}{|I|} &\leq c \left[ \frac{|\pi(\widehat{B_k}) \cap S_5|}{|E_J|} \right]^\alpha \\
&= c \left[ \frac{2^{k-Q-N-1} \cdot 2^{-1} \cdot 3^{-k}|E_J|}{|E_J|} \right]^\alpha \\
&= c\,(2^{-Q-N})^\alpha (2/3)^{k\alpha} \\
&= c\,(2^{-\alpha})^{Q+N} \beta^k,
\end{aligned}
\tag{13.25}
$$

for $k \geq Q + N + 2$. The left hand side vanishes for $0 \leq k \leq Q + N + 1$.

The set $B_1$ is exactly $R_3$, and so $F_n \equiv (1 - \varepsilon)\,|I|/|B_1|$ on $B_1$. Also, on $I$, $\mu_{n-1}$ is given by $d\mu_{n-1} = (\mu_{n-1}(I)/|I|)\,dx$. Therefore

$$
\begin{aligned}
\mu(B_1) &= \mu_n(B_1) \\
&= \int_{B_1} F_n(x)\,d\mu_{n-1}(x) \\
&= (1 - \varepsilon)\,\frac{|I|}{|B_1|}\,\frac{\mu_{n-1}(I)}{|I|}\,|B_1| \\
&= (1 - \varepsilon)\,\mu_{n-1}(I) \\
&= (1 - \varepsilon)\,\mu(I).
\end{aligned}
\tag{13.26}
$$

After these preliminaries we can estimate $EX_i$:

$$
\begin{aligned}
EX_i &= \frac{1}{\mu(I)} \int_I \log_3 \left[ \frac{|J_i(x)|}{|J|} \right] d\mu(x) \\
&= \sum_{k=1}^{\infty} (N-k) \frac{\mu(B_k)}{\mu(I)} \\
&= (N-1) \frac{\mu(B_1)}{\mu(I)} + \sum_{k=2}^{\infty} (N-k) \frac{\mu(B_k)}{\mu(I)} \\
&\geq (N-1)(1-\varepsilon) + \sum_{k=N}^{\infty} (N-k) \frac{\mu(B_k)}{\mu(I)} \\
&\geq (N-1)(1-\varepsilon) - \sum_{k=N}^{\infty} k \frac{\mu(B_k)}{\mu(I)} \\
&= (N-1)(1-\varepsilon) - \sum_{k=N}^{\infty} k \frac{\mu_n(B_k)}{\mu_{n-1}(I)} \\
&\geq (N-1)(1-\varepsilon) - \sum_{k=N}^{\infty} \frac{k}{\mu_{n-1}(I)} \int_{B_k} F_n(x) \, d\mu_{n-1}(x).
\end{aligned}
\tag{13.27}
$$

The last series converges. To show this, we split each $B_k$ into three pieces, $B_k \cap (R_2 \cup R_4)$, $B_k \cap R_1$, and $B_k \cap R_5$, and estimate the sums over the three types of terms.

First, since $F_n \equiv \delta$ on $R_2 \cup R_4$,

$$
\begin{aligned}
\sum_{k=N}^{\infty} \frac{k}{\mu_{n-1}(I)} \int_{B_k \cap (R_2 \cup R_4)} & F_n(x) \, d\mu_{n-1}(x) \\
&= \sum_{k=N}^{\infty} k \, \delta \, \frac{\mu_{n-1}\big(B_k \cap (R_2 \cup R_4)\big)}{\mu_{n-1}(I)} \\
&= \delta \sum_{k=N}^{\infty} k \, \frac{|B_k \cap (R_2 \cup R_4)|}{|I|} \\
&\leq \delta \sum_{k=N}^{\infty} k \, \frac{|B_k|}{|I|} \\
&\leq \delta \, c\,(N) \sum_{k=N}^{\infty} k \, \beta^k \\
&\leq \delta \, \frac{c\,(N)}{(1-\beta)^2},
\end{aligned}
\tag{13.28}
$$

which can be made arbitrarily small by choosing $\varepsilon$ small (which implies $\delta$ small) in the definition of $F_n$. Here the second equality holds because, on $I$, $\mu_{n-1}$ is a constant multiple of Lebesgue measure, and the second last line holds by (13.23).

Similarly, since $F_n \equiv 1$ on $R_1$,

$$\sum_{k=N}^{\infty} \frac{k}{\mu_{n-1}(I)} \int_{B_k \cap R_1} F_n(x)\, d\mu_{n-1}(x)$$

$$= \sum_{k=N}^{\infty} k\, \frac{\mu_{n-1}(B_k \cap R_1)}{\mu_{n-1}(I)}$$

$$= \sum_{k=N}^{\infty} k\, \frac{|B_k \cap R_1|}{|I|}$$

$$= \sum_{k=N+Q}^{\infty} k\, \frac{|B_k \cap R_1|}{|I|} \tag{13.29}$$

$$\leq c(N) \sum_{k=N+Q}^{\infty} k\, \left(3^{-\alpha}\right)^k,$$

which can be made arbitrarily small by choosing $Q$ large. Here the last line holds by (13.24).

Finally, since $F_n \equiv 1$ on $R_5$,

$$\sum_{k=N}^{\infty} \frac{k}{\mu_{n-1}(I)} \int_{B_k \cap R_5} F_n(x)\, d\mu_{n-1}(x)$$

$$= \sum_{k=N}^{\infty} k\, \frac{\mu_{n-1}(B_k \cap R_5)}{\mu_{n-1}(I)}$$

$$= \sum_{k=N}^{\infty} k\, \frac{|B_k \cap R_5|}{|I|}$$

$$= \sum_{k=Q+N+2}^{\infty} k\, \frac{|B_k \cap R_5|}{|I|} \tag{13.30}$$

$$\leq c(N)\, (2^{-\alpha})^{Q+N} \sum_{k=Q+N+2}^{\infty} k\, \beta^k$$

$$\leq \frac{c(N)\, (2^{-\alpha})^{Q+N}}{(1-\beta)^2},$$

which decreases to zero as $Q$ increases. Here the second last inequality holds by (13.25).

Once $N$ is fixed, we can make (13.28), (13.29), and (13.30) arbitrarily small by choosing $\varepsilon$ sufficiently small and $Q$ sufficiently large. Therefore we can make the last term in (13.27) arbitrarily small. So for some positive constant $c_1$ depending on $N$, we can ensure that $EX_i \geq c_1$ for all standard intervals $I$ in $\mathcal{H}_{n-1}$, for all $n \geq 1$.

We can estimate the second moment more simply. We showed in Section 11 that there is an upper bound $1/\delta_0$ for $F_n$ on $I$, independent of $I$ and $n$. Then

$$
\begin{aligned}
EX_i^2 &= \frac{1}{\mu(I)} \int_I \left\{ \log_3 \left[ \frac{|J_i(x)|}{|J|} \right] \right\}^2 d\mu(x) \\
&= \sum_{k=1}^{\infty} (N-k)^2 \, \frac{\mu(B_k)}{\mu(I)} \\
&= \sum_{k=1}^{\infty} (N-k)^2 \, \frac{\mu_n(B_k)}{\mu_{n-1}(I)} \\
&= \sum_{k=1}^{\infty} (N-k)^2 \, \frac{1}{\mu(I)} \int_{B_k} F_n(x) \, d\mu_{n-1}(x) \\
&\leq \frac{1}{\delta_0} \sum_{k=1}^{\infty} (N-k)^2 \, \frac{\mu_{n-1}(B_k)}{\mu_{n-1}(I)} \\
&= \frac{1}{\delta_0} \sum_{k=1}^{\infty} (N-k)^2 \, \frac{|B_k|}{|I|} \\
&\leq \frac{1}{\delta_0} \sum_{k=1}^{\infty} (N^2 + k^2) \, \frac{|B_k|}{|I|} \\
&\leq \frac{c(N)}{\delta_0} \sum_{k=1}^{\infty} (N^2 + k^2) \, \beta^k \\
&\leq \frac{c(N)}{\delta_0} \left[ \frac{N^2}{(1-\beta)} + \frac{1}{(1-\beta)^3} \right],
\end{aligned}
\tag{13.31}
$$

using (13.23) in the second last line. In the sixth line we used the fact that $\mu_{n-1}$ is a constant multiple of Lebesgue measure on $I$. Again, once $N$, $Q$, and $\varepsilon$ are fixed, this estimate gives a uniform upper bound $c_2$ on $EX_i^2$ for all standard intervals $I \in \mathcal{H}_{n-1}$, for all $n \geq 1$.

It remains to prove the same estimates for our original functions $X_i$ and $F_n$. Write $Y_i$ and $G_n$ for the modified versions of $X_i$ and $F_n$ defined at the start of the proof and used in the calculations (13.19)–(13.31) above.

To estimate $EX_i^2$, define $Z_i(x) = N^2 + k^2$ on $B_k$, $k \geq 1$, where $B_k = \{x \in I \mid Y_i(x) = N - k\}$ as above. A comparison of $Y_i$, $X_i$, and $Z_i$ shows that on each $B_k$,

$$
X_i^2 \leq (N+k)^2 \leq 2 \, (N^2 + k^2) = 2 \, Z_i.
\tag{13.32}
$$

The calculation (13.31) shows that the mean value of $Z_i G_n$ on $I$ is at most $c_2$. Also, $F_n \sim G_n$ with constants independent of $I$, since they are both bounded above and below by positive constants independent of $I$. Therefore

$$
\begin{aligned}
EX_i^2 &= \frac{1}{|I|} \int_I X_i^2 F_n \, dx \\
&\leq c \, \frac{1}{|I|} \int_I 2 \, Z_i G_n \, dx \\
&\leq 2 \, c \, c_2,
\end{aligned}
\tag{13.33}
$$

for all $I$.

Finally, we estimate $EX_i$. Let $G_b$ (respectively $F_b$) be the maximum value of $G_n$ (respectively $F_n$). Then $G_b \sim F_b$ with constants independent of $I$, by a comparison of harmonic measures in $\Lambda_J$. Let

$$
I_+ = \{ x \in I \mid F_n(x) = F_b \}
\tag{13.34}
$$

and

$$
I_s = \{ x \in I \mid F_n(x) = \delta' \},
\tag{13.35}
$$

where $\delta'$ is the minimum value of $F_n$ on $I$. In the next calculation we neglect the region where $Y_i \equiv 1$. Then, since $X_i \geq Y_i$ on $I$,

$$
\begin{aligned}
EX_i &= \frac{1}{|I|} \int_I X_i F_n \, dx \\
&\geq \frac{1}{|I|} \int_I Y_i F_n \, dx \\
&\geq \frac{1}{|I|} \int_{I_+} Y_i F_b \, dx - \frac{1}{|I|} \int_{I_s} Y_i F_n \, dx \\
&\geq c \, \frac{1}{|I|} \int_{I_+} Y_i G_b \, dx - \frac{1}{|I|} \int_{I_s} Y_i F_n \, dx \\
&\geq c \, \frac{N-1}{2} - \delta' \, \frac{1}{|I|} \int_{I_s} Y_i \, dx;
\end{aligned}
\tag{13.36}
$$

by (13.27). We showed above that $\frac{1}{|I|} \int_{I_s} Y_i \, dx$ is bounded above by a constant independent of $I$. Therefore, by choosing $\delta'$ small enough in the definition of $F_n$, we may ensure that $EX_i$ is bounded below by a positive constant independent of $I$.

This completes the proof of Lemma 13.4. $\qquad\qquad\square$

## 14. Estimates $EX_i \geq c_1$ and $EX_i^2 \leq c_2$ for non-standard intervals

Let $I \in \mathcal{H}$ be any non-standard grid interval such that $J_c = \pi(\widehat{I})$ is the non-standard central interval in some component $L$ of $[0,1] \setminus K$. In this section we prove (Lemma 14.3) the estimates $\mu(I)^{-1} \int_I X_i \, d\mu \geq c_1$ and $\mu(I)^{-1} \int_I X_i^2 \, d\mu \leq c_2$, where $c_1$ and $c_2$ are positive constants independent of $I$. Here $X_i(x) = \log_3\big(|J_i(x)|/|J_c|\big)$ is the auxiliary function for a particle which makes its $i^{\text{th}}$ jump from $J_c = \pi(\widehat{I})$ and which has not yet reached $J_\infty$. We begin with an estimate, analogous to Lemma 13.1, on the distortion in length caused by the map $P \circ \pi^{-1}$.

**Lemma 14.1.** *Let $I \in \mathcal{H}_{n-1}$, $n \geq 1$, be a non-standard grid interval such that $J_c = \pi(\widehat{I})$ is the central interval in some component $L$ of $[0,1] \setminus K$. Fix a number $\lambda > 1$, and let $q$ be any positive integer. Let $B$ and $T$ be unions of grid intervals from $\mathcal{H}_n$ such that $B \subset T \subset I$ and $\pi(\widehat{T}) \subset \{w \mid \lambda^q \leq \operatorname{dist}(w, J_c) \leq \lambda^{q+1}\}$. There are constants $c$ and $\alpha$, independent of $q$, $B$, $T$, $I$, and $n$, such that*

$$\frac{|B|}{|T|} \leq c \left[ \frac{|\pi(\widehat{B})|}{|\pi(\widehat{T})|} \right]^\alpha. \tag{14.1}$$

The proof is essentially the same as for standard intervals, but we need a stronger version of $A_\infty$-equivalence because $\partial \Lambda_J$ is unbounded. Namely, if $\Omega$ is a chord-arc domain, and $z \in \Omega$, then there are positive constants $c_1$, $c_2$, $\alpha_1$, and $\alpha_2$ such that

$$\frac{|E|}{|S|} \leq c_1 \left[ \frac{\omega(z, E, \Omega)}{\omega(z, S, \Omega)} \right]^{\alpha_1} \quad \text{and} \quad \frac{\omega(z, E, \Omega)}{\omega(z, S, \Omega)} \leq c_2 \left[ \frac{|E|}{|S|} \right]^{\alpha_2} \tag{14.2}$$

whenever $S$ is a segment of $\partial \Omega$ satisfying

$$S \subset \{w \mid \lambda^q \leq \operatorname{dist}(w, z) \leq \lambda^{q+1}\} \tag{14.3}$$

for any positive integer $q$ and for a fixed $\lambda > 1$, and $E$ is a Borel subset of $S$. The constants $c_1$, $c_2$, $\alpha_1$, and $\alpha_2$ depend on the chord-arc constant of $\partial \Omega$ and on the constant $\lambda$, but not on $q$.

This follows from the case of bounded chord-arc domains. By dilations and translations we may assume that $z = i$, and that $-i$ is at least distance $\varepsilon$ from $\Omega$. Consider the Möbius transformation $\tau : w \mapsto (w + i)^{-1}$. This maps $\Omega$ to a bounded chord-arc domain. Harmonic measure is invariant under Möbius transformations. Also, for each positive integer $q$, the map $\tau$ scales everything in the annulus in (14.3), in particular $E$ and $S$, by approximately the same factor. The estimates in (14.2) follow.

Consider the domain $D = \pi^{-1}(\Lambda_{J_c})$, where $I$ is a non-standard interval as above. The analogues of Lemmas 13.2 and 13.3 hold, and Lemma 14.1 follows in the same way as Lemma 13.1 follows from Lemmas 13.2 and 13.3.

Notice that, although the interval $I$ is no longer standard, arclength and harmonic measure are still $A_\infty$-equivalent on $\partial D$, since $D$ is a bounded chord-arc domain. Therefore, for any segment $S$ of $\partial D$ and Borel subset $E$ of $S$, we still have

$$\frac{|E|}{|S|} \leq c_1 \left[ \frac{\omega(z_J, \pi(E), \Lambda_{J_c})}{\omega(z_J, \pi(S), \Lambda_{J_c})} \right]^{\alpha_1}, \tag{14.4}$$

and the analogous inequality in the other direction, both with constants independent of $I$. We will frequently use this fact, combined with direct estimates of harmonic measures in $\Lambda_{J_c}$. We refer to such an argument as a *comparison of harmonic measure on the domain side.*

We will need the next lemma, which follows immediately from our decomposition of $\overline{\mathbb{R}} \setminus K$ into Whitney intervals.

**Lemma 14.2.** *Let $A$ be a subinterval of $[0,1]$ of length $3^{-l}$, $l \geq 0$. Then*

a) *$A$ contains a Whitney interval of length at least $3^{-2}|A|$; and*

b) *$A$ meets at most $2^k$ Whitney intervals of length $3^{-k}|A|$, for $k \geq 1$.*

The main result of this section is:

**Lemma 14.3.** *There are positive constants $c_1$ and $c_2$ such that for all $n \geq 1$, for every non-standard grid interval $I \in \mathcal{H}_{n-1}$ such that $\pi(\widehat{I}) = J_c$ is the central Whitney interval in some component $L$ of $[0,1] \setminus K$,*

$$EX_i = \frac{1}{\mu(I)} \int_I X_i(x)\, d\mu \geq c_1; \tag{14.5}$$

*and*

$$EX_i^2 = \frac{1}{\mu(I)} \int_I X_i^2(x)\, d\mu \leq c_2. \tag{14.6}$$

*Here $X_i(x) = \log_3\big(|J_i(x)|/|J_c|\big)$ is the auxiliary function for a particle which has not yet reached $J_\infty$ and which makes its $i^{\text{th}}$ jump from $J_c = \pi(\widehat{I})$.*

**Proof.** Let $J_c$ be the central Whitney interval, $J_c = J_{-N} \cup \cdots \cup J_N$, in some component $L$ of $[0,1] \setminus K$. Make the convention that $|J_c| = |L|/3$. The tips of the leaves for the standard intervals in $L$ cover $K_l$ and $K_r$, the closed construction intervals of $K$ immediately to the left and right of $L$. So the tip $E_{J_c}$ of the leaf for $J_c$ is $\overline{\mathbb{R}} \setminus (K_l \cup L \cup K_r)$. See Figure 13.

To simplify the calculations, we again assume, as in the proof of Lemma 13.4, that in each component $L$ of $E_{J_c} \setminus K$ we have not amalgamated the $2N$ central Whitney intervals but have retained them individually. As before, the result follows from this simplified version.

Figure 13. A non-standard central interval $J_c$.

We divide $E_{J_c}$ into subsets $V_j$, $j \geq 1$, of length $|V_j| = 3^{j+10}|J_c|$. Let each $V_j$ have two components of equal length, one on each side of $J_c$. Put the components of $V_1$ as close as possible to $J_c$, that is, on each side of and adjacent to $K_l \cup L \cup K_r$. Put the components of $V_2$ on each side of and adjacent to $V_1 \cup K_l \cup L \cup K_r$, and so on. For each $j \geq 1$, let $T_j$ be the subset of $I$ such that $\pi(\widehat{T_j}) = V_j$.

Fix $j \geq 1$, and for each $k \geq 1$ define the subset $B_k$ of $T_j$ by

$$
\begin{aligned}
B_k &= \{x \in T_j \mid |J_i(x)| = 3^{-k}|V_j|\} \\
&= \{x \in T_j \mid |J_i(x)| = 3^{-k} \cdot 3^{j+10}|J_c|\} \\
&= \{x \in T_j \mid X_i(x) = j + 10 - k\}.
\end{aligned}
\tag{14.7}
$$

By Lemma 14.2, the set $\pi(\widehat{B_k})$ meets at most $2^k$ Whitney intervals of size $3^{-k}|V_j|$, so

$$
|\pi(\widehat{B_k})| \leq \left(\frac{2}{3}\right)^k |V_j|.
\tag{14.8}
$$

By Lemma 14.1,

$$
\frac{|B_k|}{|T_j|} \leq c \left[\frac{|\pi(\widehat{B_k})|}{|V_j|}\right]^\alpha = c \left[\left(\frac{2}{3}\right)^k\right]^\alpha = c\,\beta^k,
\tag{14.9}
$$

where $\beta = (2/3)^\alpha$ is strictly less than one, and $c$ and $\beta$ are independent of $j$.

In fact, Lemma 14.1 in the form stated above may not apply, since $B_k$ and $T_j$ need not be unions of whole grid intervals. However, the difficulty is only that at either end of $T_j$ there may be grid intervals which meet $T_j$ (and hence some $B_k$) but are not contained in $T_j$ (or $B_k$). The distortion estimates in Section 10, and the fact that grid intervals are comparable in size to their neighbours, imply that Lemma 14.1 can be extended to cover this case.

We first prove the estimate on $EX_i^2$. In estimating $X_i^2$ we can neglect the function $F_n$, since the $F_n$ we use will be bounded above by a constant independent of $I$. So it is enough to find a bound independent of $I$ for

$$\frac{1}{|I|} \int_I X_i^2 \, dx = \sum_{j=1}^{\infty} \frac{|T_j|}{|I|} \, \frac{1}{|T_j|} \int_{T_j} \left\{ \log_3 \left[ \frac{|J_i(x)|}{|J_c|} \right] \right\}^2 \, dx. \tag{14.10}$$

First, for each $j$,

$$\begin{aligned}
\frac{1}{|T_j|} \int_{T_j} \left\{ \log_3 \left[ \frac{|J_i(x)|}{|J_c|} \right] \right\}^2 &= \sum_{k=1}^{\infty} \left\{ j + 10 - k \right\}^2 \frac{|B_k|}{|T_j|} \\
&\leq c \sum_{k=1}^{\infty} \left\{ j^2 + (k-10)^2 \right\} \beta^k \\
&\leq c' j^2 + c''.
\end{aligned} \tag{14.11}$$

Also, $|T_j|/|I|$ decays exponentially in $j$:

$$\begin{aligned}
\frac{|T_j|}{|I|} &\leq c \, \frac{|\widehat{T_j}|}{|\widetilde{I}|} \\
&\leq c \, c_1 \left[ \frac{\omega\left(z_I, \widehat{T_j}, \pi^{-1}(\Lambda_{J_c})\right)}{\omega\left(z_I, \widetilde{I}, \pi^{-1}(\Lambda_{J_c})\right)} \right]^{\alpha_1} \\
&= c \, c_1 \left[ \frac{\omega(z_J, V_j, \Lambda_{J_c})}{\omega(z_J, E_{J_c}, \Lambda_{J_c})} \right]^{\alpha_1} \\
&\leq c \left[ \omega\left(z_J, V_j, \mathbf{U}\right) \right]^{\alpha_1},
\end{aligned} \tag{14.12}$$

with constants independent of $j$. Now $\pi \cdot \omega(z_J, V_j, \mathbf{U})$ is the angle subtended at $z_J$ by $V_j$. By elementary trigonometry, this angle is less than $c \, 3^{-j}$, with $c$ independent of $j$. Therefore

$$\frac{|T_j|}{|I|} \leq c \, 3^{-\alpha_1 j}. \tag{14.13}$$

Hence

$$\frac{1}{|I|} \int_I X_i^2 \, dx \leq c \sum_{j=1}^{\infty} 3^{-\alpha_1 j} (c' j^2 + c''). \tag{14.14}$$

The right hand side is finite and independent of $I$.

It remains to prove the estimate on $EX_i$. Fix $\varepsilon > 0$, and $j \geq 2$. By Lemma 14.2, $B_2$ is not empty. (For those $V_j \not\subset [0,1]$, $B_2$ may be empty; part a) of Lemma 14.2 need not hold as stated because of the amalgamation of the large Whitney intervals far from $[0,1]$ to form $J_\infty$. However our argument can be modified to cover this case.) Let $B = B_1 \cup B_2$, and define $F_n$ on $T_j$ by

$$F(x) = F_n(x) = \begin{cases} (1-\varepsilon)\dfrac{|T_j|}{|B|}, & x \in B; \\ \delta, & x \in T_j \setminus B, \end{cases} \tag{14.15}$$

where $\delta = \varepsilon |T_j|/|T_j \setminus B|$ is chosen so that $F$ has mean value one on $T_j$. A comparison of harmonic measure shows that $|T_j|/|B|$ is uniformly bounded for all $j$ and $I$. Therefore $F$ is uniformly bounded above and below by positive constants independent of $j$ and $I$.

On $B$, $X_i \geq j + 8$. Then

$$\begin{aligned}
\frac{1}{\mu(T_j)} \int_{T_j} X_i \, d\mu &= \frac{1}{|T_j|} \int_{T_j} X_i F \, dx \\
&\geq (j+8)(1-\varepsilon) + \delta \sum_{k=3}^{\infty} (j+10-k)\frac{|B_k|}{|T_j|} \\
&\geq (j+8)(1-\varepsilon) - \delta \sum_{k=j+10}^{\infty} (k-j-10)\frac{|B_k|}{|T_j|} \\
&\geq (j+8)(1-\varepsilon) - \delta \, c(N) \sum_{k-j=10}^{\infty} (k-j-10)\,\beta^k.
\end{aligned} \tag{14.16}$$

The sum is finite and independent of $j$. Hence, by choosing $\delta$ (or equivalently $\varepsilon$) small enough in the definition of $F$, we may ensure that the mean of $X_i F$ with respect to $\mu$ on $T_j$ is at least some positive constant $c_1$, independent of $j$.

To ensure that $F \equiv 1$ on suitable intervals at each end of $I$, we define $F$ slightly differently on $T_1$. Specifically, let $I_l$ and $I_r$ be subintervals at each end of $I$. Then $\pi(\widehat{I_l})$ and $\pi(\widehat{I_r})$ are subintervals of $V_1$, one on each side of $J_c$ and as close as possible to $J_c$; in other words, adjacent to $K_l \cup L \cup K_r$. Choose the lengths of $I_l$ and $I_r$ so that

$$\frac{|\pi(\widehat{I_l})|}{|T_1|} = \frac{|\pi(\widehat{I_r})|}{|T_1|} = \frac{3}{2} \cdot 3^{-Q}, \tag{14.17}$$

for some large integer $Q$. The harmonic measure of $\pi(\widehat{I_l}) \cup \pi(\widehat{I_r})$ in $\Lambda_{J_c}$ as seen from $z_J$ is a positive constant independent of $I$, since the domains $\Lambda_{J_c}$ for different $I$ are all Euclidean dilations of each other. By a comparison of harmonic measure, there is an $\eta > 0$ such that $|I_l|$ and $|I_r| \geq \eta |I|$ for all $I$. Now let

$$F(x) = \begin{cases} (1-\varepsilon)\dfrac{|T_1|}{|B|}, & x \in B; \\ 1, & x \in I_l \cup I_r; \\ \delta, & x \in T_1 \setminus (B \cup I_l \cup I_r). \end{cases} \tag{14.18}$$

Again, we choose $\delta$ so that $F$ has mean value one on $T_1$. Then

$$\delta = \left( \varepsilon - \frac{|I_l \cup I_r|}{|T_1|} \right) \frac{|T_1|}{|T_1 \setminus (B \cup I_l \cup I_r)|} \geq \varepsilon/2, \tag{14.19}$$

say, if $Q$ is chosen large enough that $|I_l \cup I_r|/|T_1| \leq \varepsilon/2$. Therefore $F$ is bounded above and below on $T_1$ by positive constants independent of $I$.

Notice that with this new definition of $F$ on $T_1$, we can still make the mean value of $X_i F$ on $T_1$ greater than $c_1$. For, just as for standard intervals, harmonic measure estimates show that $|B_k \cap (I_l \cup I_r)|/|T_1|$ decays exponentially in $k$, and so the contribution from $I_l \cup I_r$ to the mean value goes to zero as $Q$ goes to infinity.

This completes the proof of Lemma 14.3. $\qquad\blacksquare$

## 15. $\mu$ is supported on $L_c(G)$

Define the functions $F_n$ as in Section 11 on all standard intervals $I \in \mathcal{H}_{n-1}$, $n \geq 1$, and as in Section 14 on all non-standard intervals $I \in \mathcal{H}_{n-1}$, $n \geq 1$, such that $J_c = \pi(\widehat{I})$ is the non-standard central Whitney interval in any component of $[0,1] \setminus K$. Define $F_n \equiv 1$ on the remaining grid intervals $I \in \mathcal{H}_{n-1}$, $n \geq 1$; these are exactly those $I \in \mathcal{H}$ such that $\pi(\widehat{I}) = J_\infty$. We have shown that the first two types of these functions are $(\delta, \eta)$-suitable for all $I \in \mathcal{H}_{n-1}$, $n \geq 1$, where $\delta$ and $\eta$ are constants independent of $I$ and $n$. This is also true when $F_n \equiv 1$ on $I \in \mathcal{H}_{n-1}$. Therefore, by Lemma 6.3, the measures $\mu_n$ defined by $d\mu_n = F_n(x) \cdots F_1(x)\, dx$ converge to a doubling measure $\mu$ on the circle.

Consider the random walk on the tree discussed in Section 4, where the vertices in the tree are the Whitney intervals $\widehat{I}$ in the boundary arcs $\bigcup_n \mathcal{A}_n$ of the half fundamental domains in our tiling of the disc; $V$ is the subset of vertices such that $\pi(\widehat{I}) = J_\infty$ (in other words, those intervals which contain orbit points $g(0)$, $g \in G$); and the probabilities of jumps between adjacent vertices are determined by the functions $F_n$ defined above.

In this section, we first use the estimates $EX_i \geq c_1$ and $EX_i^2 \leq c_2$ established in Sections 13 and 14 to show that a particle starting from any vertex $v$ reaches $V(v)$ with probability one. Here $V(v)$ is the set of vertices $w \in V$ below $v$ such that there are no other vertices from $V$ between $w$ and $v$. We conclude that the doubling measure $\mu$ is supported on the set $S$ of points which lie in infinitely many grid intervals $I$ such that $\pi(\widehat{I}) = J_\infty$. Finally, we show that the conical limit set of $G$ contains this set $S$, which establishes Theorem 1.2.

Let $v$ be a vertex in the tree, and $I_0 \in \mathcal{H}_s$ the corresponding grid interval. Let $V(I_0)$ be the collection of maximal grid intervals $I \subset I_0$ such that $\pi(\widehat{I}) = J_\infty$. We show that the restriction of $\mu$ to $I_0$ is supported on $V(I_0)$.

Recall that

$$X_i(x) = \begin{cases} \log_3 \left[ \frac{|J_i(x)|}{|J_{i-1}(x)|} \right], & \text{if } J_1(x), \dots, J_{i-1}(x) \neq J_\infty; \\ 1, & \text{otherwise}, \end{cases} \tag{15.1}$$

that $S_k(x) = \sum_{i=1}^{k} X_i(x)$, and that $V(I_0) = \{x \in I_0 \mid S_k(x) \to +\infty\}$. It is enough to show that $S_k(x) \to +\infty$ almost surely on $I_0$ with respect to $\mu$.

Let $E_I X_i$ denote the mean value of $X_i$ on $I$ with respect to $\mu$. We know that $|E_I X_i| < \infty$ for all grid intervals $I \in \mathcal{H}$. This is because $|E_I X_i| \leq E_I |X_i| \leq E_I X_i^2 < c_2$, since $X_i$ is integer-valued. Define new functions $\widetilde{X}_i(x)$ on $I_0$ which have mean zero on each $I \in \mathcal{H}_{s+i-1}$ which lies in $I_0$:

$$
\begin{aligned}
\widetilde{X}_1(x) &= X_1(x) - E_{I_0} X_1; \\
\widetilde{X}_i(x) &= \sum_{\substack{I \in \mathcal{H}_{s+i-1} \\ I \subset I_0}} (X_i(x) - E_I X_i)\, \chi_I(x),
\end{aligned}
\tag{15.2}
$$

for $i \geq 1$, and let

$$
\begin{aligned}
\widetilde{S}_k(x) &= \sum_{i=1}^{k} \widetilde{X}_i(x) \\
&= \sum_{i=1}^{k} \left( X_i(x) - E_{I_{i-1}} X_i \right)
\end{aligned}
\tag{15.3}
$$

where $x \in I_{i-1}$, $i \geq 1$.

Then $\left\{ \widetilde{S}_k \right\}$ is a martingale with respect to $\mu$.

Also, $E\widetilde{X}_i^2$ is bounded by $EX_i^2$, which is bounded by $c_2$. Therefore, by the strong law of large numbers for martingales [Fel], $\widetilde{S}_k/k \to 0$ almost surely on $I_0$ with respect to $\mu$.

Now let $x$ be a point in $I_0$ such that $\widetilde{S}_k(x)/k \to 0$, and let $I_i \in \mathcal{H}_{s+i}$ be the interval containing $x$, for $i \geq 0$. Then

$$
\begin{aligned}
\frac{S_k(x)}{k} &= \frac{1}{k} \sum_{i=1}^{k} X_i(x) \\
&= \frac{1}{k} \sum_{i=1}^{k} \left( \widetilde{X}_i(x) + E_{I_{i-1}}^{d\mu} X_i \right) \\
&\geq \frac{1}{k} \widetilde{S}_k(x) + c_1 \\
&\geq c_1/2,
\end{aligned}
\tag{15.4}
$$

for all $k$ sufficiently large that $|\widetilde{S}_k(x)/k| \leq c_1/2$. Therefore $S_k(x) \geq c_1 k/2$ for all such $k$, and so $S_k(x) \to +\infty$. Thus $S_k \to +\infty$ almost surely with respect to $\mu$ on $I_0$. Since $V(I_0)$ is exactly the set of points $x$ in $I_0$ such that $S_k \to +\infty$, this implies that $\mu(V(I_0)) = \mu(I_0)$.

Since a particle starting from any vertex $v$ reaches $V(v)$ with probability one, it is now clear that a particle starting anywhere in the tree will pass through infinitely many vertices in $V$, with probability one. (Note that we have now defined the only jump probabilities which were missing; these are the probabilities of the jumps from vertices in the subset $V$, and they are determined by the functions $F_n \equiv 1$ on the corresponding grid intervals.) Thus we have shown that the doubling measure $\mu$ is supported on the set $S$ of points in the circle which lie in infinitely many grid intervals $I = P(\widehat{I})$ such that $\pi(\widehat{I}) = J_\infty$.

It remains to show that $S$ is contained in the conical limit set of $G$. We show that there is a uniform constant $l$ such that for each orbit point $g_j(0) \neq 0$, the Whitney interval $\widehat{I}_j$ which contains $g_j(0)$ projects to a grid interval $P(\widehat{I}_j)$ which lies in the spherical cap $\mathrm{Cap}\,(g_j(0), l)$ on $g_j(0)$. Then each point $x \in S$ lies in infinitely many of these caps, and therefore in the conical limit set.

By definition, the tip of the leaf based at $J_\infty$ is $E_{J_\infty} = [1/9, 8/9]$. For each $g_j \in G \backslash \{\mathrm{id}\}$ let $A_j$ be the orthocircular arc in $\bigcup_n \mathcal{A}_n$ which contains $g_j(0)$, and let $\Omega_j$ be the half fundamental domain below $A_j$. Let $\mathcal{F}_j$ be the fundamental domain consisting of $\Omega_j$ and the half fundamental domain above $A_j$. Let $\widehat{E}_j$ be the segment of $\partial \Omega_j$, below $g_j(0)$, such that $\pi(\widehat{E}_j) = E_{J_\infty}$. The endpoints of $\widehat{E}_j$ lie in $\partial \mathbb{D}$. Then $P(\widehat{I}_j) = E_j$ is the arc of $\partial \mathbb{D}$ below $\widehat{E}_j$ (with the same endpoints as $\widehat{E}_j$). See Figure 14.

Figure 14. $P(\widehat{I}_j) = E_j$ lies in $\mathrm{Cap} = \mathrm{Cap}\,(g_j(0), l)$; $E = E_{J_\infty}$.

**Lemma 15.1.** *There is a constant $l > 0$ such that for all $j$,*

$$P(\widehat{I}_j) \subset \mathrm{Cap}\,(g_j(0), l). \tag{15.5}$$

**Proof.** By definition,

$$\mathrm{Cap}\,(g_j(0), l) = \{z \in \partial\mathbb{D} \mid \mathrm{dist}\,(z, g_j(0)) \le l\,(1 - |g_j(0)|)\}.$$

Let $z_j = g_j(0)$. Let $b_j \in \partial\mathbb{D}$ be the endpoint of $P(\widehat{I_j})$ which is furthest from $z_j$. We show that $|z_j - b_j|/(1 - |z_j|) \le c$.

If $b_j$ satisfies $|b_j - z_j| \le 10\,(1 - |z_j|)$, we are done. If not, let $B_j$ be the component of $\pi^{-1}([0, 1] \setminus E_{J_\infty})$ which has $b_j$ as an endpoint; and let $B = \pi(B_j)$. See Figure 14 for one possible configuration.

First, $|z_j - b_j| \le c\,\mathrm{dist}\,(z_j, B_j)$. For if $a$ is any point in $B_j \cap \partial\mathbb{D}$, then $|z_j - b_j| \le |z_j - a|$, while if $a$ is in $B_j \setminus \partial\mathbb{D}$ then $a$ lies in one of the orthocircular arcs in $\partial\mathcal{F}_j$, and $|z_j - a| \ge c\,|z_j - b_j|$.

By Beurling's Lemma (see also (10.12)),

$$
\begin{aligned}
\mathrm{dist}\,(z_j, B_j) &\le c\,\frac{\mathrm{dist}\,(z_j, \partial\mathcal{F}_j)}{\omega(z_j, B_j, \mathcal{F}_j)^2} \\
&\le c\,\frac{\mathrm{dist}\,(z_j, \partial\mathcal{F}_j)}{\omega(\pi(z_j), \pi(B_j), \pi(\mathcal{F}_j))^2} \\
&= c\,\frac{\mathrm{dist}\,(z_j, \partial\mathcal{F}_j)}{\omega(\infty, B, \overline{\mathbf{C}} \setminus K)^2} \\
&\le c\,\mathrm{dist}\,(z_j, \partial\mathcal{F}_j) \\
&\le c\,(1 - |z_j|).
\end{aligned}
\tag{15.6}
$$

$\square$

We have shown that each point in $\partial\mathbb{D}$ which is at the end of a path which passes through infinitely many orbit points is actually a conical limit point of $G$; in particular, that the orbit points on the path not only accumulate at the endpoint but lie in some non-tangential cone based at the endpoint. Therefore the set $S$ lies in the conical limit set of $G$.

This completes the proof of Theorem 1.2.

Peter W. Jones, Department of Mathematics, Yale University, New Haven, CT 06520-8283 *E-mail address*: jones@math.yale.edu.

Lesley Ward, Department of Mathematics, Rice University, Houston, TX 77251-1892 *E-mail address*: lesley@math.rice.edu